\documentclass[a4paper,12pt]{article} 
\newdimen\paperhight
\newcommand{\st}{{1\over 16}} 
\newcommand{\hf}{\frac{1}{2}}
\newcommand{\ots}{\otimes}
\newcommand{\boplus}{\bigoplus}
\newcommand{\ts}{\times}

\newcommand{\ol}{\overline} 
\newcommand{\qed}{\begin{flushright}{\bf Q.E.D.}\end{flushright}}
\newcommand{\pr}{\par \vspace{3mm} \noindent {\bf [Proof]} \qquad}
\newcommand{\prend}{\hfill \qed}

\newcommand{\1}{{\bf 1}} 
\newcommand{\bw}{{\bf w}} 
\newcommand{\al}{\alpha} 
\newcommand{\be}{\beta} 
\newcommand{\de}{\delta}

\newcommand{\ga}{\gamma} 
\newcommand{\la}{\lambda} 
\newcommand{\La}{\Lambda}

\newcommand{\CG}{{\cal G}} 
 
\newcommand{\C}{{\mathbb C}} 
\newcommand{\Z}{{\mathbb Z}} 
\newcommand{\N}{{\mathbb N}} 
\newcommand{\Q}{{\mathbb Q}} 
 
\newcommand{\R}{{\mathbb R}} 
\newcommand{\M}{{\mathbb M}} 
\newcommand{\Ind}{{\rm Ind}}

\newcommand{\Aut}{{\rm Aut}}

\newcommand{\End}{{\rm End}}
 
\newcommand{\tr}{{\rm tr}\ }
 
\newcommand{\ch}{{\rm ch}\ } 

\newtheorem{thm}{Theorem}[section]
\newtheorem{prn}{Proposition}[section]

\newtheorem{lmm}{Lemma}[section]
\newtheorem{cry}{Corollary}[section]
\newtheorem{cnj}{Conjecture}

\newtheorem{rmk}{Remark}
\usepackage{amssymb}
\usepackage[usenames]{color}
\usepackage{graphicx}

\begin{document}
\title{A new construction of the moonshine vertex operator algebra 
over the real number field. } 
\author{Masahiko Miyamoto} 
\date{\begin{tabular}{c}
Institute of Mathematics \cr
University of Tsukuba \cr
Tsukuba 305, Japan \cr
\end{tabular} }
\maketitle
\begin{abstract}
We give a new construction of the moonshine module vertex 
operator algebra $V^{\natural}$ over the real field, 
which was originally constructed in \cite{FLM2}. 
The advantage of our construction is that 
we can easily prove the facts 
that $V^{\natural}$ has a positive 
definite invariant bilinear form and 
${\Aut}(V^{\natural})$ is the Monster simple group.   
In addition, we construct a lot of conformal vectors in $V^{\natural}$ 
which give rise to 2A-involutions. 
We also construct an infinite series of holomorphic VOAs. 
Each of them has exactly one irreducible module and its 
full automorphism group is finite. At the end of the paper, 
we will calculate the character 
of a 3C element of the Monster simple group. 
\end{abstract}

\renewcommand{\baselinestretch}{1.2}\large\normalsize
\section{Introduction}
All VOAs in this paper are defined over the real number field ${\R}$ 
and ${\C}V$ denotes the VOA ${\C}\otimes_{\R}V$ for a VOA $V$. \par
The most interesting example of 
vertex operator algebra (VOA) is the 
moonshine module VOA 
$V^{\natural}=\sum_{i=0}^{\infty}V_i^{\natural}$. 
Although it has many interesting properties, 
the original construction \cite{FLM2} essentially depends on 
the actions of the centralizer $2^{1+24}Co.1$ 
of an involution called $2B$ of the 
Monster simple group and so 
it is hard to see the actions 
of the other elements explicitly. 
The Monster simple group has the other conjugacy class of involutions 
called 2A. We will construct the Moonshine VOA $V^{\natural}$ 
from the point of view of elementary abelian 2-group generated by 
2A-elements. \par
The simplest example of VOA is the rational 
Virasoro VOA $L(\hf,0)$ of 
the minimal series with central charge $\hf$. 
It has only three irreducible modules 
$L(\hf,0)$, $L(\hf,\hf)$, and $L(\hf,\st)$, 
where the first entry is the central charge and 
the second entry denotes the lowest weights. 
Its fusion rules (or fusion products) are known as, see \cite{BPZ} 
or \cite{DMZ}: \\
$$
\begin{array}{l}
(1) \qquad L(\hf,0) \mbox{ is identity}, \\
(2) \qquad L(\hf,\hf)\ts L(\hf,\hf)=L(\hf,0), \\
(3) \qquad L(\hf,\hf)\ts L(\hf,\st)=L(\hf,\st), \\
(4) \qquad L(\hf,\st)\ts L(\hf,\st)
=L(\hf,0)+L(\hf,\hf).
\end{array} \eqno{(1.1)}
$$
For a VOA $V$, we call $e\in V_2$ a rational conformal vector if a sub VOA 
$<e>$ generated by $e$ is a rational VOA 
and $e$ is the Virasoro element of $<e>$.  We are essentially interested 
in a rational conformal vector $e$ with central charge $\hf$. Under 
this assumption, $<e>$ is isomorphic to 
$L(\hf,0)$ and we can view $V$ as a $<e>$-module. 
The fusion rules of $L(\hf,0)$ will then play an important role 
in our arguments. In particular, we will use an automorphism 
$\tau_{e}$ of $V$, which is defined by the author in \cite{Mi1},  
for each rational conformal vector $e$ with 
central charge $\hf$, where $\tau_e$ is given by 
$$ \tau_{e}: \left\{ \begin{array}{rcl}
1& on&\mbox{ all $<e>$-submodules isomorphic to $L(\hf,0)$ or $L(\hf,\hf)$ }\cr
-1& on&\mbox{all $<e>$-submodules isomorphic to $L(\hf,\st)$}.
\end{array} \right. $$
Since we will treat only rational VOAs $V$, the 
tensor product of two $V$-modules $W^1$ and $W^2$ is 
well-defined \cite{Li} 
and it is equal to the fusion product $W^1\times W^2$. Therefore, 
we will also consider a fusion product as a module. From now on, 
$\otimes_{i=1}^n W^i$ means a $\otimes_{i=1}^n V^i$-module for 
$V^i$-modules $W^i$. \par
In this paper, we will consider a set of mutually orthogonal 
rational conformal vectors $\{e^i:i=1,...,n\}$ with central charge $\hf$ 
such that the sum $\sum e^i$ is Virasoro element of $V$.  
Here, "orthogonal" means $e^i_1e^j=0$ for $i\not=j$. 
We will call such a set of conformal vectors "a coordinate set". 
Thus, the sub VOA $T=<e^1,...,e^n>$ is isomorphic to $L(\hf,0)^{\ots n}$ and  
it is known that every irreducible $T$-module $W$ 
is a tensor product $\ots_{i=1}^n L(\hf,h^i)$ of irreducible 
$L(\hf,0)$-modules $L(\hf,h^i)$, see \cite{DMZ}.  
Define a binary word 
$$\tilde{\tau}_T(W)=(a_1,...,a_n)\eqno{(1.2)} $$ 
by $a_1=1$ if $h_i=\st$ and $a_i=0$ if $h_i=0$ or $\hf$.  
We call it a (binary) $\tau$-word 
since it is corresponding to the actions of automorphisms $\tau_{e^i}$.  
(The author 
once called it a word of $\st$-positions and denoted it by $\tilde{h}(W)$ 
in \cite{Mi3} and  \cite{Mi4}.)  
We note that $T$ is rational and 
the fusion product is given by 
$$ (\otimes_{i=1}^n W^i)\times (\otimes_{i=1}^n U^i)
=\otimes_{i=1}^n (W^i\times U^i)  $$
as Dong, Mason and Zhu proved in \cite{DMZ}. \par
We will construct 
the moonshine VOA $V^{\natural}$ over the real number field ${\R}$
as a direct sum of irreducible $T$-modules.  
It is not difficult to construct the underlining space $V^{\natural}$ 
as a direct sum of irreducible $T$-modules.  
Originally, it has shown by Dong, Mason and Zhu \cite{DMZ} that 
the moonshine VOA $V^{\natural}$ of rank 24 
contains 48 mutually orthogonal conformal vectors $e^i$ with the central 
charge $\hf$ such that the sum is the Virasoro element of 
$V^{\natural}$ and 
the author determined the multiplicities 
of all irreducible $T$-submodules of $V^{\natural}$ for some $T$ in 
\cite{Mi4}. \par
The reason why we will treat a VOA over ${\R}$ is that 
a positive definite 
invariant bilinear form on a VOA is very useful to determine 
an automorphism group.  For example, 
Frenkel, Lepowsky and Meurman constructed the moonshine VOA 
$V^{\natural}$ 
over the real number field 
(in fact, they constructed it over the rational number field) 
and have shown that 
$V^{\natural}$ has a positive definite invariant bilinear 
form in \cite{FLM2}.  
One of our important tools in determining the full 
automorphism group is the uniqueness theorem for a VOA 
satisfying Hypotheses I mentioned later.  
This holds for VOAs over the complex number 
field without assuming the positive definite invariant bilinear form 
(see \cite{Mi5}).  
However, it is not uniquely determined for VOAs over the real number 
field. In order to avoid this anomaly, we will treat only VOAs over ${\R}$ 
with positive definite invariant 
bilinear forms. For example, 
the code VOAs which the author defined in \cite{Mi2} 
have a positive definite 
invariant bilinear form if we construct them over ${\R}$. 
This setting offers us exactly the same situation as in VOAs over 
the complex field.  

One of our tools is a code VOA $M_D$ and its representation theory for 
an even linear binary code $D$ of length $n$ (\cite{Mi2} and \cite{Mi3}). 
We will briefly 
explain it in \S 3.  The characterization of 
a code VOA $M_D$ is    
that it is a simple VOA $V$ containing $T=L(\hf,0)^{\ots n}$ 
and $\tilde{\tau}_T(V)=(0^n)$. 
Any irreducible $M_D$-module $W$ is a direct sum of 
irreducible $T$-modules $U^i$ since $T$ is rational.  
By the fusion rules of Ising models (1.1), 
$\tilde{\tau}(U^i)$ is uniquely determined and so 
we use the same notation $\tilde{\tau}(W)$ for it.   
If a simple VOA $V$ contains a coordinate set $\{e^1,...,e^n\}$, 
then $P=<\tau_{e^i}:i=1,...,n>$ is 
an elementary abelian automorphism group. 
Decompose $V$ into a direct sum
$$ V=\oplus_{\chi\in Irr(P)}V^{\chi}    $$
of eigenspaces of $P$, 
where $V^{\chi}=\{v\in V: gv=\chi(g)v \mbox{ for } g\in P\}$ and 
$V^{1}=V^P$ is the set of $P$-invariants.  
It is known by \cite{DM2} that $V^{\chi}$ is 
a nonzero irreducible $V^P$-module. 
It follows from the definition of $\tau_{e^i}$ 
that $V^P$ contains $T=<e^1,...,e^n>$ and 
is isomorphic to a code VOA. Moreover, 
if $\tilde{\tau}_T(V^{\chi})=(a_i)$ 
then $\chi(\tau_{e^i})=(-1)^{a_i}$. Therefore, 
the representation theory of code VOA plays an essential role in 
the study of such VOAs.

Another tool is "induced VOA". 
In \cite{Mi3}, we introduced a concept of the 
induced ${\C}M_D$-module $\Ind_{E}^{D}({\C}W)$ for 
a subcode $E$ containing a maximal self-orthogonal subcode of 
$D_{\be}=\{\al\in D|Supp(\al)\subseteq Supp(\be)\}$ and  
an ${\C}M_E$-module ${\C}W$ satisfying 
$\langle \tilde{\tau}(W),D\rangle=0$.  
This is a special case of the concept of 
induced modules defined in \cite{DLi}.  
We will also define an induced module for code VOAs over ${\R}$ and  
apply it to a VOA here. Namely, 
${\Ind}_E^D(W)$ becomes a VOA if $W$ is a VOA under some conditions. 
An advantage of this construction is that 
it keeps a positive definite invariant 
bilinear form. 
As applications, we will 
construct several holomorphic VOAs from a VOA. For example, 
we will construct a lattice VOA $V_{\Lambda}$ of the Leech lattice 
from $V^{\natural}$ by restricting and defining an induced VOA.

We should note that it is possible to construct $V^{\natural}$ 
over the rational number field by our way. However, it makes us 
add several conditions to get the uniqueness theory and we will avoid 
such complications.

We will prove that our VOA $V^{\natural}$ over ${\R}$ is a 
holomorphic VOA of rank 24. 
We will also prove that 
the full automorphism group is the Monster 
simple group ${\M}$.  These information are enough to show that 
our VOA $V^{\natural}$ is isomorphic to the moonshine module VOA constructed 
in \cite{FLM2}.    

We are now in position to mention the outline of this paper. 
In this paper, we will not only construct the moonshine VOA but 
also VOAs with similar structures. 
Our essential tool is the following theorem,  
which was proved for VOAs over ${\C}$ 
by the author in \cite{Mi5}. We will show that 
this theorem is also true for VOAs over ${\R}$ with  
positive definite invariant bilinear forms.  \\

\noindent
{\bf Hypotheses I} \\
(1) $D$ and $S$ are both even linear codes of length $8k$ 
and $S\subseteq D\cap D^{\perp}$. \\
(2)  For any $\al, \be\in S, (\al\not=\be)$, there is  
a self-dual 
subcode $E=E_{\al}\oplus E_{\al^c}$ of $D$ and maximal self-orthogonal 
(doubly even) subcodes $H^{\be}$ 
and $H^{\al+\be}$ of $D_{\be}$ and $D_{\al+\be}$ 
containing $E_{\be}$ and $E_{\al+\be}$, 
respectively, such that  \\
(2.1) \quad $E_{\al}$ and 
$E_{\al^c}$ are direct sums of $[8,4,4]$-Hamming codes, \\
(2.2) \quad $H^{\be}+E=H^{\al+\be}+E,$  \\
where $\al^c$ denotes the complement $(1^{8k})-\al$ and 
$S_{\delta}$ denotes a subcode 
$\{\ga\in D:Supp(\ga)\subseteq Supp(\delta)\}$ of a code $S$ and \\
(3) There is an $S$-graded $M_D$-module $V=\boplus_{\al\in S} V^{\al}$ such 
that each $V^{\al}$ is an $M_D$-submodule with $\tilde{\tau}(V^{\al})=\al$. 
In particular, $V^{(0^{8k})}\cong M_D$ as $M_D$-modules. \\
(4)  For $ \al,\be \in S-\{(0^n)\} $ and $\al\not=\be$, 
$$V^{\al,\be}=M_D\oplus V^{\al}\oplus V^{\be}\oplus V^{\al+\be}  $$
has a simple VOA structure containing $M_D$ as a sub VOA. \\
(5)  $V^{\al,\be}$ has a positive definite invariant 
bilinear form. \\

\noindent
{\bf Theorem 3.3} \quad 
{\it 
Under the above assumptions (1)$\sim $(5) of Hypotheses I, we obtain  
the fusion product $V^{\al}\times V^{\be}=V^{\al+\be}$ for $\al,\be\in D$ and 
$$  V=\boplus_{\al\in S}V^{\al}   $$
has a structure of simple VOA with $M_D$ as a sub VOA and 
it has a positive definite invariant bilinear form. 
The structure of VOA on $V$ with a positive definite 
invariant bilinear form is uniquely 
determined up to $M_D$-isomorphisms. }\\

We note that an important assertion of Theorem 3.3 is that the fusion 
product $V^{\al}\times V^{\be}$ is irreducible, that is, if 
$I(\ast,z)\in I\pmatrix{V^{\al+\be}\cr V^{\al},\quad 
V^{\be}}$ is a nonzero intertwining operator, then 
for any VOA structure $(V,\tilde{Y})$ there is a scalar $\la$ such that 
$\tilde{Y}(v,z)_{|V^{\be}}=\la I(v,z)$.
As we will show that the uniqueness of 
the VOA structure on $V$ comes from this property. 
 
The assumptions (1) and (2) are conditions on the codes $D$ and $S$. 
So our construction is just 
to collect a set $\{V^{\al}:\al\in S\}$ of $M_D$-modules 
satisfying (4) and (5). 
In order to prove the condition (4), we will use the following theorems. 
These are also essentially based on the 
irreducibility of fusion products. \\

\noindent
{\bf Theorem 3.2} \quad {\it Assume that (1) and (2) of Hypotheses I 
hold for $(S,D)$. Choose $\al,\be\in S$ so that $\dim <\al,\be>=2$, 
where $<\al,\be>$ is the 
code generated by $\al$ and $\be$.  
Let $F$ be an even linear code containing $D$ and 
assume $\al,\be\in F^{\perp}$. 
If $U=M_D\oplus W^{\al}\oplus W^{\be}\oplus W^{\al+\be}$ 
has a simple VOA structure satisfying $\tilde{\tau}(W^{\ga})=\ga$ for 
$\ga\in <\al,\be>$, then 
$${\rm Ind}_D^F(U)=M_F\oplus {\Ind}_{M_D}^{M_F}(W^{\al})
\oplus {\Ind}_{M_D}^{M_F}(W^{\be})\oplus 
{\Ind}_{M_D}^{M_F}(W^{\al+\be})$$ 
has a simple VOA structure. }  \\

{\bf Theorem 4.1} \quad {\it  Under the assumptions in Theorem 3.2, 
if a VOA $U$ has a positive definite invariant 
bilinear form, then ${\rm Ind}_D^F(U)$ has a VOA structure with a 
positive definite invariant bilinear form. Furthermore, 
such a VOA is uniquely determined up to $M_F$-isomorphisms.}\\

In order to construct a VOA by using Theorem 3.3, 
it is sufficient to 
collect $M_D$-modules satisfying (4) and (5) for a small code $D$ 
as we showed in Theorem 4.1.  
We will gather such modules from the lattice VOA $\tilde{V}_{E_8}$ 
with a positive definite invariant bilinear form  
constructed from a even unimodular lattice of type $E_8$.  
We will also prove that $\tilde{V}_{E_8}$ has a structure satisfying 
Hypotheses I. Namely, $\tilde{V}_{E_8}$ contains 
16 mutually orthogonal conformal vectors $\{e^1,...,e^{16}\}$ such that \\
(1) \quad the order of $P=<\tau_{e^i}:i=1,...,8>$ is $32$, \\
(2) \quad $(\tilde{V}_{E_8})^P$ is 
isomorphic to a code VOA $M_{D_{E_8}}$, where $D_{E_8}$ is a 
Reed Muller code $RM(4,2)$ and  \\
(3) \quad $\tilde{V}_{E_8}=\boplus_{\al\in S_{E_8}}\tilde{V}_{E_8}^{\al}$, 
where 
$S_{E_8}=D^{\perp}\cong RM(4,1)$, $\tilde{V}_{E_8}^{(0^8)}\cong M_{RM(4,2)}$ 
and $\tilde{V}_{E_8}^{\al}$ are irreducible $M_{RM(4,2)}$-modules.  
Note that 
$$S_{E_8}=<(1^{16}), \ (0^81^8), \ (\{0^41^4\}^2), \ (\{0^21^2\}^4), 
\ (\{01\}^8) > \eqno{(1.3)}$$ 
and the weight enumerator of $S_{E_8}$ is $x^{16}+30x^8y^8+y^{16}$. 
Moreover, the 
minimal weight of $D_{E_8}$ is 4 and the pair $(D_{E_8}, S_{E_8})$ 
satisfies the conditions (1) and (2) of Hypotheses I, see Lemma 5.1. 
Therefore, a VOA structure on the $M_{D_{E_8}}$-module 
$\tilde{V}_{E_8}=\oplus_{\al\in S_{E_8}}\tilde{V}_{E_8}^{\al}$ is uniquely 
determined by Theorem 3.3.  We also have 
a fusion product $\tilde{V}_{E_8}^{\al}\times \tilde{V}_{E_8}^{\be}
=\tilde{V}_{E_8}^{\al+\be}$ 
of $M_{D_{E_8}}$-modules for any $\al,\be\in S_{E_8}$.  \par

We will next explain how to construct the moonshine VOA. 
In order to define the moonshine VOA $V^{\natural}$, we will set 
$$S^{\natural}=\{(\al,\al,\al),(\al,\al,\al^c),(\al,\al^c,\al),
(\al^c,\al,\al): \al\in S_{E_8}\}, \eqno{(1.4)} $$ 
where $\al^c$ is the complement of $\al$.  
Set $D^{\natural}=(S^{\natural})^{\perp}$ and call it a moonshine 
code. It is of dimension $41$ and contains 
$D_{E_8}^3=D_{E_8}\oplus D_{E_8}\oplus D_{E_8}$.  
We note that $S^{\natural}$ and 
$D^{\natural}$ 
are even linear codes of length 48. 
Clearly, the pair $(D_{E_8}^3, S^{\natural})$ satisfies the conditions (1) and 
(2) of Hypotheses I. \par

Our construction consists of the following three steps. \\
First, since $\tilde{V}_{E_8}^{\al}\times \tilde{V}_{E_8}^{\be}
=\tilde{V}_{E_8}^{\al+\be}$ for $\al,\be\in S_{E_8}$, 
$$V^1=\boplus_{(\al,\be,\ga)\in S^{\natural}}
(\tilde{V}_{E_8}^{\al}\ots \tilde{V}_{E_8}^{\be}\ots \tilde{V}_{E_8}^{\ga}) 
\eqno{(1.5)}$$
is a sub VOA of $\tilde{V}_{E_8}\ots \tilde{V}_{E_8}\ots \tilde{V}_{E_8}$. 
Clearly, $V^1$ has a positive definite invariant bilinear form. 
Our second step is to twist it. Namely, set $\xi_1=(10^{15})$ and 
let $R=M_{D_{E_8}+\xi_1}$ be a coset module. To simplify the 
notation, denote $R\times \tilde{V}_{E_8}^{\al}$ by $R\tilde{V}_{E_8}^{\al}$. 
Set
$$Q=<(\xi_1\xi_1 0^{16}),(0^{16}\xi_1\xi_1)>
\subseteq {{\Z}}_2^{48}.  $$
We induce $V^1$ to
$$ V^2={\Ind}_{D_{E_8}^3}^{D_{E_8}^3+Q}(V^1).  $$ 
Although $V^2$ is not a VOA, we can find the following $M_{D^3}$-submodules 
in $V^2$: 
$$ \begin{array}{l}
W^{(\al,\al,\al)}=\tilde{V}_{E_8}^{\al}\ots \tilde{V}_{E_8}^{\al}
\ots \tilde{V}_{E_8}^{\al} \cr
W^{(\al,\al,\al^c)}=(R\tilde{V}_{E_8}^{\al}) \ots (R\tilde{V}_{E_8}^{\al})
\ots \tilde{V}_{E_8}^{\al^c} \cr
W^{(\al,\al^c,\al)}=(R\tilde{V}_{E_8}^{\al}) \ots \tilde{V}_{E_8}^{\al}
\ots (R\tilde{V}_{E_8}^{\al}) \cr 
W^{(\al^c,\al,\al)}=(\tilde{V}_{E_8}^{\al}) \ots (R\tilde{V}_{E_8}^{\al})
\ots (R\tilde{V}_{E_8}^{\al}). 
\end{array} $$
for $\al\in S_{E_8}$.  
At the end, we set
$$(V^{\natural})^{\chi}
={\rm Ind}_{D^3}^{D^{\natural}}(W^{\chi}) $$
for $\chi\in S^{\natural}$. 
We will show that these $M_{D^{\natural}}$-modules $(V^{\natural})^{\chi}$ 
satisfy the condition (4) of Hypotheses I. Therefore, 
we obtain a VOA 
$$  V^{\natural}=\boplus_{\chi\in S^{\natural}}(V^{\natural})^{\chi}  $$
which possesses a positive definite 
invariant bilinear form. 
Since we construct $V^{\natural}$ under the 
condition $S^{\natural}=(D^{\natural})^{\perp}$, 
$V^{\natural}$ is the only irreducible $V^{\natural}$-module by Theorem 6.1. 
From the construction, we will see that 
$\dim (V^{\natural})_0=1$ and $(V^{\natural})_1=0$.  
It comes from the structure of $V^{\natural}$ and 
the multiplicity of irreducible $M_{D^{\natural}}$-submodules that 
$q^{-1}ch_{V^{\natural}}=J(q)=q^{-1}+196884q+...$ is the J-function. 
We will also see that the full automorphism group of $V^{\natural}$ 
is the Monster simple group.  Although it is not easy 
to determine the full automorphism groups of VOAs in general, 
our construction 
has certain advantages.  For example, it is easy to prove that 
the full automorphism group of a VOA satisfying Hypotheses I 
is finite if $V_1=0$ (Theorem 9.2). Furthermore, 
if $S$ is a subcode of $\{(\al,\al): \al\in {\Z}_2^{n/2}\}$ 
by rearranging the order of coordinates, then we will show that 
our VOA is a sub VOA of some lattice VOA with rank $n$ by 
the uniqueness of VOA structures.  
Also since our VOA $V$ contains a lot of rational 
conformal vectors $\{e^i:i\in I\}$ 
with central charge $\hf$, $V$ has a large automorphism group 
generated by $\{\tau_{e^i}:i\in I\}$, which is clearly a normal 
subgroup of ${\Aut}(V)$.  
Using these properties, 
we will prove that 
the space $(V^{\natural})^{<\de>}$ of $\de$-invariant is 
isomorphic to $V_{\La}^{\theta}$ for a lattice VOA $V_{\La}$ 
of the Leech lattice and an automorphism $\theta$ of $V_{\La}$ 
induced from $-1$ on $\La$ for $\de=\tau_{e^1}\tau_{e^2}$. 
For a conformal vector 
$e\in (V^{\natural})^{<\de>}\cong V_{\La}^{\theta}$, 
we can define automorphisms 
$\tau_e\in {\Aut}(V^{\natural})$ 
and $\tilde{\tau}_e\in {\Aut}(V_L)$. By this correspondence, 
we can calculate $C_{{\Aut}(V^{\natural})}(\delta)$. 
Also, we can calculate 
$N_{{\Aut}(V^{\natural})}(<\tau_{e^1}\tau_{e^2}, \tau_{e^1}\tau_{e^3}>)$ 
and 
$N_{{\Aut}(V^{\natural})}(<\tau_{e^1}\tau_{e^2}, \tau_{e^1}\tau_{e^3}, 
\tau_{e^1}\tau_{e^5}>)$.  
By this information, we can conclude that ${\Aut}(V^{\natural})$ 
is the Monster simple group and 
$V^{\natural}$ coincides with the moonshine module VOA constructed 
in \cite{FLM2}. 
Thus, this is a new construction of the moonshine VOA and the monster 
simple group.  \\
 
\noindent
{\bf Remark} \qquad 
It is possible to induce $V^1$ in (1.5) into a VOA  
$$ \tilde{V}={\Ind}_{D^3}^{D^{\natural}}(V^1) $$ 
directly. It follows from a direct calculation 
and the fusion rule (1.1) that 
$\tilde{V}_1$ is a commutative Lie algebra of dimension 24.  Since 
$\tilde{V}$ is a holomorphic VOA by Theorem 6.1,  
$\tilde{V}$ is isomorphic to the lattice VOA $V_{\Lambda}$ of 
Leech lattice $\Lambda$ by \cite{Mo}, (see Section 9). \\

Another important theorem in this paper is that if $S=D^{\perp}$ then 
a simple VOA $V$ satisfying Hypotheses I 
has the exactly one irreducible $V$-module $V$, see Theorem 6.1. 
Since Dong, Griess and H\"ohn \cite{DGH} 
have proved that a simple VOA satisfying 
Hypotheses I is rational, 
the VOAs $V=\oplus_{m=0}^{\infty} V_m$ satisfying $S=D^{\perp}$ 
are holomorphic and so 
$q^{-n/48}\sum (\dim V_m)q^m$ is a modular function of $SL_2({\Z})$ 
with a linear character by \cite{Z}.  \\

In \S 4, we construct a VOA $V_{E_8}$ with a positive definite 
invariant bilinear form. In \S 5, 
we investigate the structure of $V_{E_8}$. In \S 7, 
we construct the moonshine VOA 
$V^{\natural}$. In \S 8, we will construct a lot of rational 
conformal vectors of $V^{\natural}$ explicitly. 
In \S 9, we prove that ${\Aut}(V^{\natural})$ is the Monster simple group 
and $V^{\natural}$ is equal to the one constructed in \cite{FLM2}.  In \S 10, 
we will construct an infinite series of holomorphic VOAs with finite 
full automorphism groups.  
In \S 11, we will calculate the characters of some elements of the 
Monster simple group. \\
\noindent
{\bf  Acknowledgment} \\
The author wishes to thank K. Harada, T.~Kondo and H.~Yamaki 
for their helpful advises.
The author also would like to express the appreciation to J.~Lepowsky  
for his useful comments.

\section{Notation and preliminary results}
We adopt all notation and results from \cite{Mi3} and recall 
the construction of a lattice VOA. 
\subsection{Notation}
\begin{tabular}{ll}
$\al^c$ &The complement $(1^n)-\al$ of a binary word $\al$. \cr
$D, \ D(m)$  &Even binary linear codes, also see \S 4. \cr
$D_{\be}$ &$=\{\al\in D: Supp(\al)\subseteq Supp(\be)\}$. \cr 
$D^3$ &$=\{(\al,\be,\ga):\al,\be,\ga\in D\}$. \cr
$D^{\natural}, S^{\natural}$ &The moonshine codes. See (1.4). \cr
$D_{E_8}, S_{E_8}$ &See (1.3). \cr
$\{ e^i\ |\ i=1,...,n\}$ &A set of mutually 
orthogonal rational conformal vectors \cr
&with central charge $\hf$. \cr
$e^{\pm}(x)$ &
$={\st}x(-1)^1{\bf 1}\pm {1\over 4}(\iota(x)+\iota(-x))\in V_L$: \cr 
&the conformal vectors defined by $x\in L$ with $\langle x,x\rangle=4$. \cr
$E_8$, $E_8(m)$  &An even unimodular lattice of type $E_8$, also see (5.1).\cr
$\{f^i:i\},\{d^i:i\}$ &The other sets of mutually orthogonal 
eight conformal \cr
&vectors in a Hamming code VOA $M_{H_8}$, see \cite{Mi5}. \cr
$H_8$ &The $[8,4,4]$-Hamming code. \cr
$H(\hf,\al)$, $H(\st,\be)$ &The irreducible $V_{H8}$-modules, see Def.13 
in \cite{Mi5}. \cr
$\Ind_{M_E}^{M_D}(U)$ &The induced $M_D$-module from an $M_E$-module $U$, \cr
&see Sec.5.2 in Sec.6.2 in \cite{Mi5}. \cr
 $\iota(x)$ &  A vector in a lattice VOA 
 $V_L=\boplus_{x\in L}M(1)\iota(x)$, see \cite{FLM2}. \cr
 $L$ &A lattice. \cr
 $M_{\be+D}$ &A coset module 
 $\boplus_{(a^i)\in \be+D}\left( 
 (\ots_{i=1}^n M_{a^i})\ots e^{(a^i)}\right)$. \cr
 $M_D$ &A code VOA, see \S 3. \cr
 $Q$   &$=<(10^{15}10^{15}0^{16}), (10^{15}0^{16}10^{15}) >$. \cr
 $R$ &$M_{(10^7)+D}$. \cr
 $RV_{E_8}^{\al}$ &$R\times V_{E_8}^{\al}$. \cr
$\tilde{\tau}(W)$ &A $\tau$-word $(a_1,...,a_n)$, see (1.2). \cr
$T$  &$=\ots_{i=1}^nL(\hf,0)$. \cr
$\times$ &A fusion rule or a tensor product. \cr
$A(x,z)\sim B(x,z)$ &$(x-z)^n(A(x,z)-B(x,z))=0$ for an $n\in \N$. \cr
$\theta$ &An automorphism of $V_L$ defined by $-1$ on $L$. \cr
 $V_L$ &A lattice VOA $\boplus_{x\in L}M(1)\iota(x)$, 
 see \cite{FLM2} and \S 2.2. \cr
$\xi_i$&A word which is $1$ in the $i$-th entry and $0$ 
everywhere else, \cr
 &for example, $(0^{i-1}10^{n-i})$, 
$(0^{i-1}10^{8-i})$, $(0^{i-1}1^{16-i})$. \cr
$(1^m0^n)$ &$=(1\cdots 10\cdots 0)$. \cr
$(\{abc\}^n\ast)$ &$=(abcabc\cdots abc\ast)$. \cr
\end{tabular}

\subsection{Lattice VOA}
Let $L$ be a lattice with a bilinear form $\langle\cdot,\cdot\rangle$. 
Viewing $H={\R}\otimes_{\Z}L$ as a commutative 
Lie algebra with a bilinear form 
$<,>$, we define the affine Lie algebra 
$$\begin{array}{l}
\hat{H}=H[t,t^{-1}]+{\R}C\cr
[C,\hat{H}]=0, \quad [ht^n,h't^m]=\delta_{m+n,0}\langle h,h'\rangle C 
\end{array}$$ 
associated with $H$ 
and the symmetric tensor algebra $M(1)=S(\hat{H}^-)$ of $\hat{H}^-$, 
where $\hat{H}^-=H[t^{-1}]t^{-1}$.  As in 
\cite{FLM2}, we shall define 
the Fock space $V_L=\oplus_{x\in L}M(1)\iota(x)$ 
with the vacuum ${\bf 1}=\iota(0)$ 
and the vertex operators $Y(\ast,z)$ as follows: 
The vertex operator of $\iota(a)$ is given by 
$$ Y(\iota(a),z)={\exp}
\left(\sum_{n\in {\Z}_+}{a(-n)\over n}z^{n}\right){\exp}
\left(\sum_{n\in {\Z}_+}{a(n)\over -n}z^{-n}\right)
e^{a}z^{a} $$
and that of $a(-1)\iota(0)$ is 
$$ Y(a(-1)\iota(0),z)=a(z)=\sum a(n)z^{-n-1}.$$
The vertex operators of other elements are defined by the normal product:\\
$$  Y(a(n)v,z)=a(z)_nY(v,z)
={\rm Res}_x\{(x-z)^na(x)Y(v,z)-(z-x)^nY(v,z)a(x)\}.  $$
Here the operator of $a\otimes t^n$ on $M(1)\iota(b)$ are 
denoted by $a(n)$ and 
$$ \begin{array}{l}
a(n)\iota(b)=0  \mbox{ for } n>0 \cr
a(0)\iota(b)=<a,b>\iota(b) \cr
e^a\iota(b)=c(a,b)\iota({a+b}) \mbox{ for some cocycle } c(a,b)\cr
z^a\iota(b)=\iota(b)z^{<a,b>}. \end{array}$$
We note that the above definition of vertex operator is very general and 
so we may think 
$$  Y(v,z)\in {\End}(V_{\R\otimes L})\{z\}   $$
for $v\in \R\otimes_{\Z} L$, where 
$V_{\R\otimes_{\Z} L}=\sum_{a\in \R\otimes_{\Z} L}M(1)\iota(a)$. Set 
$\1=\iota(0)$. 
It is worthy to note that if we set $Y(v,z)=\sum_{n\in \R} v_nz^{-n-1}$, 
then $v_{-1}\iota(0)=v$ for any $v\in \R\otimes_{\Z} L$.

\subsection{$L(\hf,\st)\otimes L(\hf,\st)$}
In this subsection, we assume $L=\Z x$ with $\langle x,x\rangle=1$ and 
we don't use a cocycle $c(a,b)$ since $\iota(mx)$ is generated by one 
element $\iota(x)$ and $\iota(x)\in (V_L)_{\hf}$.  
As mentioned in \cite{DMZ}, we can find 
two mutually orthogonal conformal vectors 
$$\begin{array}{ll}
e^{+}(2x)={1\over 4}x(-1)^2{\bf 1}+{1\over 4}(\iota({2x})+\iota({-2x}))  
&\mbox{ and }\cr
e^-(2x)={1\over 4}x(-1)^2{\bf 1}-{1\over 4}(\iota({2x})+\iota({-2x}))  & 
\end{array} \eqno{(2.1)}$$
with central charge $\hf$ such that 
${\bw}=e^+(2x)+e^-(2x)={1\over 2}x(-1)^2{\bf 1}$ is the Virasoro 
element of $V_{2{\Z}x}$.
Let $\theta$ be the automorphism of $V_L$ induced from the automorphism 
$-1$ on $L$, which is given by $\theta(x(-n_1)\cdots x(-n_i)\iota(v))
=(-1)^ix(-n_1)\cdots x(-n_i)\iota(-v)$. We should note that $\theta$ 
is usually defined by $\theta(x(-n_1)\cdots x(-n_i)\iota(v))
=(-1)^{i+k}x(-n_1)\cdots x(-n_i)\iota(-v)$ for $\iota(v)\in (V_L)_k$, but 
we here have a half integer weight $k$. 
Take the fixed point space $(V_L)^{\theta}$ of $V_L$ by $\theta$. 
We note that each $e^{\pm}(2x)$ generates a simple vertex operator
subalgebra $<e^{\pm}(2x)>$ isomorphic to $L(\hf,0)$ since it is 
contained in $(V_{2{\Z}x})^{\theta}$, which has 
a positive definite invariant bilinear form as we will see 
in the next subsection.
As we mentioned in the introduction, 
$<\!e^{\pm}(2x)\!>\cong L(\hf,0)$ has only three irreducible modules 
$L(\hf,0), L(\hf,\hf),L(\hf,\st)$. 
By calculating the dimensions of weight spaces, 
there are no $L(\hf,\st)$ in $V_{L}$ 
since all elements $v\in V_L$ have 
integer or half integer weights. Since $\dim (V_L)_0=1$, 
$\dim (V_L)_1=1$, and $\dim (V_L)_{1/2}=2$, 
we conclude that $V_L$ is isomorphic to the direct sum of 
the tensor products 
$$ \left( L(\hf,0)\!\otimes \!L(\hf,0)\right)\oplus 
\left( L(\hf,0)\!\otimes\! L(\hf,\hf)\right)\oplus 
\left( L(\hf,\hf)\!\otimes\! L(\hf,0)\right)\oplus 
\left( L(\hf,\hf)\!\otimes\! L(\hf,\hf)\right) $$
as $<\!e^+(2x)\!> \otimes <\!e^-(2x)\!>$-modules by the actions of 
$e^{\pm}(2x)$ on $(V_L)_{\hf}$. 
Since $\theta$ fixes 
$e^{\pm}(2x)$ and $x(-1)(\iota(x)-\iota(-x))$, 
it keeps the above four irreducible 
$<\!e^+(2x)\!>\otimes<\!e^-(2x)\!>$-submodules invariant. 
Hence we obtain the decomposition:
$$(V_L)^{\theta}\cong \left(L(\hf,0)\otimes L(\hf,0)\right)\quad \oplus \quad 
\left(L(\hf,\hf)\otimes L(\hf,0)\right) $$
as $<\!e^+(2x)\!>\otimes <\!e^-(2x)\!>$-modules, see (4.11).
Take the subspace 
$M=\{ v\in (V_L)^{\theta}\ |\  {e^-(2x)}_1v=0\}$. Since $V_L$ is a SVOA, 
$M$ is a SVOA with the Virasoro element $e^+(2x)$ and we see 
$$M=M_0\oplus M_1, \quad M_0\cong L(\hf,0) \mbox{ and }
M_1\cong L(\hf,\hf) \eqno{(2.2)}$$
as $<\!e^+(2x)\!>$-modules. 
We note that 
$q=\iota(x)+\iota(-x)$ is a 
lowest degree vector of $M_1$ and $q_0q=2\iota(0)$. \par

It follows from the definition of vertex operators that 
$V_{2\Z x+\hf x}$ and $V_{2\Z x-\hf x}$ are irreducible $V_{2\Z x}$-modules. 
Hence, we have the following correspondence: 
$$ \begin{array}{|ll|l|}
\hline
 & &\theta \cr
\hline 
x(-1)\1 &\in L(\hf,\hf)\otimes L(\hf,\hf) & -1 \cr 
\iota(x)-\iota(-x) &\in L(\hf,0)\otimes L(\hf,\hf) & -1 \cr
\iota(x)+\iota(-x) &\in L(\hf,\hf)\otimes L(\hf,0) & +1 \cr
\iota(\pm x/2)&\in L(\hf,\st)\ots L(\hf,\st)
\oplus L(\hf,\st)\ots L(\hf,\st)  & \cr
\hline
\end{array} \eqno{(2.3)}$$  
Fix the lowest weight vectors $\iota({\hf x})$ and $\iota({-\hf x})$ 
of $V_{2{\Z}x+x/2}$ and $V_{2{\Z}x-x/2}$, respectively. 
By restricting $v$ in $M_{\bar{1}}\cong L(\hf,\hf)$ and taking 
the eigenspace $W$ of $e^-(2x)_1$ with an eigenvalue $\st$,  
$Y(v,z)$ defines the following three intertwining operators: 
$$\begin{array}{l}
I^{\hf,0}(\ast,z)\in I\pmatrix{L(\hf,\hf)\cr L(\hf,\hf) \quad L(\hf,0)},   \cr 
I^{\hf,\hf}(\ast,z)\in I\pmatrix{L(\hf,0) \cr L(\hf,\hf) 
\quad L(\hf,\hf)} \mbox{ and }\cr
I^{\hf,\st}(\ast,z)\in I\pmatrix{L(\hf,\st) \cr L(\hf,\hf) 
\quad L(\hf,\st)}. 
\end{array} \eqno{(2.4)} $$
Also, the restriction to $M_{\bar{0}}\cong L(\hf,0)$ 
defines the following intertwining operators: \\
$$ \begin{array}{ll}
I^{0,0}(\ast,z)\in I\pmatrix{L(\hf,0) \cr L(\hf,0) \quad L(\hf,0)}, &  \cr 
I^{0,\hf}(\ast,z)\in I\pmatrix{L(\hf,\hf) \cr L(\hf,0) 
\quad L(\hf,\hf)} \mbox{ and }\cr
I^{0,\st}(\ast,z)\in I\pmatrix{L(\hf,\st) \cr L(\hf,0) 
\quad L(\hf,\st)}, 
\end{array}   
\eqno{(2.5)}$$
which are actually module vertex operators of $<\! e^+(2x) \!>$.
We fix these intertwining operators throughout this paper. \par

We recall their properties from \cite{Mi3}. 

\begin{prn}
(1) The powers of $z$ in $I^{0,\ast}(\ast,z)$, 
$I^{\hf,0}(\ast,z)$ and $I^{\hf,\hf}(\ast,z)$ are all 
integers and those of $z$ in $I^{\hf,\st}(\ast,z)$ are half-integers, 
that is, in $\hf+{\Z}$. \\
(2) $I^{\ast,\ast}(\ast,z)$ satisfies the $L(-1)$-derivative property. \\
(3) $I^{\ast,\st}(\ast,z)$ satisfies 
the supercommutativity: \\
$$\begin{array}{l}
I^{0,\st}(v,z_1)I^{0,\st}(v',z_2) 
\sim I^{0,\st}(v',z_2)I^{0,\st}(v,z_1), \cr 
I^{0,\st}(v,z_1)I^{\hf,\st}(u,z_2) 
\sim I^{\hf,\st}(u,z_2)I^{0,\st}(v,z_1), \cr 
I^{\hf,\st}(u,z_1)I^{\hf,\st}(u',z_2) 
\sim -I^{\hf,\st}(u',z_2)I^{\hf,\st}(u,z_1), 
\end{array}  \eqno{(2.6)}$$
for $v,v'\in M_{\bar{0}}$ and $u,u'\in M_{\bar{1}}$. 
\end{prn}

\subsection{A lattice VOA with a positive definite invariant bilinear 
form} 
As we will see, 
we will gather the pieces from $\tilde{V}_{E_8}$ 
to construct $V^{\natural}$.  
In order to construct $V^{\natural}$ with a positive 
definite invariant form, we will show that 
there is a VOA $V_{E_8}$ over ${\R}$ with a positive 
definite invariant bilinear form. \\

We should note that $\tilde{V}_{E_8}$ in (1) is slightly different 
from an ordinary lattice VOA $V_{E_8}$ 
constructed from a lattice of type $E_8$. 
If we construct a lattice VOA $V_{E_8}$ over ${\R}$ by the construction 
in \cite{FLM2}, then 
$\iota(v)_{2k-1}\iota(v)\in S(\hat{H}^-)\iota(2v)\cap (V_{E_8})_0=\{ 0\}$ 
for any element $0\not=v\in L$ and $\langle v,v\rangle
=2k$ and so    
$\langle \iota(v),\iota(v)\rangle
=\langle \1,(-1)^k\iota(v)_{2k-1}\iota(v)\rangle
=0$. Namely, $V_{E_8}$ does not 
have a positive definite invariant bilinear form.

\begin{prn}  Let $L$ be an even lattice. Then there is a VOA $\tilde{V}_L$ 
which 
has a positive definite invariant bilinear form such that 
${\C}\ots \tilde{V}_L\cong {\C}V_{L}$. 
\end{prn}

\pr 
A lattice VOA $V_L=\boplus_{v\in L}S({\R}\otimes_{\Z}L^+)\iota(v)$
constructed by the lattice construction in \cite{FLM2} 
has an invariant bilinear form $\langle \ ,\ \rangle$. That is, 
it satisfies 
$$ \langle Y(a,z)u, v\rangle=\langle u, Y(e^{zL(1)}(-z^{-2})^{L(0)}a,z^{-1})v
\rangle  $$
for $a,u,v\in V_L$, see \cite{FHL}. 
$Y^{\dagger}(a,z)=Y(e^{zL(1)}(-z^{-2})^{L(0)}a,z^{-1}) 
=\sum a^{\dagger}_nz^{-n-1}$ is called the adjoint vertex operator. 
For $v\in {\R}\ots L$, identify it with $v(-1)\iota(0)\in (V_L)_1$. 
Since $L(1)v(-1)\iota(0)=0$ and 
$L(0)v(-1)\iota(0)=v(-1)\iota(0), Y^{\dagger}(v,z)=-z^{-2}Y(v,z^{-1})$ and so 
we have $v^{\dagger}(n)=-v(-n)$.  In the definition of $V_L$ 
in \cite{FLM2}, 
they used a group extension satisfying 
$ \iota(u')\iota(u)=(-1)^{<u',u>}\iota(u)\iota(u') \mbox{ and } 
\iota(v)\iota(-v)=\iota(0)$ for $\iota(v)\in (V_L)_k$.
Namely, $\iota(v)_{2k-1}\iota(-v)=\iota(-v)_{2k-1}\iota(v)=\iota(0)$.  

By definition, 
$Y^{\dagger}(\iota(v),z)=(-z^{-2})^{\langle v,v\rangle/2}Y(\iota(v),z^{-1})$. 
Hence, for $\iota(v)\in V_k$, we have  
$(\iota(v))^{\dagger}_n=(-1)^{k}(\iota(v))_{2k-n-2}$ 
and so  
$$\begin{array}{l}
\langle \iota(v)+\iota(-v),\iota(v)+\iota(-v)\rangle\iota(0) 
=(-1)^k(\iota(v)+\iota(-v))_{2k-1}(\iota(v)+\iota(-v)) \cr
=(-1)^k(\iota(v)_{2k-1}\iota(-v)+\iota(-v)_{2k-1}\iota(v))
=(-1)^k2\iota(0) 
\end{array}$$ 
and  
$$\langle \iota(v)-\iota(-v),\iota(v)-\iota(-v)\rangle=(-1)^{k+1}2\iota(0).$$ 
Let $\theta$ be an automorphism of $V_L$ induced from $-1$ on $L$, 
which is given by  
$$\theta(v^1(-i_1)\cdots v^m(-i_m)\iota(x))
=(-1)^{k+m}v^1(-i_1)\cdots v^m(-i_m)\iota(-x).$$ 
Hence, the space $V^0=(V_L)^{\theta}$ of $\theta$-invariants 
is spanned by the elements of the forms 
$$\begin{array}{l}
v^1(-n_1)...v^{2m}(-n_{2m})(\iota(v)+(-1)^k\iota(-v)) \qquad \mbox{ and } \cr
v^1(-n_1)...v^{2m+1}(-n_{2m+1})(\iota(v)-(-1)^k\iota(-v))
\end{array}$$ 
for $\iota(v)\in V_k$ and so $V^0$ has a positive definite invariant form. 
Similarly, $V^1=(V_L)^-$ has a negative definite invariant bilinear form,  
where $(V_L)^-=\{v\in V_L: \theta(v)=-v\}$.  
Since $V_L=V^0\oplus V^1$ has a ${\Z}_2$-grade, 
it is possible to denote the vertex operator of $v\in V^0$ by 
$\pmatrix{Y^{11}(v,z) &0 \cr 0&Y^{22}(v,z)}$ and the vertex operator 
of $u\in V^1$ by $\pmatrix{0& Y^{21}(u,z)\cr Y^{12}(u,z)&0}$, where 
$Y^{ij}(v,z)\in {\rm Hom}(V^i,V^j)[[z,z^{-1}]]$.  
Define new vertex operators by 
$$\tilde{Y}(v,z)=\pmatrix{Y^{11}(v,z) &0 \cr 0&Y^{22}(v,z)}$$ 
for $v\in V^0$ and 
$$\tilde{Y}(u,z)=\pmatrix{0& -Y^{21}(u,z)\cr Y^{12}(u,z)&0}$$ 
for $u\in V^1$.   
Then $(V,\tilde{Y})$ is a VOA with a positive definite 
invariant bilinear form.  This is the desired VOA. 
\prend

In the remaining of this paper, $\tilde{V}_{E_8}$ denotes the above 
VOA $(V_{E_8},\tilde{Y})$ 
with a positive definite invariant bilinear form. 
Since we mainly treat a VOA with a positive definite invariant 
bilinear form, we sometimes denote $V_L$ by 
$(\tilde{V}_L)^{\theta} \oplus \sqrt{-1}\tilde{V}_L^{-}$, 
where $\tilde{V}_L^{-}=\{v\in \tilde{V}_L: \theta(v)=-v\}$.

\section{Code VOAs with positive definite invariant bilinear forms}
In this section, we recall and prove several results from 
\cite{Mi2}$\sim$ \cite{Mi5}. 
We will first construct a code VOA $M_D$ 
with a positive definite invariant bilinear form for 
an even linear binary code $D$ of length $n$.  
Set $M_0=L(\hf,0)$ and $M_1=L(\hf,\hf)$.  
It is known that $F=M_0\oplus M_1$ has a super VOA structure $(F,Y^F)$, 
see (2.2). 
Although a SVOA structure on ${\C}F$ is uniquely determined, 
a SVOA structure on $F$ is not unique.  
Since $F$ has a $\Z_2$-grade, we can express a vertex operator 
$Y(v,z)$ by a $2\times 2$-matrix: 
$$ \begin{array}{lll}
 Y(v,z)=\pmatrix{Y^{00}(v,z)&0\cr 0&Y^{11}(v,z)} &for &v\in M_0, \cr
 Y(v,z)=\pmatrix{0& Y^{10}(v,z)\cr Y^{01}(v,z)&0 } &for &v\in M_1. 
 \end{array} $$
If we define new vertex operators $Y'(v,z)$ by  
$$ \begin{array}{lll}
 Y'(v,z)=\pmatrix{Y^{00}(v,z)&0\cr 0&Y^{11}(v,z)} &for &v\in M_0, \cr
 Y'(v,z)=\pmatrix{0& -Y^{10}(v,z)\cr Y^{01}(v,z)&0 } &for &v\in M_1,  
 \end{array} $$
then $(F,Y')$ is also a SVOA and it is not isomorphic to $(F,Y)$.  
So we choose one of them satisfying $q_{0}q\in {\R}^+\1$, where   
$q$ is a highest weight vector of $M_1$ and ${\R}^+=\{r\in \R|r>0\}$. \par
  
An essential property is a super-commutativity:
$$ Y^F(v,z_1)Y^F(u,z_2)\sim (-1)^{|v||u|}Y^F(u,z_2)Y^F(v,z_1)
\eqno{(3.1)}$$
for $|u|,|v|=0,1$ and 
$v\in M_{|v|}$ and $u\in M_{|u|}$. 
Here $A(z_1,z_2)\sim B(z_1,z_2)$ means  
$ (z_1-z_2)^NA(z_1,z_2)=(z_1-z_2)^NB(z_1,z_2)$ 
for a sufficiently large integer $N$. 
For a binary word $\al=(a_1,...,a_n)\in Z_2^n$, set 
$\hat{M}_{\al}=\ots_{i=1}^n M_{a_i}$, which is a subspace of 
$$F^{\ots n}=(M_0\oplus M_1)^{\otimes n}=
\boplus_{\al\in Z_2^n}\hat{M}_{\al}. $$ 
Define a vertex operator $Y^{\ots n}(u,z)$ of $u\in F^{\ots n}$ by 
$$  Y^{\ots n}(\ots_{i=1}^n v^i,z)(\otimes_{i=1}^n u^i)
=\ots_{i=1}^n (Y^F(v^i,z)u^i) \eqno{(3.2)} $$
for $u^i, v^i\in F$ and extend it to the whole space $F^{\ots n}$ linearly. 
It follows from (3.1) that for $v\in \hat{M}_{\al}$ and 
$u\in \hat{M}_{\be}$, we have the super commutativity: 
$$ Y^{\ots n}(v,z_1)Y^{\ots n}(u,z_2)
\sim (-1)^{\langle \al,\be\rangle}Y^{\ots n}(u,z_2)Y^{\ots n}(v,z_1).
\eqno{(3.3)}$$
Viewing $D$ as an elementary abelian 2-group with an invariant form, 
we shall use a central extension $\hat{D}=\{\pm e^{\al}:\al\in D\}$ 
of $D$ by $\pm 1$ in order to 
modify the supercommutativity (3.3). 
Let $\xi_i\ (i=1,...,n)$ denote a word $(0^{i-1}10^{n-i})$ and
$e^{\xi_i}$ a formal element satisfying 
$ e^{\xi_i}e^{\xi_i}=1$ and $e^{\xi_i}e^{\xi_j}=-e^{\xi_j}e^{\xi_i}$ 
for $i\not=j$.  For a word $\al=\xi_{j_1}+\cdots +\xi_{j_k}$ 
with $j_1<\cdots <j_k$, set 
$$e^{\al}=e^{\xi_{j_1}}e^{\xi_{j_2}}\cdots e^{\xi_{j_k}}. 
\eqno{(3.4)}$$
It is straightforward to check the following:
\begin{lmm}[\cite{Mi3}] 
For $\al,\be$, 
$$ \begin{array}{l}
e^{\al}e^{\be}=(-1)^{\langle\al,\be\rangle+|\al||\be|}e^{\be}e^{\al} \cr
e^{\al}e^{\al}=(-1)^{{k(k-1)\over 2}} \mbox{ for } |\al|=k 
\end{array} .
\eqno{(3.5)} $$
\end{lmm}

In order to combine (3.3) and (3.5), set
$$  M_{\delta}=\hat{M}_{\delta}\otimes e^{\delta} \eqno{(3.6)}$$  and
$$M_D=\boplus_{\de \in D} M_{\de}. \eqno{(3.7)}$$
Define a new vertex operator $Y(u,z)$ of $u\in M_D$ by setting 
 $$  Y(v\ots e^{\be},z)=Y^{\otimes n}(v,z)\ots{e^{\be}} \eqno{(3.8)}$$ 
for $v\ots e^{\be} \in M_{\be}=\hat{M}_{\be}\ots e^{\be}$ and 
extending it linearly. 
We then obtain the desired commutativity: 
$$ Y(v,z_1)Y(w,z_2)\sim Y(w,z_2)Y(v,z_1) \eqno{(3.9)}$$
for $v,w\in M_D$.
It is not difficult to see that  
$$ {\bw}
=\sum_{i=1}^n({\bf 1}^1\ots ...\ots {\bf 1}^{i-1}\ots {\bw}^i\ots 
{\bf 1}^{i+1}\ots ...\ots {\bf 1}^n)\ots e^0  \eqno{(3.10)}$$
is Virasoro element of $M_D$ and 
$$ {\bf 1}=({\bf 1}^1\ots ...\ots {\bf 1}^n)\ots e^0 \eqno{(3.11)} $$
is the vacuum of $M_D$, where ${\bf w}^i$ and 
${\bf 1}^i$ are Virasoro element and the vacuum of $M^i$, respectively. 
So we have proved the following theorem in \cite{Mi2}. \\

\begin{thm} 
If $D$ is an even binary linear code, then 
$(M_D, Y, {\bw}, {\1})$ 
is a simple VOA.
\end{thm}

It follows from the construction that $M_{\be+D}$ is an irreducible 
$M_D$-module 
and we will call it a {\it coset module} of $M_D$. 
From the choice of our cocycle, we can easily prove the following lemma. 
\begin{lmm}  If $g\in {\Aut}(D)$, there is an 
automorphism $\tilde{g}$ of a code VOA $M_D$ such that 
$\tilde{g}(e^i)=e^{g(i)}$ and $\tilde{g}(M_{\al})=M_{g(\al)}$. 
\end{lmm}

\pr  
For $g\in {\Aut}(D)$, we define a permutation 
$g_1$ on $\{\hat{M}_{\al}:\al\in D\}$ by 
$g_1(\otimes {M}_{a_i})=\otimes {M}_{a_{g(i)}}$ and 
an automorphism $g_2$ of $\hat{D}$ by 
$g_2(e^{\xi_{i_1}}\cdots e^{\xi_{i_t}})
=e^{\xi_{g(i_1)}}\cdots e^{\xi_{g(i_t)}}$.  
Combining the both action on $M_D=\sum_{\al\in D} 
\hat{M}_{\al}\otimes e^{\al}$, 
$\tilde{g}=g_1\otimes g_2$ becomes an automorphism of $M_D$.
\prend

We will next construct an invariant bilinear form on $M_D$. 
Let $(M, Y_M)$ be a module of $(V,Y)$.  
A bilinear form $\langle \cdot,\cdot\rangle$ on $M$ is said to be invariant 
\cite{FHL} if 
$$  \langle Y_M(a,z)u,v\rangle
=\langle u,Y_M(e^{zL(1)}(-z^{-2})^{L(0)}a,z^{-1})v\rangle
\mbox{ for } a\in V, u,v\in M, \eqno{(3.12)} $$
where $L(n)={\bw}_{n+1}$. 
It was proved in \cite{Li} that any invariant bilinear form 
on a VOA is automatically symmetric and there is a one-to-one correspondence
between invariant bilinear forms and elements of 
${\rm Hom}(V_0/L(1)V_1,{\R})$. 
Since $\dim V_0=1$ and $L(1)V_1=0$ for a code VOA $V=M_D$, 
there is a unique invariant bilinear 
form $\langle\cdot,\cdot\rangle$ satisfying $\langle\1,\1\rangle=1$. 
Using (3.12), it is given by 
$$
\langle u,v\rangle{\1}=\langle u_{-1}{\1},v\rangle{\1}
= Res_zz^{-1}(Y(e^{zL_1}(-z^{-2})^{L_0}u,z^{-1})v. 
\eqno{(3.13)}$$
Set $B=<L(1),L(0),L(-1)>$. 
Since $B\cong sl_2({\R})$ and $L(1)(M_D)_1=0$, 
$M_D$ is a direct sum of irreducible $B$-modules. 
Let $U$ be an irreducible $B$-submodule of $M_D$. Then 
there is an element $u\in (M_D)_k$ satisfying $L(1)u=0$ such that  
$U$ is spanned by $\{L(-1)^su:  s=0,1,...\}$. 
For any $v\in V_k$, 
$$
\langle u,v\rangle{\1}=\langle u_{-1}{\1},v\rangle{\1}
= Res_z(Y(((-1)^kz^{-2k})u,z^{-1})z^{-1}v=(-1)^{k}u_{2k-1}v. 
\eqno{(3.14)}$$
Also we note
$$\langle L(-1)^iu, L(-1)^jv\rangle=\langle L(-1)^{i-1}u, 
L(1)L(-1)^jv\rangle=(2kj+j^2-j)\langle L(-1)^{i-1}u, 
L(-1)^{j-1}v\rangle \eqno{(3.15)} $$
and $(2kj+j^2-j)>0$.  
Thus, $\langle\ ,\ \rangle$ is positive definite if and only if 
$$u_{2k-1}u\in (-1)^k{\R}^+{\1}\eqno{(3.16)}$$ 
for $0\not=u\in V_k$ satisfying $L(1)u=0$.

We first prove the ${\R}$-version of Theorem 4.5 in \cite{Mi3}. \\

\begin{prn}  
2) \quad Let $V=\boplus_{m=0}^{\infty} V_m$ be a simple VOA over 
${\R}$ with $\dim V_0=1$.  Assume that $V$ 
contains a set of mutually orthogonal 
conformal vectors $\{e^1,...,e^{n}\}$ 
so that the sum of them is the Virasoro element of $V$ and 
$\{e^1,...,e^{n}\}$ generates $T=L(\hf,0)^{\ots n}$. 
Assume further that $V$ has a positive definite invariant bilinear form and 
$\tilde{\tau}(V)=(0^n)$.  Then there is an even linear code $D$ such that 
$V$ is isomorphic to a code VOA $M_D$. 
\end{prn}

\pr  
Since $\tilde{\tau}(V)=(0^{n})$, $\tau_{e^i}=1$ and so we can 
define automorphism $\sigma_{e^i}$ for $i=1,...,n$. 
Set $Q=<\sigma_{e^i}: i=1,...,n>$. $Q$ is an elementary abelian 
2-group and let 
$$V=\oplus_{\chi\in Irr(Q)}V^{\chi}  $$
be the decomposition of $V$ into the direct sum of eigenspaces of $Q$. 
Since $\dim V_0=1$ and $V^{\chi}$ is an irreducible $V^Q$-module 
by \cite{DM2}, we have 
$V^{Q}=T$ and $V^{\chi}\cong \otimes L(\hf,h^i/2)$ 
as $T$-modules. 
Here $h^i\in \{0,1\}$ is given by $\chi(\sigma_{e^i})=(-1)^{h^i}$. 
Let $q$ denote a highest weight vector of $M_1$ 
such that $q_0q=\1\in M_0$. 
For a binary word $\al=(a^i)$, $q^{(a^i)}$ denotes 
$\otimes q^{a^i}\in M_{\al}$, where $q^{0}={\bf 1}$ and $q^1=q$. 
Identifying $\chi$ and $(h^i)$, 
$V^{\chi}\cong M_{\chi}\otimes \tilde{e}^{\chi}$ as $T$-modules such that 
$\langle q^{\chi}\otimes \tilde{e}^{\chi}, q^{\chi}\otimes \tilde{e}^{\chi})
\rangle=1$. 

Assume $|\chi|=2k$.
By the choice of $q^{\chi}\otimes \tilde{e}^{\chi}$ and 
$q^{\chi}_{2k-1}q^{\chi}={\bf 1}$, we have 
$$ \begin{array}{rl}
{\bf 1}=&\langle q^{\al}\otimes \tilde{e}^{\al}, 
q^{\al}\otimes \tilde{e}^{\al}\rangle\1\cr
=&\langle \1,
(-1)^k(q^{\al}\otimes \tilde{e}^{\al})_{2k-1}
q^{\al}\otimes \tilde{e}^{\al}\rangle\1\cr
=&\langle \1,(-1)^k\tilde{e}^{\al}\tilde{e}^{\al}\rangle \1. 
\end{array} $$ 
Hence, $\tilde{e}^{\al}\tilde{e}^{\al}=(-1)^{k}\tilde{e}^0$, which 
uniquely determine a cocycle that coincides with (3.17). 
This completes the proof of Proposition 3.1. 
\prend

As a corollary, we have 

\begin{cry}
For an even linear code $D$, 
$M_D$ has a positive definite invariant bilinear form. 
In particular, if $\al$ is even, 
then the coset module $M_{D+\al}$ also has a 
positive definite invariant bilinear form.  \\
\end{cry}

\pr
Recall that for a word $\al$ with $|\al|=2k$, 
say $\al=(1^{2k}0^{n-2k})$, 
$$ e^{\al}e^{\al}=e^{\xi_1}\cdots e^{\xi_{2k}}e^{\xi_1}\cdots e^{\xi_{2k}}
=(-1)^{k(2k-1)}=(-1)^k. \eqno{(3.17)}  $$
Let $S^n$ be the set of all even words of length n. 
Since all code VOAs are subVOAs of the code VOA $M_{S^n}$, 
it is sufficient to 
prove the assertion for the code $S^n$. 
Also, since $M_{S^n}\cong 
M_{S^n}\otimes ({\R}\1)^{\otimes n}\subseteq M_{S^{2n}}$ 
as sub VOAs, 
we may assume that $D$ is the set of all even words of length 2n. 
Let $\{x^1,...,x^n\}$ be an orthonormal basis of an Euclidian 
space of dimension $n$ and set 
$$L=\{\sum a_ix^i: a_i\in \Z, \sum a_i\equiv 0 \pmod{2}\}. 
\eqno{(3.18)} $$
Let $V_L$ be a lattice VOA constructed from $L$, (see \S 2.2).  
Let $\theta$ be an automorphism of $V_L$ induced from $-1$ on $L$ and 
decompose $V_L$ into $(V_L)^{\theta}\oplus (V_L)^-$, where 
$(V_L)^-=\{v\in V_L| \theta(v)=-v\}$. 
$(V_L)^{\theta}$ contains $2n$ mutually orthogonal rational conformal vectors 
$$  e(2x^i)^{\pm}
={1\over 4}x^i(-1)^2\1\pm {1\over 4}(\iota(2x^i)+\iota(-2x^i))
\eqno{(3.19)}$$
with central charge $\hf$ by (2.1). 
Set $\tilde{V}_L=(V_L)^{\theta}\oplus \sqrt{-1}(V_L)^-$. 
Then $\tilde{V}_L$ is a VOA with a  
positive definite invariant bilinear form containing 
$$T=<e(2x^i)^{\pm}:i=1,...,n>$$ 
by Proposition 2.2. 
Since $\langle v,2x^j\rangle\in 2\Z$ for $v\in L$, 
(2.3) implies $\tilde{\tau}(\tilde{V}_L)=(0^{2n})$. 
By Proposition 3.1, a code VOA $M_{S^{2n}}$ is isomorphic 
to $\tilde{V}_L$ which has a positive definite invariant bilinear 
form. 
\prend

If $\al\in D$ is a codeword of weight 2, say $\al=(110^{n-2})$, then 
$(M_{\al})_1\not=0$. Set $E=\{(00),(11)\}$, then 
$M_E$ is isomorphic to $V_{2\Z x}$ with $\langle x,x\rangle=1$ 
given in \S 2.3. We note $V_{\Z x}\cong M_{{\Z}_2^2}$ and 
${\exp}(\pi ix(0))$ keeps $V_{\Z x}$ invariant. We also note 
that $x(-1){\bf 1}\in (V_{2\Z x})_1$ and $x(0)=(x(-1){\bf 1})_0$. 
It follows from a direct calculation that 
$$ {\exp}(\pi ix(0))=(-1)^{\langle \be,(11)\rangle} \mbox{ on } M_{\be}  $$
for $\be\in {\Z}_2^2$.  
As long as a VOA $V$ contains a vector $v$ of 
weight 1, we can define automorphism $\exp(v_0)$ of ${\C}V$. 
Hence, we have the following lemma. 

\begin{lmm}  If a VOA contains a code $M_D$ and 
$D$ contains a codeword $\xi_i+\xi_j$ of weight 2, then 
${\C}V$ contains an automorphism $g$ such that 
$$ g=(-1)^{\langle \be,\xi_i+\xi_j\rangle} \mbox{ on } M_{\be}.  $$
In particular, it coincides with $\sigma_{e^i}\sigma_{e^j}$ on $M_D$.
\end{lmm}
  
\begin{cnj}  If a simple VOA $V$ contains a code VOA 
$M_D$ and $\be\in D$, 
then there is an automorphism $g$ of $V$ such that 
$g=\prod_{i\in Supp(\be)}\sigma_{e^i}$ on $M_D$. 
\end{cnj}

An important property of our cocycle is that if a maximal
self-orthogonal subcode $H$ of $D$ is doubly even, 
(for example, a Hamming code), 
then $\hat{H}=\{\pm e^{\al}:\al\in H\}$ 
is a maximal normal (elementary) abelian 2-subgroup of $\hat{D}$ and 
so every irreducible 
${\R}\hat{D}$-module is induced from a linear ${\R}\hat{H}$-module.  
In the remainder of this section, we assume that for 
an $M_D$-module $W$, one of maximal 
self-orthogonal subcode of $D_{\tilde{\tau}(W)}$ is doubly even and 
we denote it by $H$. 

We recall the structures of irreducible ${\C}M_D$-modules from \cite{Mi3}.  
Let ${\C}W$ be an irreducible $M_D$-module with $\tilde{\tau}({\C}W)=\mu$. 
If $H$ is a maximal self-orthogonal subcode of $D_{\mu}$ and 
$U'$ is an irreducible ${\C}M_H$-submodule of ${\C}W$, then   
the author showed in \cite{Mi3} that 
${\C}W={\Ind}_{H}^{D}(U')$ and every irreducible ${\C}M_H$-module 
is irreducible as a ${\C}T$-module.
We will show that these results also hold for an irreducible $M_D$-modules 
under the assumption that $H$ is doubly even. 

\begin{thm}  
Let $(X,Y^X)$ be an $M_D$-module with $\tilde{\tau}(X)=\mu$ 
and 
$\{ X^i:i=1,...,m\}$ the set of all non-isomorphic irreducible 
$T$-submodules of $X$.  Then there are 
representations $\phi^i:\hat{D}_{\mu}\to {\End}(Q^i)$ with 
$\phi^i(-e^0)=-I$ for $i=1,...,m$ such that 
$X\cong \oplus_{i=1}^m (X^i\otimes Q^i)$ as $M_{D_{\mu}}$-modules. 
Moreover, if $X$ is irreducible, then all $\phi^i$ are irreducible.  
For $\al\in D_{\mu}$, the module vertex operator $Y^X(q^{\al},z)$ of 
$q^{\al}=(\otimes q^{a_i})\otimes e^{\al}\in M_{\al}$ 
on $X^j\otimes Q^j$ is given by 
$$  \otimes_{i=1}^n I^{a_i/2,\ast}(q^{a_i},z)\otimes \phi^j(e^{\al}) $$
for $\al=(a_1,...,a_n)$. Here $q^0$ denotes the vacuum of $M_0$ and 
$q^1$ denotes the lowest degree vector $q$ in $M_1$. 
See \S 2.3 for $\otimes_{i=1}^n I^{a_i/2,\ast}(q^{a_i},z)$.   
\end{thm}

\pr
Let $U$ be a homogeneous component of $X$ generated by all $T$-submodules 
isomorphic to $X^1$ and 
let $U=\oplus_{i=1}^k U^i$ be a decomposition of $U$ into a 
direct sum of irreducible $T$-submodules $U^i$. 
By (1.1), $U$ is an $M_{D_{\mu}}$-module. Let $Q^1$ be the lowest 
degree space of $U$.  Since the dimension of 
the lowest degree space of $X^1$ is one, $\dim(U^i\cap Q^1)=1$ and 
$U\cong Q^1\otimes X^1$ 
as vector spaces.  Let $u^i$ be a nonzero lowest degree vector of $U^i$, then 
$\{u^1,...,u^k\}$ is a basis of $Q^1$. 
Let $\pi_j:U\to U^j=u^j\otimes X^1$ be a projection of $U$, that is, 
$\pi_j(u^i\otimes v)=\delta_{i,j}u^i\otimes v$ for $v\in X^1$.  
By (1.1), 
$Y^X(q^{\al},z)U\subset U[[z,z^{-1}]]$ for 
$q^{\al}\in M_{\al}$ and 
$\al\in D_{\mu}$.  
Since $\pi_j(Y^X(q^{\al},z)|_{u^i\otimes X^1})$ is an intertwining 
operator of type $\pmatrix{X^1 \cr M_{\al}\quad X^1}$ for $\al\in D_{\mu}$, 
the vertex operator $Y^X(q^{\al},z)_{|U}$ of $q^{\al}$ has an expression 
$$Y^X(q^{\al},z)=A(e^{\al})\otimes 
\left( (\otimes I)(\hat{q}^{\al},z)\right),$$ 
where $A(e^{\al})$ is 
a $k\times k$-matrix acting on $Q^1={\R}u^1\oplus \cdots {\R}u^k$ 
and $(\otimes I)(\hat{q}^{\al},z)$ is 
the tensor product $\otimes I^{({a_i\over 2},\ast)}(q^{a_i},z)$ 
of the fixed intertwining operators in \S 2.3 for 
$\hat{q}^{\al}=\otimes q^{a_i}$. 
We note that $Y^X$ is uniquely determined by $\{Y^X(q^{\al},z):\al\in D\}.$ 
Since $Y^X(q^{\al},z)$ satisfies the commutativity
$$  Y^X(q^{\al},z)Y^X(q^{\be},w)\sim Y^X(q^{\be},w)Y^X(q^{\al},z) $$
and 
$(\otimes I)(\hat{q}^{\al},z)$ satisfies the super-commutativity
$$ 
(\otimes I)(\hat{q}^{\al},z)(\otimes I)(\hat{q}^{\be},w)
\sim (-1)^{\langle\al,\be\rangle}
(\otimes I)(\hat{q}^{\be},w)(\otimes I)(\hat{q}^{\al},z), $$
we obtain 
the supercommutativity : \\
$$  A(e^{\al})A(e^{\be})=(-1)^{\langle \al,\be\rangle}A(e^{\be})A(e^{\al}). $$
Moreover, since $Y^X(\ast,z)$ satisfies the associativity
$$  Y^X(q^{\al}_mq^{\be},z)
=Y^X(q^{\al},z)_mY^X(q^{\be},z) $$ 
and 
$(\otimes I)(\hat{q}^{\al},z)$ satisfies 
the superassociativity by \cite{Mi3}, we have 
$$\begin{array}{rl}
A(e^{\al}e^{\be})(\otimes I)(\hat{q}^{\al}_m\hat{q}^{\be},z) 
=&Y^X(q^{\al}_mq^{\be},z) \cr
=&Y^X(q^{\al},z)_mY^X(q^{\be},z) \cr
=&Res_{x}\{(x-z)^mA(e^{\al})(\otimes I)(\hat{q}^{\al},x)A(e^{\be})(\otimes I)
(\hat{q}^{\be},z)\cr
&-(-z+x)^mA(e^{\be})(\otimes I)(\hat{q}^{\be},z)A(e^{\al})
(\otimes I)(\hat{q}^{\al},x)\} \cr
=&Res_{x}\{(x-z)^mA(e^{\al})A(e^{\be})(\otimes I)(\hat{q}^{\al},x)(\otimes I)
(\hat{q}^{\be},z)\cr
&-(-z+x)^mA(e^{\be})A(e^{\al})(\otimes I)(\hat{q}^{\be},z)
(\otimes I)(\hat{q}^{\al},x)\} \cr
=&A(e^{\al})A(e^{\be})Res_{x}\{(x-z)^m(\otimes I)(\hat{q}^{\al},x)(\otimes I)
(\hat{q}^{\be},z)\cr
&-(-1)^{\langle \al,\be\rangle}
(-z+x)^m(\otimes I)(\hat{q}^{\be},z)(\otimes I)(\hat{q}^{\al},x)\} \cr
=&A(e^{\al})A(e^{\be})(\otimes I)(\hat{q}^{\al},z)_{(m)}
(\otimes I)(\hat{q}^{\be},z) \cr
=&A(e^{\al})A(e^{\be})(\otimes I)(\hat{q}^{\al}_m\hat{q}^{\be},z). 
\end{array} $$ 
Hence we have the associativity : 
$$A(e^{\al})A(e^{\be})=A(e^{\al}e^{\be}) $$ and  
$A(e^\al)A(e^\al)=(-1)^{|\al|/2}I$ for all $\al, \be\in D_{\mu}$,  
where $I$ is the identity matrix.  
Hence $A$ is a matrix representation of the central 
extension $\hat{D}_{\mu}$ on $Q^1$.   
We next assume that $X$ is irreducible. 
Let $Q^0$ be an irreducible $\hat{D}_{\mu}$-submodule and  
$W$ the subspace spanned by $\{ v_nw: v\in M_D, w\in Q^0\otimes X^1, 
n\in {\Z}\}$. 
Proposition 4.1 in \cite{DM2} implies $X=W$. 
On the other hand, the tensor product  
$M_{\be+D_{\mu}}\times (Q^0\otimes X^1)$ does not contains a submodule 
isomorphic to $X^1$ for $\be\not\in D_{\mu}$ by (1.1) and so 
$U=W\cap U=Q^0\otimes X^1$. 
Hence, $Q^1$ is an irreducible $\hat{D}_{\mu}$-module on which 
$-e^0$ acts as $-1$.
\prend
  
As a corollary, we have the followings:

\begin{cry}  If $D$ is a doubly even code and $W$ is an irreducible 
$M_D$-module, then $W$ is also irreducible as a $T$-module. 
\end{cry}

Note that if $H$ is a maximal self-orthogonal subcode (doubly even) of $D$, 
then every irreducible ${\R}\hat{D}$-module $W$ is 
induced from a ${\R}\hat{H}$-module and $W$ is a direct sum 
of distinct irreducible ${\R}\hat{H}$-modules.
Hence, we have the following corollary. 

\begin{cry}  If $(W, Y^W)$ is an irreducible $M_D$-module with 
$\tilde{\tau}(W)=(1^n)$, then there is an irreducible 
representation $\phi:\hat{D}\to {\End}(Q)$ satisfying $\phi(-e^0)=-I$ 
such that $W\cong L(\hf,\st)^{\otimes n}\otimes Q$ as $M_D$-modules. 
Here the module vertex operator $Y^X(q^{\al},z)$ of 
$q^{\al}=(\otimes q^{a_i})\otimes e^{\al}\in M_{\al}$ 
on $L(\hf,\st)^{\otimes n}\otimes Q$ is given by 
$$  \otimes_{i=1}^n I^{a_i/2,\st}(q^{a_i},z)\otimes \phi(e^{\al}) $$
for $\al=(a_1,...,a_n)$. 
In particular, $Y^W$ is uniquely determined by an irreducible 
$M_H$-submodule. 
\end{cry} 

Conversely, we will prove the following proposition:

\begin{prn}
Let $\mu$ be a word such that $\langle D,\mu\rangle=0$. 
Assume 
that $H$ is a maximal self-orthogonal (doubly even) subcode of $D_{\mu}$ 
and $U$ is an irreducible $M_H$-module with $\tilde{\tau}(U)=\mu$.  
Then there is an irreducible $M_D$-module $W$ containing $U$ 
as an $M_H$-submodule. 
\end{prn}

\pr
We may assume $\mu=(0^{n-m}1^{m})$.
By the above lemmas, there is a binary word $(a_1,...,a_{n-m})$ such that 
$U=(\otimes_{i=1}^{n-m} L(\hf,{a^i\over 2}))\otimes 
(L(\hf,\st)^{\otimes m})\otimes {\R}_{\chi}$. 
Since $D\subseteq <\mu>^{\perp}$ and $D$ is even, 
$D \subseteq S_{n-m}\oplus S_m$, 
where $S_r$ denotes the set of all even words of length $r$.  
If $n=m$, then  $L(\hf,\st)^{\otimes m}\otimes 
{\Ind}_{\hat{H}}^{\hat{D}}({\R}_{\chi})$ 
is the desired $M_D$-module.  
If $m=0$, then a coset module $M_{(a^i)+D}$ is the desired $M_D$-module. 
For general cases, let $K$ be a maximal self-orthogonal 
subcode of $S_m$ containing $H$ and choose an irreducible $\hat{K}$-module 
$Q$ containing ${\R}_{\chi}$.   the tensor product 
$M_{(a^i)+S_{n-m}}\otimes {\Ind}_{\hat{K}}^{\hat{S}_m}(
L(\hf,\st)^{\otimes m}\otimes Q)$ is an $M_{S_{n-m}\oplus S_m}$-module 
containing $U$. By Theorem 3.2, there is an irreducible $M_D$-submodule 
containing $U$, which is the desired $M_D$-module.  
\prend

Our next aim is to prove that an $M_D$-module $W$ satisfying the 
above condition is uniquely determined. We will call it an 
{\it induced module} and denote it by ${\rm Ind}_H^D(U)$. 
Applying Proposition 11.9 in \cite{DL} into our case, 
we have the following lemma (see \cite{Mi3}).

\begin{lmm}   Let $E$ be a subcode of $D$. 
Let $W^1, W^2, W^3$ be irreducible $M_D$-module 
and $U^1, U^2$ irreducible $M_E$-submodules of $W^1$ and $W^2$, respectively, 
then there is an injection map:
$$ \phi\ :\ I_{M_D}\pmatrix{W^3 \cr W^1 \quad W^2}
\to I_{M_E}\pmatrix{W^3\cr U^1\quad U^2}. $$
\end{lmm}

Using the above lemma, we will prove the 
following fusion rule, which was proved in \cite{Mi3} for 
a code VOA ${\C}M_D$.  

\begin{thm}
If $X$ is an irreducible $M_D$-module, 
then the fusion product 
$$ M_{\al+D} \times X  $$
is an irreducible $M_D$-module for any $\al$.   
\end{thm}

Set $\mu=\tilde{\tau}(X)$.  We will first prove the following lemmas, 
(Lemma 3.5$\sim$ 3.7). 

\begin{lmm}  Assume that $D$ is a doubly even code and 
$Supp(D)\subseteq Supp(\mu)$.  
Then $M_{\al+H}\times X$ is irreducible.  
\end{lmm}

\pr
By Corollary 3.2 and 3.3, 
$W\cong (\otimes L(\hf,h^i))\otimes {\R}_{\chi}$ as $M_D$-modules, 
where $\chi$ is a linear representation of $\hat{D}$ on ${\R}_{\chi}$. 
Let $U$ be an irreducible $M_H$-module so that 
$0\not=I_{M_H}\pmatrix{U \cr M_{\al+H} \quad W}$.  Clearly, 
$\tilde{\tau}(U)=\tilde{\tau}(W)$ and so 
$U\cong \otimes L(\hf,h^i+{a_i\over 2})\otimes {\R}_{\phi}$ 
for some linear representation $\phi$ of $\hat{D}$.  
By Lemma 3.4,
there is an injective map 
$$ \pi:I_{M_H}\pmatrix{U \cr M_{\al+H} \quad W} \to 
I_T\pmatrix{U \cr M_{\al}\quad W}.$$  Since $M_{\al}=\otimes L(\hf,a_i/2)$,   
$M_{\al}\times W$ is irreducible as a $T$-module by (1.1) and so 
$M_{\al}\times W\cong U$ as $T$-modules. 
We fix a nonzero intertwining operator 
$J(\ast,z)\in I_T\pmatrix{U\cr M_{\al} \quad W}$. 
Then for any intertwining operator 
$I(\ast,z)\in I\pmatrix{U\cr M_{\al+H}\quad W}$ we may assume 
$I(v,z)=J(v,z)$ for $v\in M_{\al}$ by multiplying a scalar. 
Since $I$ satisfies the commutativity: 
$$ I^{\otimes n}(q^{\al},x)\phi(e^{\al})I(v,z)
\sim I(v,z)I^{\otimes n}(q^{\al},x)\chi(e^{\al})  $$
with the module vertex operators, 
$\phi$ is uniquely determined by $\chi$ 
and so $M_{\al+D}\times W$ is irreducible. 
\prend

\begin{lmm}  Assume that $Supp(D)\subseteq Supp(\mu)$.  
Then $M_{\al+D}\times W$ is irreducible.  
\end{lmm}

\pr
Let $H$ be a maximal orthogonal (doubly even) subcode of $D$. 
By Corollary 3.2 and 3.3,  
$W\cong \otimes L(\hf,h^i)\otimes Q_{\chi}$ as $M_D$-modules, 
where $\chi$ is a representation of $\hat{D}$ on $Q_{\chi}$ such that 
$\chi(-e^0)=-I$.  Since $\hat{H}$ is a maximal normal abelian subgroup 
of $\hat{D}$, there is a $\hat{H}$-submodule $Q^0$ such that 
${\Ind}_{\hat{H}}^{\hat{D}}(Q^0)=Q_{\chi}$.  
Let $U$ be an irreducible $M_D$-module such that 
$0\not=I_{M_D}\pmatrix{U \cr M_{\al+D} \quad 
W}$.  Clearly, $\tilde{\tau}(U)=\mu$ and so $U\cong \otimes L(\hf,k^i)\otimes 
Q_{\phi}$ for some $\hat{D}$-module $Q_{\phi}$.  
By Lemma 3.6, 
$M_{\al+H}\times (\otimes L(\hf,h^i)\otimes Q^0)$ 
is irreducible and so $U$ contains an $M_H$-module 
$M_{\al+H}\times (\otimes L(\hf,h^i)\otimes Q^0)$.  
Therefore, $U$ is uniquely determined. 
Since $Q_{\phi}$ is a direct sum of distinct irreducible 
$\hat{H}$-modules, 
 $\dim I\pmatrix{U \cr M_{\al+D}\quad W}=1$ and so 
 $M_{\al+D}\times W$ is irreducible. 
\prend

\begin{lmm}  Let $(W,Y^W)$ be an irreducible $M_D$-module with 
$\tilde{\tau}(W)=\mu$ and let $W=\oplus_{i=1}^r U^i$ be an 
decomposition of $W$ into 
the direct sum of distinct homogeneous $M_{D_{\mu}}$-submodules $U^i$.  
Then $U^i$ is irreducible and $Y^W$ is uniquely determined 
by an $M_{D_{\mu}}$-module $U^i$ for any $i$. 
\end{lmm}

\pr 
Let $X$ be an irreducible $M_{D_{\mu}}$-submodule of $U^1$ and set 
$X\cong \otimes L(\hf,h_i)$. 
By (1.1), $U^0$ is homogeneous as a $T$-module, that is,   
every irreducible $T$-submodule of $U^0$ is isomorphic to $X$. 
By \cite{DM2}, $\{ v_mu: u\in X, v\in M_{\al}, \al\in D \}$ spans 
$W$.  On the other hand, 
if $\al=(a_i)\not\in D_{\be}$, then irreducible $T$-submodule 
generated by $v_mu$ is isomorphic to 
$\otimes L(\hf,h_i+{a_i\over 2})$ and so 
$<v_mu:u\in X, v\in M_{\al},\al\in D>\cap U^0=X$, which proves $U^0=X$.
Clearly, $<v_mu: u\in U^0, v\in M_{\al+D_{\mu}}>$ is an 
irreducible $M_{D_{\mu}}$-module $U^{\al}$ by the same argument.  
Lemma 3.8 implies that $M_{\al+D_{\mu}}\times U^0$ is irreducible. 
Since the restriction $Y(v,z):U^0 \to U^{\al}[[z,z^{-1}]]$ 
for $v\in M_{\al+D_{\mu}}$ 
is a nonzero intertwining operator, we conclude  
$M_{\al+D_{\mu}}\times U^{\be}=U^{\al+\be}$. 
Namely, if one of $\{U^i:i=1,...,r\}$ is given, then the other $U^j$ are 
uniquely determined as $M_{D_\mu}$-modules. 
By Proposition 3.2, there is at least one $M_D$-module $S$ such that 
$S=\oplus_{\be\in D/D_{\mu}} U^{\be}$. Let $Y^S$ be the module 
vertex operator of $S$ and set 
$I^{\al,\be}(\ast,z)=Y^W(\ast,z):U^{\be}\to U^{\al+\be}$ for 
$v\in M_{\al+D_{\be}}$ by restriction.  
Since 
$\dim I\pmatrix{U^{\al+\be}\cr M_{D_{\mu}+\al}\quad U^{\be}}
=1$, $J^{\al,\be}(v,z)=\la_{\be,\be+\al}I^{\al,\be}(v,z)$ for 
$v\in M_{\al+D_{\mu}}$.  Let $A(\al)$ be a matrix 
$( \la_{\be,\be+\al})$ whose $(\be,\be+\al)$-entry is  
$\la_{\be,\be+\al}$ and 0 otherwise. 
Since $I$ and $J$ satisfy the mutually commutativity and 
the associativity, respectively, 
$A: D/D_{\mu} \to M(n\times n, \R)$ is a regular representation and 
so we can reform $A(\al)$ into a permutation matrix by 
changing the basis.  Thus, $J^{\al,\be}=I^{\al,\be}$ 
and so $W$ is isomorphic to $S$ as an $M_D$-module. 
\prend

\noindent
$[${\bf Proof of Theorem 3.3}$]$  \qquad 
Set $\al=(a_i)$ and let $\mu=\tilde{\tau}(W)$. 
Let $H$ be a maximal self-orthogonal (doubly even) 
subcode of $D_{\mu}$. 
Let $W^1$ be an irreducible $M_{D_{\mu}}$-submodule of $W$ and  
let $U$ be an irreducible $M_D$-module such that 
$I\pmatrix{U\cr M_{\al+D}\quad W}\not=0$.  Clearly, $\tilde{\tau}(U)=\mu$. 
By \cite{DM2}, there is an injective map:
$$ \pi: I\pmatrix{U \cr M_{\al+D}\quad W}\to 
I_T\pmatrix{U\cr M_{\al+D_{\mu}} \quad 
W^0}.  $$
By Lemma 3.7, $W^1=M_{\al+D_{\mu}}\times W^0$ is irreducible and so 
$U$ contains $W^1$. Since $W^1$ determines $U$ uniquely and $U$ 
contains only one irreducible $M_{D_{\mu}}$-submodule isomorphic to 
$W^1$,  
$M_{\al+D_{\mu}}\times W=U$. 
\prend

Combining the above arguments, we have the following theorem:

\begin{thm}
Let $W$ be an irreducible $M_D$-module with $\tilde{\tau}(W)=\mu$.  
Let $E$ be an even linear code containing $D$ and assume 
$\langle E,\mu\rangle=0$. 
Assume that there is a maximal self-orthogonal 
(doubly even) subcode $H$ of $D_{\mu}$ and $H$ is also a maximal 
self-orthogonal in $E_{\mu}$. 
Then there is a unique irreducible $M_E$-module $X$ containing $W$ 
as an $M_D$-submodule.  
\end{thm}

We will call $X$ as an {\it induced $M_E$-module} and denote it by 
${\rm Ind}_D^E(W)$. 

We next quote the results about Hamming code VOA from \cite{Mi2}. 
In this paper, a Hamming code means a $[8,4,4]$-Hamming code. 
Let $H$ be a Hamming code and $\{e^1,...,e^8\}$ be a set of 
coordinate conformal vectors of a Hamming code VOA $M_{H_8}$. 
Let $W$ be an irreducible $M_{H_8}$-module.  If $\tilde{\tau}(W)=(0^8)$, 
then $W$ is isomorphic to a coset module $M_{H_8+\al}$.  We denote it 
by $H(\hf,\al)$. If $\tilde{\tau}(W)=(1^8)$, then there is a linear 
representation $\chi:\hat{H_8}\to \{\pm 1\}$ such that 
$W$ is isomorphic to $(L(\hf,\st)^{\otimes 8})\otimes {\R}_{\chi}$. 
If we fix a basis $\{\al^1,\al^2,\al^3,\al^4\}$ of $H_8$, then 
there is a word $\be$ such that 
$\chi(\al^i)=(-1)^{\langle \be,\al^i\rangle}$.  
We denote $W$ by $H(\st,\be)$.  We should note that $H(\st,\be)$ 
depends on the choice of a basis of $H_8$.  So, we will 
fix a basis $\{(1^8),(1^40^4),(1^20^21^20^2),((10)^4)\}$ in this paper.  
Namely, we have the following result. 

\begin{thm}
Let $W$ be an irreducible $M_{H_8}$-module. 
If $\tilde{\tau}(W)=(0^8)$, then 
$W$ is isomorphic to one of 
$$\{ H(\hf,\al): \al\in {\Z}_2^8 \}.$$
If $\tilde{\tau}(W)=(1^8)$, then $W$ is isomorphic to one of 
$$\{ H(\st,\al): \al\in {\Z}_2^8 \}.$$
$H(\hf,\al)\cong H(\hf,\be)$ if and only if $\al+\be\in H_8$ and 
$H(\st,\al)\cong H(\st,\be)$ if and only if $\al+\be\in H_8$. 
$H(\hf,\al)$ is a coset module $M_{\al}$ and 
$H(\st,\be)$ is isomorphic to $L(\hf,\st)^{\ots 8}$ as 
an $L(\hf,0)^{\ots 8}$-module. 
\end{thm}

In \cite{Mi5}, the author obtain the following fusion rules. 

\begin{lmm}\cite{Mi2}
$$\begin{array}{l}
H(\hf,\al)\times H(\hf,\be)=H(\hf,\al+\be) \cr
H(\st,\al)\times H(\hf,\be)=H(\st,\al+\be) \cr
H(\st,\al)\times H(\st,\be)=H(\hf,\al+\be) 
\end{array}. $$
\end{lmm}

The proof is based on the 
nice properties of the Hamming code VOA $M_{H_8}$. 
To simplify the notation, 
we will choose another cocycle of $\hat{H_8}$ for a while. 
Set $\bar{e}^{\al}=e^{\la_1\al_1}\cdots e^{\la_4\al_4}$ for 
$\al=\la_1\al_1+\cdots+\la_4\al_4\in H_8$, where $\{\al_1,...,\al_4\}$ 
is a fixed basis of $H_8$.  
In ${\cal H}_8$, there are 14 words of weight 4. 
For such a codeword (or a 4 points set) $\alpha$, set
$$\bar{q}^{\al}={1\over 4}(\otimes_{i=1}^8 q^{a_i})\otimes \bar{e}^{\al}.  $$
It follows from a direct calculation that 
$$s^{\al}={1\over 8}(e^1+...+e^8)
+{1\over 8}\sum_{\be\in H_8, \ |\be|=4}(-1)^{(\al,\be)}\bar{q}_{\be}  $$
is a conformal vector with central charge 
$\hf$ for a word $\alpha$ in \cite{Mi2}. 
Clearly, $s^{\al}=s^{\be}$ if and only if $\al+\be\in H_8$. 
It is also straightforward to check that $\langle s^{\al},s^{\be}\rangle=0$ 
if and only if $\al+\be$ is even word. 
Therefore, there are two other sets of coordinate 
conformal vectors $\{d^1,...,d^8\}$ and $\{f^1,...,f^8\}$ in $M_{H_8}$. 
By the definition of the set of coordinate conformal vectors, 
$T_d=<d^1,...,d^8>$ and $T_f=<f^1,...,f^8>$ are coordinate sets of 
conformal vectors. 
Viewing an $M_{H_8}$-module as a $T_d$-module and a $T_f$-modules, 
we have the following correspondence:
(see Proposition 2.2 and Lemma 2.7 in \cite{Mi5}). 

\begin{lmm}
There are the other two sets of coordinate 
conformal vectors $\{d^1,...,d^8\}$ and $\{f^1,...,f^8\}$ in $M_{H_8}$ 
such that  
$$ \begin{array}{lclcl}
H(\hf,(0^8)) \mbox{ w.r.t.} <e^i> &\leftrightarrow 
& H(\hf,(0^8)) \mbox{ w.r.t.} <d^i> 
&\leftrightarrow& H(\hf,(0^8)) \mbox{ w.r.t.} <f^i>\cr
H(\hf,\xi_1) \mbox{ w.r.t.} <e^i> &\leftrightarrow 
& H(\st,(0^8)) \mbox{ w.r.t.} <d^i> 
&\leftrightarrow& H(\st,\xi_1) \mbox{ w.r.t.} <f^i> \cr
H(\st,(0^8)) \mbox{ w.r.t.} <e^i>  &\leftrightarrow 
& H(\st,\xi_1) \mbox{ w.r.t.} <d^i> 
&\leftrightarrow& H(\hf,\xi_1) \mbox{ w.r.t.} 
<f^i> \cr
H(\st,\xi_1) \mbox{ w.r.t.} <e^i> &\leftrightarrow 
& H(\hf,\xi_1) \mbox{ w.r.t.} <d^i> 
&\leftrightarrow& H(\st,(0^8)) \mbox{ w.r.t.} <f^i>, 
\end{array} $$
where $\xi_1$ denotes $(10^7)$.
\end{lmm}

As a corollary, we have: 

\begin{cry}
If $W$ is an irreducible $M_{H_8}$-module, then 
$W\times H(\hf,\al)$ and $W\times H(\st,\al)$ are irreducible 
for any $\al\in {\Z}_2^8$.  
\end{cry}

We next recall the following important theorem from \cite{Mi5} 
and prove it as a corollary of the above results.  

\begin{thm}
Let $W^1$ and $W^2$ be irreducible $M_D$-modules 
and assume that the pair $(D,<\tilde{\tau}(W^1),\tilde{\tau}(W^2)>)$ 
satisfies (1) and (2) of Hypotheses I. Then 
$W^1\times W^2$ is irreducible.  
\end{thm}

\pr 
Set $\tilde{\tau}(W^1)=\al$ and $\tilde{\tau}(W^2)=\be$. 
If $\al=0$ or $\be=0$, then the assertion follows from Theorem 3.3. 
We may assume $\al=(1^{8r}0^{8s})$. 
Let $U$ be an irreducible $M_D$-module so that 
$0\not= I\pmatrix{U \cr W^1\quad W^2}$. 
Clearly $\tilde{\tau}(U)=\al+\be$. 
By Hypotheses I, 
there is a self-dual subcode $E=E_{\al}\oplus E_{\al^c}$ of $D$ such that 
$E$ is a direct sum of Hamming codes. Assume that 
$E_{\be}$ is a direct factor of $E$. Then $E=E_{\be}\oplus E_{\be^c}$.
Let $U^i$ be irreducible $M_{E}$-submodules 
of $W^i$ for $i=1,2$. By Theorem 3.5, 
$U^1\cong (\otimes_{i=1}^r H(\st,\al^i))\otimes (
\otimes_{j=1}^s H(\hf,\be^j))$ as $M_E$-modules and so 
$U^1\times U^2$ is irreducible.  
Since $U$ contains $U^1\times U^2$, $U$ is uniquely determined. 
Since $U$ is a direct sum of distinct irreducible $M_E$-submodules,   
we have that $W^1\times W^2=U$ is irreducible. 

So we may assume that $E_{\be}$ is not a direct factor of $E$. 
By Hypotheses I, there are maximal self-orthogonal 
subcodes $H_{\be}$ and $H_{\al+\be}$ of $D_{\be}$ and $D_{\al+\be}$ 
containing $E_{\be}$ and $E_{\al+\be}$, respectively, such that 
$H_{\be}+E=H_{\al+\be}+E$. 
Since $H_{\be}+E$ satisfies Hypotheses I for $<\al,\be>$ and 
$W^i$ and $U$ are direct sums of distinct irreducible 
$M_{H_{\be}+E}$-modules, we may assume $D=H_{\be}+E=H_{\al+\be}+E$.  
We first assert the following: \\

\noindent
{\bf Claim:}\qquad {\it 
$W^2$ and $U$ are irreducible as $M_E$-modules. }\\

Since the proofs are almost the same, we will prove the assertion for $W^2$. 
Set $E=E_1\oplus E_k$, where $E_i\cong H_8$.   
Assume first that $E_{\be}$ contains a direct factor of $E$, say $E_1$. 
Namely, assume $\be=(1^8...)$. Then $\al+\be=(0^8...)$. Let 
$\pi_{\be}:(a^i)\to (a^i)_{i\in Supp(\be)}$ be a projection. Since 
$\pi_{(1^80^{8k-8})}(D)=\pi(H_{\al+\be}+E)=E_1$, 
$D=E_1\oplus D_{(0^81^{8k-8})}$ and so 
it is sufficient to prove the assertion for $D_{(0^81^{8k-8})}$.  
By the induction and $\langle \be,\al\rangle=0$ for $\al\in D$, 
we may assume $\be=(1^40^41^40^4...1^40^4)$. 
Since $H^{\be}$ contains $E_{\be}$ and $D=H^{\be}+E$, $D_{\be}=H^{\be}$.  
Let $X$ be an irreducible $M_{H^{\be}}$-submodule of $W^2$, then 
$W^2={\Ind}_{H^{\be}}^D(X)$ and $X$ is irreducible as an $T$-module. 
In particular, $X$ is irreducible as an $M_{E_{\be}}$-module with 
$\tilde{\tau}(X)=\mu$.   
Hence ${\Ind}_{E_{\be}}^E(X)$ is an irreducible $M_E$-submodule of $W^2$.   
On the other hand, since $D/H^{\be}\cong E/E_{\be}$, 
$\dim {\Ind}_{E_{\be}}^E(X)=\dim W^2$, which proves the claim.\\ 

We now go back to the proof of Theorem 3.6.  
Set $\ga=\al+\be$. Let $X$ be an irreducible $M_E$-submodule of $W^1$. 
Since $W^2$ and $U$ are both irreducible $M_E$-modules by 
the above claim,   
$$\dim I_{M_D}\pmatrix{U\cr W^1\quad W^2}
\leq \dim I_{M_E}\pmatrix{U\cr X \quad W^2}=1 $$
and so $U\cong X\times W^2$ as $M_E$-modules.  
Fix a nonzero intertwining operator 
$$I^1(\ast,z)\in I_{M_E}\pmatrix{W^3\cr X\quad W^2}.$$ 
For 
$I(\ast,z)\in I_{M_D}\pmatrix{W^3\cr W^1\quad W^2}$, 
there is a scalar $\la$ 
such that $I(v,z)=\la I^1(v,z)$ for $v\in X$. 
$Y^U(u,z)I(v,z)\sim I(v,z)Y^2(u,z)$ and so 
$Y^U(u,z)I^1(v,z)=I^1(v,z)Y^2(u,z)$ for $u\in M_D$ and $v\in X$. 
Since $<I^1(v,z)w:v\in X, w\in W^2>=U$, 
$Y^U(u,z)$ is uniquely determined by $Y^2(u,z)$ and so 
$W^1\times W^2=W^3$. 
\prend

When we want to prove the condition $(4)$ in Hypotheses I, 
the notion of induced VOA will be 
very useful as we mentioned in the introduction. 
Set $S=<\al,\be>$.  
We assume that a pair $(D,S)$ satisfies the conditions (1) and (2) 
in Hypotheses I for a while. 
An important tool is Theorem 3.1. Namely, 
for any irreducible $M_D$-modules $V^{\al'}$ and 
$V^{\be'}$ with $\tilde{\tau}(V^{\al'})=\al'$ 
and $\tilde{\tau}(V^{\be'})=\be'$, 
an $M_D$-module $V^{\al'}\times V^{\be'}$ is irreducible 
for $\al',\be'\in S$. 
Set $V^{\al+\be}=V^{\al}\times V^{\be}$. 
Then by the property of fusion rules, we have 
$V^{\ga+\de}=V^{\ga}\times V^{\de}$ for $\ga,\de\in S$.  
This implies that there is a unique nonzero 
intertwining operator of type $\pmatrix{V^{\ga+\de} \cr
V^{\ga}\quad V^{\de}}$ up to scalar multiple for $\de,\ga\in S$. 
So if we have an algebraic structure (like a VOA) on    
$$(M_D\oplus V^{\al}\oplus V^{\be}\oplus V^{\al+\be}, Y), $$
then $Y$ is uniquely determined up to an $M_D$-isomorphism. 
 
The purpose of this section is to show the following theorem, which is 
a ${\R}$-version of Theorem 6.5 in \cite{Mi5}:

\begin{thm} 
Set $S=<\al,\be>$ and assume the pair $(D,S)$ satisfies the conditions 
(1) and (2) of Hypotheses I. 
Let $F$ be an even linear code containing $D$ such that 
$\langle F,S\rangle=0$. 
Assume that $W=M_D\oplus W^{\al}\oplus W^{\be}\oplus W^{\al+\be}$ 
has a simple VOA structure. Then 
$$V=M_F\oplus {\Ind}_{M_D}^{M_F}(W^{\al})
\oplus {\Ind}_{M_D}^{M_F}(W^{\be})\oplus {\Ind}_{M_D}^{M_F}(W^{\be})$$ 
also has a simple VOA structure up to an $M_F$-isomorphism. 
\end{thm}  

In order to prove the above theorem, we will prove the following lemma:

\begin{lmm}  ${\Ind}_D^F(W^{\al})\times {\Ind}_D^F(W^{\be})
={\Ind}_D^F(W^{\al+\be}).$ 
\end{lmm}

\pr
Let $U$ be an irreducible $M_F$-module such that 
$I\pmatrix{U\cr {\Ind}_D^F(W^{\al}) \quad {\Ind}_D^F(W^{\be})}\not=0$.   
Since  
${\Ind}_D^F(W^{\al})$ and ${\Ind}_D^F(W^{\be})$ are 
irreducible modules satisfying Hypotheses I,  
${\Ind}_D^F(W^{\al})\times {\Ind}_D^F(W^{\be})$ 
is irreducible. 
By Theorem 11.9 in \cite{DL},  we have an injective map 
$$\phi: I\pmatrix{U 
\cr {\Ind}(W^{\al}) \quad {\Ind}(W^{\be})}\to I\pmatrix{U\cr W^{\al}
\quad W^{\be}}$$ 
and so $U$ contains $W^{\al}\times W^{\be}=W^{\al+\be}$, which 
implies $U={\Ind}_D^F(W^{\al+\be})$. 
\prend

\noindent
{\bf Proof of Theorem 3.7.} \qquad 
Let $Y^W(v,z)\in {\End}(W)[[z,z^{-1}]]$ be the given 
vertex operator of $v\in U$ and let 
$$J^{\al',\be'}(v,z)\in I\pmatrix{W^{\al'+\be'}\cr W^{\al'}\quad 
W^{\be'}}$$ be the restriction of $Y(v,z)$ for $v\in W^{\al'}$ and 
$\al',\be'\in S=<\al,\be>$. 
Since Theorem 11.9 in \cite{DL} implies that 
$\phi:I\pmatrix{{\Ind}(W^{\ga'})\cr {\Ind}(W^{\al'})\quad {\Ind}(W^{\be'})} 
\to I\pmatrix{W^{\ga'}\cr W^{\al'}\quad W^{\be'}}$
is injective and the multiplicity of $W^{\al'}\times W^{\be'}$ in 
${\Ind}(W^{\al'+\be'}$ is one, we can choose 
$$I^{\al',\be'}(\ast,z)
\in I\pmatrix{{\Ind}(W^{\al'+\be'})\cr {\Ind}(W^{\al'})\quad {\Ind}(W^{\be'})}$$ 
such that $I^{\al',\be'}(v,z)u=J^{\al',\be'}(v,z)u$ for $v\in W^{\al'}$ 
and $u\in W^{\be'}$.  
Define $Y(v,z)\in {\End}(V)[[z,z^{-1}]]$ by 
$I(v,z)u=I^{\al',\be'}(v,z)u$ for $v\in {\Ind}_D^F(W^{\al'})$ 
and $u\in {\Ind}_D^F(W^{\be'})$. Note $Y(v,z)u=Y^U(v,z)u$ for $u,v\in U$. 
Moreover, the powers of $z$ in $Y(v,z)$ are all integers since 
$\langle \tilde{\tau}({\Ind}(U)), F\rangle=0$.  For $u,v \in U$,  
$Y(u,z)$ and $Y(v,z)$ satisfy the commutativity on $U$.  
For $v\in V$, $Y(v,z)|_{{\Ind}(W^{\al})}$ is at least an intertwining operator 
and  
$Y(v,z)$ satisfies the commutativity with a vertex operator 
$Y(u,z)$ of $u\in M_F$. 
By this commutativity, 
$$\{ w\in {\Ind}(U): I(u',z)I(u,x)w\sim I(u,x)I(u',z)w \}\eqno{(3.9)}$$ 
is a $M_F$-module for $u,u'\in U$. 
Since it contains $U$, it coincides with $V$. 
Namely, $\{Y(u,z):u\in U\cap M_D\}$ satisfies 
the mutual commutativity on $V$.  
Clearly, $\{ Y(v,z): v\in M_D\cup U \}$ 
generates all intertwining operators by the normal products and so 
all $I(v,z)$ for $v\in V$ satisfy the mutually commutativity by Dong's 
lemma.  The other required conditions are also easy to check and 
so we have a VOA structure on ${\Ind}(U)$.  
\prend

\begin{lmm}
Let $V=\oplus_{\al\in S} V^{\al}$ be a VOA satisfying Hypotheses I 
and $W$ be an irreducible $V$-module. Then there is a word $\ga$ and 
irreducible $M_D$-modules with $\tilde{\tau}(W^{\be})=\be$ for $\be\in S+\ga$  
such that $W=\oplus_{\be\in S+\ga}W^{\be}$. 
\end{lmm}

\pr Since $T$ is rational, $W$ is a direct sum of irreducible $T$-modules 
and so we have $W=\oplus_{\be\in S'}W^{\be}$ 
for some $S'$, where $W^{\be}$ is the sum of all irreducible $T$-submodules 
$X$ with $\tilde{\tau}(X)=\be$.  By (1.1), $W^{\be}$ is 
an $M_D$-module. By the similar arguments as in the proof of Theorem 3.2, 
$W^{\be}$ is irreducible. 
Since $\tilde{\tau}(V^{\al}\times W^{\be})=\al+\be$, 
$S'=S+\ga$ for some $\ga$.  
\prend

\begin{lmm}
Let $V=\oplus_{\al\in S} V^{\al}$ be a VOA satisfying Hypotheses I 
and $W=\oplus_{\be\in S+\ga}W^{\be}$ be an irreducible 
module. Assume that $<S+{\Z}_2\ga, D>$ satisfies Hypotheses I. 
Then $W$ is uniquely determined by a $W^{\be}$ for some $\be$. 
\end{lmm}

\pr
Since $V^{\al}\times W^{\be}=W^{\al+\be}$, 
$M_D$-module structure on $W$ is uniquely determined by $W^{\be}$. 
By the similar arguments as in the proof of Theorem 3.4, 
we have the desired assertion. 
\prend

Since the fusion rules (1.1) are all well-defined over ${\R}$ (even 
over ${\Q}$), 
we can rewrite Theorem 4.1 in \cite{Mi5} into the following 
theorem.

\begin{thm}
Under the assumptions (1)$\sim $(4) of Hypotheses I, we obtain 
a fusion product $V^{\al}\times V^{\be}=V^{\al+\be}$ for $\al,\be\in S$. 
Moreover, there is a simple VOA structure on 
$$  V=\boplus_{\al\in S}V^{\al}   $$
such that it contains 
$M_D$ as a sub VOA $V^{(0^n)}$ and 
has a positive definite invariant bilinear form. 
A simple VOA structure on $V$ with a positive definite 
invariant bilinear form is uniquely 
determined up to $M_D$-isomorphisms. 
\end{thm}

\pr
First, we fix module vertex operators $Y^V(v,z)$ for $v\in M_D$. 
Let $Y^{\al,\be}$ be the vertex operator of the VOA 
$V^{\al,\be}=M_D\oplus V^{\al}\oplus V^{\be}\oplus V^{\al+\be}$ such that 
$Y^{0,\be}(v,z)u=Y^V(v,z)u$ for $v\in M_D$ and $u\in V^{\al,\be}$. 
Since $V^{\al}\times V^{\al}=M_D$, there are two possible 
simple VOA structures on $M_D\oplus V^{\al}$.  
Moreover, since we assumed that 
$M_D\oplus V^{\al}$ has a positive definite invariant bilinear form, 
there is a unique VOA structure on $M_D\oplus V^{\al}$ up to 
$M_D$-isomorphisms. Namely, 
if we fix an orthonormal basis 
$\{u^{\al}_i:i\in I_{\al}\}$ of $V^{\al}$. 
then $Y^{\al,\be}(u,z)v$ for $u,v\in V^{\al}$ does not depend on the 
choice of $\be$. 
Define a nonzero intertwining operator 
$$I^{\al,\be}(\ast,z)\in I\pmatrix{V^{\al+\be}\cr V^{\al}\quad V^{\be}}$$ 
for $\al,\be\in S$ by  
$I^{\al,\be}(v,z)u=Y^{\al,\be}(v,z)u$ for $v\in V^{\al}$, $u\in V^{\be}$.
 
Our next step is to choose suitable scalars $\la^{\al,\be}$ 
and define a new vertex operator $Y(v,z)\in {\End}(V)[[z,z^{-1}]]$ by 
$$Y(v,z)u=\la^{\al,\be}I^{\al,\be}(v,z)u  \eqno{(3.a)}$$ 
for $v\in V^{\al}$ and $u\in V^{\be}$ 
such that $\{Y(v,z):v\in V\}$ satisfies the mutual commutativity.
Since intertwining operators satisfy the $L(-1)$-derivative property 
and the other conditions except the mutual commutativity, 
$(V,Y)$ becomes a simple VOA with a positive definite invariant bilinear 
form. 

Set $\dim S=t$ and let $\{\al_1,...,\al_t\}$ be a basis of $S$. 
Set $S_i=<\!\al_1,...,\al_i\!>$ for $i=0,1,...,t$ and  
$V^{i}=\oplus_{\al\in S_i}V^{\al}$. 
We will choose $\la^{\al,\be}$ inductively. 
Since $V^{\al}$ are all $M_D$-modules, 
the module vertex operators $Y^V(v,z)$ of $v\in V^0=M_D$ on $V$ 
satisfy the mutual commutativity if we choose $\la^{0,\al}=1$.
We next assume that there is an integer $r$ such that 
the vertex operators $\{Y(v,z):v\in V^r\}$ satisfy 
the mutual commutativity by choosing $\la^{\al,\be}$ for $\al\in S_r$.   
In particular, $V^r$ is a sub VOA 
and $V$ is a $V^r$-module by these vertex operators. 
It is 
clear that $V^{S^r+\delta}=\oplus_{\ga\in S^r}V^{\delta+\ga}$ are 
irreducible $V^r$-modules 
for any $\delta\in S$ by the fusion rules and $V$ 
decomposes into the direct sum of irreducible $V^r$-modules. 
By the fusion rule of $M_D$-modules $V^{\be}$ and Lemma 3.12, 
we obtain a fusion rule:
$$ V^{\delta+S^r}\times V^{\ga+S^r}=V^{\delta+\ga+S^r}  $$ 
as $V^r$-modules. 

Decompose $V^{r+1}=V^r\oplus V^{\al_{r+1}+S^r}$ as $V^r$-modules.
To simplify the notation, set $\al=\al_{r+1}$ 
and $I^{\be}(v,z)=I^{\al,\be}(v,z)$ for a while. 
Let $\{\ga^i\in S:i\in J\}$ be a set of representatives of cosets 
$S/S^{r+1}$.  
Since there is an injection 
$$ \pi: I\pmatrix{V^{S^r+\al+\ga^i} \cr V^{S^r+\al} \quad V^{S^r+\ga^i}} 
\to I\pmatrix{V^{S^r+\al+\ga^i} \cr V^{\al} \quad V^{\ga^i}} $$
and 
$\dim I\pmatrix{V^{S^r+\al+\ga^i} \cr V^{\al} \quad V^{\ga^i}}=1$,
we can choose a nonzero intertwining operator \\
$I^{\al+S^r,\ga^i+S^r}(\ast,z)\in 
I\pmatrix{V^{S^r+\al+\ga^i} \cr
V^{S^r+\al}\quad V^{S^r+\ga^i}}$
such that 
$I^{\al+S^r,\ga^i+S^r}(v,z)u=Y^{\al,\ga^i}(v,z)u$ for 
$v\in V^{\al}, u\in V^{\ga^i}$. 
Restricting $I^{\al+S^r,\ga^i+S^r}(\ast,z)$ into $V^{\al+\be,\ga^i+\de}$ 
for $\be,\de\in S^r$, we have a scalar $\la_{\al+\be,\ga^i+\de}$ such that 
$I^{\al+S^r,\ga^i+S^r}(v,z)u
=\la_{\al+\be,\ga^i+\de}Y^{\al+\be,\ga^i+\de}(v,z)u$ 
for $v\in V^{\al+\be}$ and $u\in V^{\ga^i+\de}$. 
We will show that $V^{r+1}$ is a sub VOA and $V$ is a $V^{r+1}$-module 
by the above intertwining operators $I^{\al+S^r,\ga^i+S^r}(\ast,z)$.

Set 
$$Q=\{w\in V| Y(u,z)Y(u',x)w\sim Y(u',x)Y(u,z)w 
\mbox{ for }u,u'\in V^{\al}\}.$$
Since $Y(\ast,z)$ is an intertwining operator of $V^r$-modules, 
$Q$ is a $V^r$-module.  On the other hand, by the choice of $Y$, 
$Q$ contains $V^{\ga^i}$ for all $i$. Hence, $Q$ coincides with $V$.  
In particular, all vertex operators in 
$\{Y(u,z): u\in V^r\cup V^{\al}\}$ satisfy the mutual 
commutativity.  Since $V^{r+1}$ is generated by $V^r$ and $V^{\al}$, 
we have 
the desired result. 
This completes the construction of our VOA. 

We next show that the VOA structures on $V$ is unique. 
Assume that there are two VOA structures $(V,Y)$ and $(V,Y')$ on $V$.  
We may assume that all $V^{\al,\be}$ are sub VOAs of $(V,Y)$.  
Since $\dim I\pmatrix{V^{\al+\be} \cr V^{\al}\quad V^{\be}}=1$, 
there are $\la_{\al,\be}$ such that 
$Y'(v,z)u=\la_{\al,\be}Y(v,z)u$ for $v\in V^{\al}$, $u\in V^{\be}$.  
We may assume $\la_{0^n,\be}=1$ and so $\la_{\be,0^n}=1$ by the skew 
symmetry. We will show that by changing the sign of an orthonormal 
basis of 
$V^{S^r +\al_{r+1} }$ if necessary we can get $Y=Y'$.  
Define $f=(-1)^{\be}\in {\End}(V)$ by $(-1)^{<\be,\al>}$ on $V^{\al}$.  
It is clear that $f$ is an automorphism.  Assume 
$Y|_{V^{S^r}}=Y'|_{V^{S^r}}$ and $Y|_{V^{S^{r+1}}}\not= Y'|_{V^{S^{r+1}}}$.  
We assert $\la_{\al}\in \{\pm 1\}$. By changing the sign of 
an orthogonal basis of $V^{\al+S^r}$, we may assume 
$Y'|_{V^{\al}}=Y|_{V^{\al}}$.  By Lemma 3.12, 
$(V^{\al+S^r},Y)$ is isomorphic to $(V^{\al+S^r},Y')$ as a $V^{S^r}$-module.  
On the other hand, since $V^{\al+S^r}\times V^{\al+S^r}=V^{S^r}$, 
$Y'|_{V^{\al+S^r}}=\la_{\al}Y|_{V^{\al+S^r}}$. Hence, $Y=Y'$ on $S^{r+1}$. 
Namely, $\la_{\al}=-1$.  Let 
$\be\in <S^r>^{\perp}$ and $\langle \be,\al\rangle=1$. 
Then by using an automorphism $(-1)^{\be}$, we may assume $Y=Y'$ on $S^{r+1}$. 
By induction, we have $Y=Y'$ on $V$.  
\prend

We will next show a relation between automorphisms of $M_D$ and 
fusion product modules $M_{\al+D}\times W$.  For a word $\al$, 
we can define an automorphism $\sigma_{\al}$ 
of $M_D$ by 
$$ \sigma_{\al}: (-1)^{\langle \be,\al\rangle} \mbox{ on } M_{\be},  $$
which coincides with $\prod_{i\in Supp(\al)}\sigma_{e^i}$, where 
$\sigma_{e^i}$ is an automorphism given in \cite{Mi1} of type 2. 

\begin{lmm}  Suppose $\be=\tilde{\tau}(W)$ and $D_{\be}$ contains a 
maximal self orthogonal subcode $H$ which is doubly even and 
is orthogonal to $\al$, 
then $\sigma_{\al}W$ is isomorphic to $W$. 
\end{lmm}

\pr   
Decompose $M_D$ into $M_D^{+}\oplus M_D^{-}$, where 
$M_D^{\pm}=\{v\in M_D: \sigma_{\al}(v)=\pm v\}$.
Set $E=\{\be\in D:\langle \be,\al\rangle=0\}$. Clearly, 
$M_D^+=M_E$. 
Since $E$ contains a maximal self-orthogonal subcode $H$ of $D_{\be}$  
which is doubly even,  
there is an $M_E$-module $U$ such that 
${\Ind}_{M_E}^{M_D}(U)=W$ by Proposition 3.2. 
It follows from the definition of the induced modules that 
${\Ind}_{M_E}^{M_D}(U)\cong U\oplus (M_D^-\times U)$ as $M_E$-modules.  
The actions of $M_D^-$ switch $U$ and $M_D^-\times U$, that is, 
$u_n(U)\subseteq M_D^-\times U$ and $u_n(M_D^-\times U)\subseteq U$ 
for any $n\in \Z$ and $u\in M_D^-$. Moreover, 
$u_n \sigma_{\al}v=-u_nv$ for $u\in M_D^-$ and $v\in {\Ind}_E^D(U)$.  
It is easy to check that $(1_U,-1_{M_D^-\times U})$ on 
$U\oplus M_D^-\times U$ 
is an isomorphism from $\sigma_{\al}({\Ind}_E^D(U))$ to ${\Ind}_E^D(U)$.  
\prend

For an irreducible $M_D$-module $W$, $\sigma_{\al}W$ is also an irreducible 
$M_D$-module.  
Clearly, $W$ and $\sigma_{\al}W$ are isomorphic as $T$-modules and 
$\sigma_{\al}=\sigma_{\be}$ if and only if $\al+\be\in D^{\perp}$.  
We next investigate an irreducible 
$M_D$-module $M_{D+\al}\times W$ for $\al$ 
satisfying $Supp(\al)\subseteq Supp(\tilde{\tau}(W))$. 
In this case, $M_{D+\al}\times W$ is isomorphic to $W$ as a $T$-module.  
The following lemma is important. 

\begin{lmm}  Let $W$ be an irreducible $M_D$-module and assume $Supp(\al)
\subseteq Supp(\tilde{\tau}(W))$.  
Then $M_{D+\al}\times W$ is isomorphic to 
$\sigma_{\al}W$ as an $M_D$-module. 
\end{lmm}

\pr
Set $U=M_{\al+D}$ and $\be=\tilde{\tau}(W)$. Clearly, 
$\tilde{\tau}(M_{D+\al}\times W)=\tilde{\tau}(\sigma_{\al}W)=\be$. 
By Theorem 3.3, $W'=U\times W$ is irreducible. 
Let $H_{\be}$ be a maximal self-orthogonal 
(doubly even) subcode of $D_{\be}$.  Since an $M_D$-module $W$ 
with $\tilde{\tau}(W)=\be$ is 
uniquely determined by an $M_{H_{\be}}$-submodule, we may assume 
that $D$ is a self-orthogonal doubly even code and 
$Supp(D)\subseteq Supp(\be)$. 
In particular, we may also assume that $W$ and $W'$ are both isomorphic to 
$L(\hf,\st)^{\ots n}$ as $T$-modules.
Since 
$1\leq \dim I_{M_D}\pmatrix{W' \cr U,\quad W}\leq 
\dim I_{T}\pmatrix{L(\hf,\st)^{\ots n} \cr 
M_{\ga}\quad L(\hf,\st)^{\ots n}}=1$, 
an intertwining operator of type 
$\pmatrix{W'\cr M_{\ga+D},\quad W}$ is \\
uniquely determined up to scalar multiples for $\ga\in D+\al$. 
As shown in \S 2.3 or in \cite{Mi5}, we can choose 
a nonzero intertwining operator 
$I(\ast,z)\in I_T\pmatrix{L(\hf,\st)^{\ots n} \cr
M_{\ga}\quad L(\hf,\st)^{\ots n}}$ by 
$$  I(q^{\ga},z)=I(\hat{q}^{\ga},z)=\ots I^{g_i,\st}(q^{g_i},z),   $$
where $I^{g_i,\st}(\ast,z)$ is a fixed intertwining operator 
of type $\pmatrix{L(\hf,\st)\cr L(\hf,{g_i\over 2})\quad L(\hf,\st)}$, 
see \S 2.3. 

By Theorem 3.2, there are linear representations $\chi$ and $\phi$ of 
$\hat{D}$ such that \\
$W\cong L(\hf,\st)^{\otimes n}\otimes Q_{\chi}$ and 
$W'\cong L(\hf,\st)^{\otimes n}\otimes Q_{\phi}$. 
By the associativity property of intertwining operators, 
$$\begin{array}{l}
I(v_nq^{\al},z)=Res_{x}\{(x-z)^nY^{W'}(q^{\be},x)I(q^{\al},z)
 -(-z+x)^nI(q^{\al},z)Y^W(q^{\be},x)\} \cr
=Res_{x}\{(x-z)^nI^{\otimes n}(\hat{q}^{\be},x)\phi(e^{\be})I(q^{\al},z)
 -(-z+x)^nI(q^{\al},z)I^{\otimes n}(\hat{q}^{\be},x)\chi(e^{\be})\} 
\end{array} $$
for $q^{\be}\in M_{\be}\subseteq M_D$ and $u\in M_{\al}$. In particular, 
for a sufficiently large $N$, we obtain 
$$0=Res_{x}\{(x-z)^NI^{\otimes n}(\hat{q}^{\be},x)\phi(e^{\be})I(q^{\al},z)
-(-z+x)^NI(q^{\al},z)I^{\otimes n}(\hat{q}^{\be},x)\chi(e^{\be})\}. $$
On the other hand, as we showed in \S 2.3, $I(\ast,z)$ 
satisfies the super-commutativity:
$$(x-z)^NI^{\otimes n}(\hat{q}^{\be},x)I^{\otimes n}(\hat{q}^{\al},z)
-(-1)^{\langle \al,\be\rangle}(-z+x)^NI^{\otimes n}(\hat{q}^{\al},z)
I^{\otimes n}(\hat{q}^{\ga},x)=0.$$ 
Therefore, 
$$Res_{x}\{(x-z)^N\phi(e^{\be})-(-1)^{\langle \al,\be\rangle}(-z+x)^N
\chi(e^{\be})\}=0 $$
and so $\phi(e^{\be})=(-1)^{\langle \al,\be\rangle}\chi(e^{\be})$ 
for $\be\in D$. Hence, $W'$ is isomorphic to $\sigma_{\al}W$ as 
$M_D$-module. 
\prend

\begin{rmk}  The above lemma may 
look a little strange since we 
usually obtain relations 
$\sigma(W^1)\times \sigma(W^2)=\sigma(W^1\times W^2)$
and 
$(M_{\al+D}\times W^1)\times (M_{\al+D}\times W^2)=(W^1\times W^2)$ 
for an automorphism $\sigma$ and a coset module $M_{\al+D}$, 
respectively.  However, if we have $\sigma(W^i)\cong M_{D+\al}\times W^i$ 
for $i=1,2$, then 
$W^1\times W^2$ does not satisfy the condition of the 
above lemma by (1.1) and so $\sigma(W^1\times W^2)=W^1\otimes W^2$.  
\end{rmk}

\section{Positive definite invariant bilinear form}
In order to construct $V^{\natural}$, we will use 
"induced VOAs". So, we will prove the following theorem. 

\begin{thm}  
Assume that $W^{\al}$ is an irreducible $M_D$-module with 
$\tilde{\tau}(W^{\al})=\al$ and $(D,<\al>)$ satisfies the 
conditions (1) and (2) of Hypotheses I. 
Let $F$ be an even linear code containing $D$ satisfying 
$\langle F,\al\rangle=0$. 
If a simple VOA $U=M_D\oplus W^{\al}$ has a positive definite 
invariant bilinear form, then
${\Ind}_{M_D}^{M_F}(U)\ (=M_F\oplus {\Ind}_{M_D}^{M_F}(W^{\al}))$ also has 
a positive definite invariant bilinear form. 
\end{thm}  

\pr
We note that if an irreducible $M_D$-module $W$ is not isomorphic to $M_D$, 
then the lowest degree of $W$ is greater than $0$.  
Clearly, it is sufficient to prove the lemma for $F=<\al, (1^n)>^{\perp}$.  
Since $<\al,(1^n)>^{\perp}$ is generated by words of weight 2, 
we may assume $F=D+{\Z}_2\be$ such that $|\be|=2$ by induction.  
Say $\be=(110^{n-2})$.  
Since $|\be|=2$ and $\langle \be,\al\rangle=0$, $Supp(\be)\subseteq Supp(\al)$ 
or $Supp(\be)\cap Supp(\al)=\emptyset$. 

Since $D_{\al}$ contains a direct sum $E$ of Hamming codes such that 
$Supp(E)=Supp(\al)$, ${\Ind}_D^F(W^{\al})$ is irreducible. 
Set 
$$ V=M_F\oplus {\Ind}_D^F(W^{\al}). $$
By Theorem 3.7, $V$ has a VOA structure. 
Since ${\Ind}_D^F(W^{\al})\times {\Ind}_D^F(W^{\al})=M_F$ by Lemma 3.10, 
there are only two possible VOA structures on $V$. 
Namely, if one is $(M_F\oplus {\Ind}_D^F(W^{\al}), Y)$, then the other 
is $(M_F\oplus \sqrt{-1}{\Ind}_D^F(W^{\al}), Y)$.  
Since $W^{\al}\times W^{\al}=M_D$, we may assume 
$(M_F\oplus {\Ind}_D^F(W^{\al}), Y)$ contains U as a sub VOA.   
Let $E$ be a maximal self orthogonal doubly even subcode of $D_{\al}$.
Then $W^{\al}$ is a direct sum $\oplus W^i$ of distinct 
irreducible $M_E$-modules $W^i$ and 
$V^i=M_E\oplus M_{E+\be}\oplus W^i\oplus (M_{E+\be}\times W^i)$ is a sub VOA 
of $V$ for each $i$.  
Since $(M_E\oplus W^i, Y)$ is a sub VOA of $M_D\oplus W^{\al}$, 
$(M_E\oplus W^i, Y)$ has a positive definite invariant bilinear form. 

We will later show that a VOA structure $(V^i,Y)$ 
on $V^i$ has a positive definite invariant 
bilinear form. 
In particular, $W^i\oplus (M_{E+\be}\times W^i)$ has an orthonormal 
basis with respect to $Y$. 
Then since $M_{D+\be}\times W^{\al}$ coincides with 
$\oplus (M_{E+\be}\times W^i)$, we have the desired result. 
Therefore, we may assume that $Supp(D)=Supp(\al)$ and 
$D$ is a direct sum $D=E^1\oplus\cdots \oplus E^s$ 
of Hamming codes $E^i$ by Hypotheses I. In particular, 
$W^{\al}$ is irreducible as a $T$-module.  
Since a VOA structure $(V,Y)$ on $V$ containing $U$ is uniquely 
determined, it is sufficient to show that there exists a VOA structure on 
$(V,Y)$ with a positive definite invariant bilinear form. 
For if $(V,Y')$ is the other simple VOA structure on $V$, then 
$(W^{\al},Y')$ has a negative definite invariant bilinear form and so 
$(V,Y')$ does not contain $U$. 
If $Supp(\be)\cap Supp(\al)=\emptyset$, then $<D,\be>$ is self-orthogonal. 
Let $D^0$ be the code of length $n-2$ consisting of the codewords $\ga$ 
such that $(00\ga)\in D$.  Then 
$M_D=L(\hf,0)\otimes L(\hf,0)\otimes M_{D^0}$ and 
$M_{D+\be}=L(\hf,\hf)\otimes L(\hf,\hf)\otimes M_{D^0}$.  
By the above decompositions, we can write 
$$W^{\al}\cong L(\hf,h^1)\otimes L(\hf,h^2)\otimes W'$$ 
and 
$$M_{D+\be}\times W^{\al}\cong 
L(\hf,h^1+\hf)\otimes L(\hf,h^2+\hf)\otimes W'. $$
for some irreducible $M_{D^0}$-module $W'$ and $h^1,h^2=0,\hf$ and 
$hi+\hf=0$ if $h^i=\hf$ and $h^i+\hf=\hf$ if $h^i=0$. 
Since $L(\hf,0)^{\otimes 2}\oplus L(\hf,\hf)^{\otimes 2}\cong 
\tilde{V}_{2\Z x}=(V_{2\Z x})^{\theta}\oplus \sqrt{-1}(V_{2\Z x})^-$ 
for $\langle x,x\rangle=1$, $\sqrt{-1}x(0)$ 
is an isomorphism from $L(\hf,h^1)\otimes L(\hf,h^2)$ to 
$L(\hf,h^1+\hf)\otimes 
L(\hf,h^2+\hf)$ and $x(0)^2$ acts diagonally on 
$L(\hf,h^1)\otimes L(\hf,h^2)$ 
with a positive eigenvalues.  
Let $\{v^i:i\in I\}$ be an orthogonal normal basis such that each $v^i$ is in 
an eigenspaces of $x(0)^2$.  Then $\{\sqrt{-1}x(0)v^i:i\in I\}$ is a basis of 
$L(\hf,h^1+\hf)\otimes L(\hf,h^2+\hf)$ and 
$$\langle \sqrt{-1}x(0)v^i, \sqrt{-1}x(0)v^j\rangle=
\langle v^i, x(0)^2v^j\rangle=\delta_{ij}\langle v^i,x(0)^2v^j \rangle 
\geq 0.$$ 
Hence, ${\Ind}_D^F(U)$ has a positive definite invariant bilinear form. 

We next assume $Supp(\be)\subseteq Supp(\al)$. Since $D$ is a direct 
sum of Hamming codes and the weight of $\be$ is $2$, it is sufficient 
to treat the following two cases: \\
(1) $Supp(\be)\subseteq Supp(E^1)$.   \\
(2) $D=E_8\oplus\cdots \oplus E_8$ and $\be=(10^710^70^{n-16})$. \\

Case (1).  By Lemma 3.8, there is another set of coordinate conformal 
vectors $\{d^i\}$ of $M_D$ such that $W$ is a coset module $M_{D+\ga}$ 
w.r.t. $<d^i>$.  
Since $Supp(\be)\subseteq Supp(E^1)$ and $\be$ has an even weight, 
$M_{\be+D}$ is also a coset module $M_{\delta+D}$. 
Namely, ${\Ind}_D^F(U)$ is a code VOA $M_{<D,\delta,\ga>}$ w.r.t. $<d^i>$.  
Hence, it has a positive definite invariant bilinear form. \\

Case (2).  By taking another set of coordinate conformal vectors, 
we may assume that $\al=(1^{16}0^{n-16})$ and $\be=(10^710^70^{n-16})$. 
Since $L(\hf,\hf)\otimes L(\hf,\hf)$ has a positive definite invariant 
bilinear form and the lowest weight is an integer, we may also assume that 
$n=16$ and $\al=(1^{16})$.  We will find such a VOA in $V_{E_8}$ in 
the next section. This will complete the proof of Theorem.   
\prend

\section{$E_8$-lattice VOA}
As we mentioned in the introduction, we will gather the 
parts of $V^{\natural}$ from $\tilde{V}_{E_8}$. Hence, 
the main aim of this section is to study the structure of $V_{E_8}$ and 
$\tilde{V}_{E_8}$. 
In particular, we will show that $\tilde{V}_{E_8}$ satisfies 
the conditions (1)$\sim$(5) of Hypotheses I. Incidentally, 
we will see that the orbifold construction from VOA $V_{E_8}$ 
coincides with the change of a set of coordinate 
conformal vectors of a Hamming code sub VOA of $V_{E_8}$.  

Let $E_8$ denote the root lattice of type $E_8$.  It is known that $E_8$ 
is the unique positive definite unimodular even lattice of rank 8. 
We first define lattices $E_8(m):m=1,2,3,4,5$ and $L(1)$. 
Let $\{x^1,...,x^8\}$ be an orthonormal basis and 
set $E_8(1)=< \hf(\sum_{i=1}^8x^i), x^i \pm x^j: i,j=1,...,8 >$ 
and $L(1)=<x^i : i=1,...,8>$, where $<u^i:i\in I>$ denotes 
a lattice generated by $\{u^i:i\in I\}$.
It is easy to check that $E_8(1)$ is isomorphic to $E_8$. 
We can define the other $E_8$-lattices as follows:
$$\begin{array}{ll}
E_8(2)=&< \hf(x^1-x^2-x^3-x^4)\!+\!x^5, \hf(x^5+x^6+x^7+x^8)\!+\!x^1, \cr
&x^i\pm x^j :i,j\in\{1,2,3,4\}, \mbox{ or } i,j\in\{5,6,7,8\}>. \cr
E_8(3) =&< \hf(x^1-x^2-x^5-x^6)+x^3, \hf(x^1+x^2-x^3-x^4)-x^7, \cr
&\hf(-x^5-x^6+x^7+x^8)+x^1, x^1+x^3+x^5+x^7, 
x^{2i-1}+x^{2i},\ (i=1,2,3,4)>  \cr
E_8(4) =&<\hf(x^1-x^3-x^5-x^7)+x^2, \hf(x^1-x^2+x^5-x^6)-x^3, \cr
&\hf(-x^1+x^2-x^3-x^4)-x^7, \hf(x^1+x^3-x^6+x^8)+x^5, 2x^1,...,2x^8> 
\end{array} \eqno{(5.1)}$$
Fix $m=1,2,3,4$ and set $L=E_8(m)$.  
Let $V_L$ be a lattice VOA constructed as in \cite{FLM2} and 
$\theta$ an automorphism of $V_L$ induced from $-1$ on $L$. 
Since $E_8(m)$ contains $\{2x^1,...,2x^8\}$,   
we obtain a set $I=\{e^i:i=1,...,16\}$ 
of 16 mutually orthogonal conformal vectors of $V_L$, where 
$$  e^{2i-j}={1\over 4}x^i(-1)^2{\bf 1}-(-1)^j{1\over 4}
(\iota(2x^i)+\iota(-2x^i)) 
\eqno{(5.2)}$$
for $i=1,...,8, \ j=1,0$ as given in \cite{DMZ}.  
Since they are all in $V_L^{\theta}$, we can also 
take this set as a set of coordinate conformal vectors of 
$\tilde{V}_{E_8}$.
Hence, the decompositions of $V_L$ and $V_{E_8}$ into 
the direct sum of irreducible $T$-submodules are the same, where 
$T=<e^1,...,e^{16}>$, (see the proof of Proposition 2.2).   

Let $P(m)=<\tau_{e^i}:i=1,...,16>$ and 
$L(m)=E_8(m)\cap L(1)$. 
By (2.3), $\tilde{V}_{L(m)}$ contains $<e^1,...,e^{16}>$ and it is 
straightforward to check that 
$(\tilde{V}_L)^{P(m)}$ coincides with $\tilde{V}_{L(m)}$. 
Define a code $D(m)$ of length 16 by  
$$ M_{D(m)}\cong (V_{E_8})^{P(m)}.$$ 
It is also not difficult to check that $(\tilde{V}_L)^{P(m)}$ 
has a decomposition satisfying Hypotheses I with respect to 
$(D(m), D(m)^{\perp})$. 
However, this is not what we want because   
$D(m)$ has a root and so $(M_{D(m)})_1\not=0$ for 
$m=1,2,3,4$.  We are going to get a code $D$ without roots. 
In order to find such a decomposition, we will change 
the set of coordinate conformal vectors.  Incidentally, 
this process coincides with an orbifold construction as we will 
show. 

Let's explain the relation between the orbifold construction and changing 
the coordinate sets of conformal vectors. 
It is known that any orbifold construction from $V_L$ 
is isomorphic to itself. Let's explain the orbifold construction. 
Let $\theta$ be an automorphism of $V_L$ induced from $-1$ on $L$. 
$\theta$ fixes 
$\iota(x^i)+\iota(-x^i)$ and acts as  
$-1$ on $x^i(-1){\bf 1}$ and $\iota(x^i)-\iota(-x^i)$.  
Hence, $\theta$ acts on 
$M_{\al}$ as $(-1)^{\langle \al,(\{01\}^8)\rangle}$ and so 
the fixed point space $M_{D(m)}^{\theta}$ is equal to the 
direct sum $\boplus_{\al\in D(m,+)}M_{\al}$, where 
$D(m,+)=\{\al\in D(m): \langle \al,(\{01\}^8)\rangle=0\}$. 
Assume that the twisted part of the orbifold construction does not contain
any coset modules.  
Suppose that $V=\oplus_{\al\in S}V^{\al}$ is a VOA satisfying 
Hypotheses I such that 
$\tilde{\tau}(V^{\al})=\al$ and $V^{(0^{2n})}\cong M_D$, 
where $D$ is a code of length $2n$ containing 
$(0^{2i}110^{2n-2i-2})$ for all $i=1,...,m$. 
Set $\be=(\{01\}^{n})$. \\
Then the orbifold construction 
is corresponding to the following three steps as we will see in the 
next example. \\
(1)  Take an half $M_{D(+)}$ of $M_D$, where 
$D(+)=\{\al\in D:\langle \al,\be\rangle=0\}$. \\
(2)  Take an $M_{D(+)}$-module $V^{\be}$ with $\tilde{\tau}(V^{\be})=\be$ 
and generate $M_{D(+)}$-modules $V^{\be+\ga}$ with 
$\tilde{\tau}(V^{\be+\ga})=\be+\ga$ by $V^{\be+\ga}=V^{\be}\times V^{\ga}$ 
for $\ga\in S$. \\
(3)  Construct a VOA structure on $\tilde{V}=\oplus_{\al\in <S,\be>}V^{\al}$.\\
 
In the case of $E_8(1)$, $\tilde{\tau}(V_{L(1)+v})=(1^{16})$ 
for $v=\hf(\sum_{i=1}^{16}x^i)$ and so $S(1)=<(1^{16})>$ and 
$D(1)$ is the set of all even words of 
length 16.  
$D(1)$ contains a self dual subcode $H=H^1_8\oplus H^2_8$, where 
$H^i_8$ are Hamming codes and $Supp(H_8^1)=\{1,2,...,8\}$ and 
$Supp(H_8^2)=\{9,...,16\}$. Since $\langle ((10)^{8}),\be\rangle=0$ for 
any $\be\in H$, 
we have $M_H \subseteq V_L^{\theta}$. Therefore, the decompositions of 
$V_L$ and $\tilde{V}_L$ as $M_H$-modules are exactly the same.  
Since $D(1)$ consists of all even words, the center $Z(\widehat{D(1)})$ 
is $<\pm e^{(0^{16})}, \pm e^{(1^{16})}>$ and so 
there are exactly 2 irreducible $M_{D(1)}$-modules 
${\Ind}_{M_H}^{M_{D(1)}}(H(\st,(0^8))\ots H(\st,(0^8)))$ and  
${\Ind}_{M_H}^{M_{D(1)}}(H(\st,(0^8))\ots H(\st,\xi^1))$ 
by Theorem 3.2. The difference between them is judged by the action of 
$q^{(1^{16})}=(q^{\otimes 16})\otimes e^{(1^{16})}$.  
By (2.3) and the proof of Proposition 2.2, we have 
$q^{(1^{16})}=x^1(-1)\cdots x^8(-1){\bf 1}$ and 
$x^i(-1){\bf 1}=\sqrt{-1}(q^{\otimes 2})e^{\xi_{2i-1}}e^{\xi_{2i}}$. 
Since the eigenvalue of $q^{(1^{16})}$ on $\iota(\hf \sum x^i)$ is positive, 
$$V_{E_8}=M_{D(1)}\oplus {\Ind}_{M_H}^{M_{D(1)}}(H(\st,(0^8))\otimes H(\st,(0^8))) 
\eqno(5.4) $$
by the choice of $E(1)$. We should note that 
the difference between the above two modules is given by 
the action of $q^{(1^{16})}=(\otimes_{i=1}^{16}q^1)\otimes e^{(1^{16})}$.  
By Lemma 3.8, 
$M_{H_8}^{\theta}$ contains another set 
of coordinate conformal vectors $\{f^1,...,f^8\}$ 
such that 
$H(\st,(0^8))$ w.r.t. $<e^i>$ is 
isomorphic to $H(\hf,\xi_1)$ w.r.t. $<f^i>$. 
We note that $H(\hf,\al)$ is a coset module $M_{H_8+\al}$. 
Take the set $J=\{f^1,...,f^8,e^9,...,e^{16}\}$ as a new set of 
coordinate conformal vectors. Then 
for $\be\in D(1)$ satisfying $\langle \be,(1^80^8)\rangle=1$, 
the $\tilde{\tau}(M_{H+\al})$ is also a coset module w.r.t. $J$ 
and 
$\tilde{\tau}(H(\st,(0^8))\ots H(\st,(0^8)))$ w.r.t. $J$ is $(0^81^8)$.  \\
Hence, the set $\tilde{\tau}(V_L)$ w.r.t. $J$ is 
$S^2=\{(0^{16}),(1^80^8),(0^81^8),(1^{16})\}$. 
Set $P^2=<\tau_{f^i},\tau_{e^j}:i=1,...,8,j=9,...,16>$ 
and define a linear code $D_2$ by $(V_L)^{P^2}\cong M_{D_2}$ w.r.t. $J$, 
then $D_2$ splits into a direct sum $D_2^1\oplus D_2^2$ such that 
$D_2^1$ and $D_2^2$ are the sets of all 
even words in $\{1,2,...,8\}$ and $\{9,...,16\}$, respectively.  
We note that this process is corresponding to an orthogonal 
transformation 
$$ \hf\pmatrix{
1 &-1&-1&-1 \cr
1 & -1& 1& 1 \cr
1 & 1& -1& 1 \cr
1 & 1& 1& -1 } \eqno{(5.5)} $$
by (2.3).  Therefore, this decomposition coincides with 
the decomposition given by $E_8(2)$. We note that 
$(1^{16})\in D_2$ and 
$M_{(1^{16})+E}$ w.r.t. $<e^i>$ is still equal to $M_{(1^{16})+E}$ w.r.t. $J$. 

We next consider the case of $E_8(2)$ and $S^2=<(1^80^8),(0^81^8)>$. 
We use the above decomposition again by renaming 
$\{f^1,...,f^8,e^9,...,e^{16}\}$, $J$ and $D_2$ by $\{e^1,...,e^{16}\}$, 
$I$ and $D(2)$, respectively.  
Set 
$$ \begin{array}{l}
I_1=\left\{\al\in D(2): Supp(\al)\subset \{1,...,4,9,...,12\} \right\}  \cr
I_2=\left\{\al\in D(2): Supp(\al)\subset \{5,...,8,13,...,16\} \right\}. 
\end{array} $$  It is clear that 
$I_i$ contains Hamming code $H_i$ for $i=1,2$.  Take a new coordinate set 
$\{f^1,...,f^4,f^9,...,f^{12}\}$ of $H_1$ and define a new set 
$$J=\{f^1,...,f^4,e^5,...,e^8,f^9,...,f^{12},e^{13},...,e^{16}\} $$ 
as a set of 
coordinate conformal vectors of $V_L$.  Then 
if an $M_{H_1}\ots M_{H_2}$-module $U$ has a $\tau$-word $(\al,\be)\in 
\{1,...,4,9,...,12\}\oplus \{5,...,8,13,...,16\}$ w.r.t. $I$, then the 
$\tau$-word w.r.t. $J$
is either $(\al,\be)$ or $(\al^c,\be)$.  Moreover, there is a 
submodule with a $\tau$-word $(1^40^41^40^4)$ w.r.t. $J$.  
An example is $M_{H_1\oplus H_2+\al}$, where $\al$ is a word with 
$\langle \al,(1^40^41^40^4)\rangle=1$.  Therefore, we have 
$$D_3=<D_3^1\oplus D_3^2\oplus D_3^3\oplus D_3^4, \{1,5,9,13\}>
\eqno{(5.6)}$$ where 
$D_3^i$ is the set of all even words in $\{4i-3,4i-2,4i-1,4i\}$ for 
$i=1,...,4$.  We also obtain 
$$ S^3=<(1^{16}),(1^80^8),(1^40^41^40^4)>.\eqno{(5.7)} $$  
This corresponds to 
the decomposition with respect to $E_8(3)$ and $D_3=D(3)$. 
$D(3)$ also contains two orthogonal Hamming codes $H_1(3)$ and $H_2(3)$ 
whose supports are 
$$\{1,2,5,6,9,10,13,14\} \quad \mbox{ and }\quad \{3,4,7,8,11,12,15,16\}.$$  
Repeating the above arguments, we have 
$$ S^4=<(1^{16}),(1^80^8),(1^40^41^40^4),(\{1^20^2\}^4)> \eqno{(5.8)}$$
and $D(4)=(S^4)^{\perp}$.  $D(4)$ still contains a direct sum of 
2 Hamming codes whose supports are $(\{10\}^8)$ and $(\{01\}^8)$.  
Repeating the same arguments again, we obtain 
$$ S^5=<(1^{16}), (1^80^8), (1^40^41^40^4), (\{1100\}^4), (\{10\}^8) > 
\eqno{(5.9)}$$
and $D(5)=(S^5)^{\perp}$. \\

\begin{rmk}
Since $D(5)$ does not contains a subcode of rank 8 consisting 
of the form $\{(\al,\al):\al\in {\Z}_2^8\}$ 
for any splits of coordinates into 8 and 8, it is impossible 
to assign $x^i(-1){\bf 1}$ to $L(\hf,\hf)\ots L(\hf,\hf)$ for all 
$i=1,...,8$. Thus, we cannot construct 
$D(5)$ and $S^5$ from a lattice directly. 
\end{rmk}

Let's finish the proof of Theorem 4.1. 
Set $D=D(1)$ and $\be=(1^80^8)$. Set $H=H_8\oplus H_8$ as in (5.5). 
Viewing $V_{E_8}$ as an $M_{H}$-module, $V_{E_8}$ is a direct sum of 
distinct irreducible $M_H$-modules. 
Since $D$ is the set of all even words, 
$M_D$ contains $H(\hf,\xi_1)\ots H(\hf,\xi_1)$ and so 
$V_{E_8}$ has a sub VOA isomorphic to  
$$\begin{array}{rl}
&(H(\hf,(0^8))\ots H(\hf,(0^8)))\oplus (H(\hf,\xi_1)\ots H(\hf,\xi_1)) \cr
\oplus &(H(\st,(0^8))\ots H(\st,(0^8)))\oplus (H(\st,\xi_1)\ots H(\st,\xi_1)),
\end{array} \eqno{(5.10)}$$
where $\xi_1=(10^7)$. This is the desired VOA in the proof of Theorem 4.1. \\

Set $D_{E_8}=D(5)$ and $S_{E_8}=S^5$.  
We will show that this pair $(D_{E_8}, S_{E_8})$ satisfies the 
conditions $(1)$ and $(2)$ of Hypotheses I.   
We note that $D_{E_8}$ is a Reed Muller code $RM(4,3)$ and $S_{E_8}$ is 
a Reed Muller code $RM(4,1)$.  

\begin{lmm}  The pair $(RM(4,3), RM(4,1))$ satisfies the condition 
$(1)$ and $(2)$ of Hypotheses I. 
\end{lmm}

\pr
Set $D=RM(4,3)$ and $S=RM(4,1)$.  
Condition $(1)$ is clear. 
The weight enumerator of $RM(4,1)$ is $x^{16}+30x^8y^8+y^{16}$. 
We note that for any $\be\in RM(4,1)$ with weight $8$, 
$D_{\be}$ and $D_{\be^c}$ are $[8,4,4]$-Hamming codes. 
We always set $H_{\ga}=E_{\ga}=D_{\ga}$ for $\ga\in RM(4,1)$ with weight $8$.
If $\ga=(0^{16})$, then set $H_{\ga}=E_{\ga}=\{(0^{16})\}$. 
Let $\al,\be\in RM(4,1)$.  
We can always choose $H_{\be}$, $H_{\al+\be}$, 
$E_{\al}$ and $E_{\al^c}$ satisfying the conditions $(2),(2.1)$ and $(2.2)$ 
of Hypotheses I.   

Assume $\al=(0^{16})$ or $(1^{16})$, then 
$D_{\be}=D_{\be+\al}$ or $D_{\be}\oplus D_{\be+\al}\subseteq D$. 
In particular, there is a direct sum $H$ of 2 Hamming codes containing 
some maximal self orthogonal subcodes $H_{\be}$ and $H_{\al+\be}$. 
Set $E_{(1^{16})}=H$. 
Clearly, since $E_{\al}+E_{\al^c}=H$ and $H_{\be},H_{\be^c}\subseteq H$, 
they satisfy the condition $(2.3)$ of Hypotheses I. 

We next assume that the weight of $\al$ is $8$. 
If $\be=(0^{16}),(1^{16}), \al$ or $\al^c$, 
then set $H_{(1^{16})}=E_{\al}\oplus E_{\al^c}$. 
Then they satisfy (2.3). 

The remaining case is that $\al, \be, \al+\be$ have weight 8. 
Say $\al=(1^80^8)$ and $\be=(1^40^41^40^4)$. 
We use expressions 
$${\Z}_2^{16}=\{(\delta_1,\delta_2,\delta_3,\delta_4): \delta\in \Z_2^4 \}.$$
Clearly, since $E_{\ga}=H_{\ga}=D_{\ga}$ is a Hamming code for 
$\ga\in S$ with $|\ga|=8$, we have 
$$\begin{array}{rl}
E_{\al}=&\{(\delta\delta 0^40^4),\ (\delta\delta^c 0^40^4)
:\delta\in \Z_2^4 \mbox{ even}\}, \cr
E_{\al^c}=&\{(0^40^4 \delta\delta),\ (0^40^4 \delta\delta^c)
:\delta\in \Z_2^4 \mbox{ even}\}, \cr
H_{\be}=&\{(\delta 0^4\delta 0^4),\ (\delta 0^4\delta^c 0^4):
\delta\in \Z_2^4 \mbox{ even}\}, \cr
H_{\al+\be}=&\{(0^4\delta\delta 0^4),(0^4\delta\delta^c 0^4)
:\delta\in \Z_2^4 \mbox{ even}\}. \end{array} $$
Since 
$$ \begin{array}{l}
(0^4\delta\delta 0^4)-(\delta 0^4\delta 0^4)
=(\delta\delta 0^40^4) \quad \mbox{ and }\cr
(0^4\delta\delta^c 0^4)-(\delta 0^4\delta^c0^4)
+(\delta\delta 0^40^4),
\end{array} $$ 
we obtain $H_{\al+\be}+E_{\al}=H_{\be}+E_{\al}$ and so $(2.3)$. 
\prend 

\begin{prn}  
There are 16 mutually orthogonal conformal vectors $\{e^1,...,e^{16}\}$ 
in $\tilde{V}_{E_8}$ 
such that the decomposition 
$$V_{E_8}=\boplus_{\chi\in S_{E_8}}V_{E_8}^{\chi}$$
given by $\{e^1,...,e^{16}\}$ satisfies Hypotheses I, where \\
(1) \quad the order of $P=<\tau_{e^i}:i=1,...,8>$ is $32$, \\
(2) \quad $D_{E_8}\cong RM(4,1)$, $S_{E_8}=D_{E_8}^{\perp}$, \\
(3) \quad $(V_{E_8})^P=V_{E_8}^{(1^{16})}$ is isomorphic to a code VOA 
$M_{D_{E_8}}$, \\
(4) $\tilde{\tau}(V_{E_8})^{\chi}=\chi$. 
\end{prn}

\pr 
We have already shown that there are 16 mutually orthogonal 
conformal vectors in $\tilde{V}_{E_8}$ satisfying the conditions (1)
$\sim$ (3).  
By Lemma 5.1, $(D_{E_8}, S_{E_8})$ satisfies the conditions (1) and 
(2) of Hypotheses I.  Hence, they satisfy all conditions of Hypotheses I. 
\prend

We next talk about the reverse of the above process. It is clear that 
we can reverse the process.  However, there is another 
important step. Namely, let 
$$  \tilde{V}_{E_8}=\oplus_{\al\in S^n}V^{\al}   $$
be the decomposition such that $V^{(0^{16})}\cong M_{D^n}$.  
Let $\be$ be an even word so that $<\be>^{\perp}\cap S^n\cong S^{n-1}$.  
Set $\tilde{S}^{n-1}=<\be>^{\perp}\cap S^n$ and 
$\tilde{D}^{n-1}=(\tilde{S}^{n-1})^{\perp}$. Then 
$V^+=\oplus_{\al\in \tilde{S}^n}V^{\al}$ is a sub VOA and  
the induced VOA 
$$  \tilde{V}^{n-1}={\Ind}_{D^n}^{\tilde{D}^{n-1}}(V^+)   $$
is also a VOA containing $M_{\tilde{D}^{n-1}}$. 

At the end of this section, we will explain properties of 
the automorphisms 
of a lattice VOA $V_L$ for an even lattice $L$.  
Let $L_2$ denote the set of all elements of $L$ with squared length 4.  
As we showed, for any $a\in L_2$, we can define two conformal vectors \\
$$ \begin{array}{l}
e^+(a)={1\over 16}a(-1)^2{\bf 1}+{1\over 4}(\iota(a)+\iota(-a)) \cr
e^-(a)={1\over 16}a(-1)^2{\bf 1}+{1\over 4}(\iota(a)+\iota(-a)).  
\end{array}$$
Then we have :

\begin{lmm} 
$\tau_{e^+(a)}=\tau_{e^-(a)}$ on $V_L$. By setting $\tau_a=\tau_{e^+(a)}$, 
we obtain $[\tau_a,x(m)]=0$ and 
$$ \tau_{a}: \iota(x) \to (-1)^{\langle x,a\rangle}\iota(x).  $$
In particular, $<\tau_a: a\in L_2>$ is an elementary abelian 2-subgroup of 
${\Aut}(V_L)$.  
If $\langle a,b\rangle$ is odd for $a, b\in L_2$, then 
$\tau_{b}(e^{\pm}(a))=e^{\mp}(a)$.  
\end{lmm} 
 
\pr
Since $\langle a,L\rangle\in \Z$ and 
$\langle a,a\rangle=4$, 
$L\subseteq {1\over 4}{\Z}a\oplus {1\over 4}<a>^{\perp}$. 
In particular, we may view $V_L\subseteq V_{{1\over 4}{\Z}a}\oplus 
V_{{1\over 4}<a>^{\perp}}$.
Recall (2.3) 
$$\tau_{e^i_a}:\left\{ \begin{array}{rcl}
1 &on & a(-1){\1}, \  \iota((\hf+{\Z})a), \ \iota({\Z}a) \cr
-1&on& \iota((\pm {1\over 4}+{\Z})a)
\end{array} \right.  $$
for $i=1,2$.  Then, $[\tau_{e^{\pm}(a)},x(m)]=1$ and 
$$ \tau_{e^{\pm}(a)}: \iota(x) \to (-1)^{\langle x,a\rangle}\iota(x).  $$
Therefore, we obtain the desired results.  
\prend

\begin{thm}  For $g\in {\Aut}(S_{E_8})$, there is an automorphism 
$\tilde{g}$ of $\tilde{V}_{E_8}$ such that $\tilde{g}(e^i)=e^{g(i)}$. 
\end{thm} 

\pr
First we note that $S_{E_8}$ is isomorphic to the Reed Muller code $RM(4,1)$, 
which is defined as follow: \\
Let $F={\Z}_2^4$ be a vector space over ${\Z}_2$ of dimension $4$ 
and denote $(1000)$, $(0100)$, $(0010)$, $(0001)$ by $v^1$, $v^2$, 
$v^3$, $v^4$, respectively. 
Define $\langle (a_i),(b_i)\rangle=\sum a_ib_i$. 
The coordinate set of Reed Muller code $RM(4,1)$ is 
the set of all $16$ vectors of $F$ and 
the codewords of $RM(4,1)$ are given by hyperplanes. 
It is easy to see that 
$${\Aut}(RM(4,1))\cong GL(5,2)_1=\{ g\in GL(5,2) |
g{}^t(10000)={}^t(10000)\}$$ 
and it is generated by 
$$ \begin{array}{ll}
g(i): & v\in F \to v+v^i \cr
g(i,j): & v\in F \to v+\langle v,v^j\rangle v^i 
\end{array} $$
for $i\not=j$. 

By reversing the above processes, we have a set of mutually 
orthogonal conformal vectors $\{\bar{e}^1,...,\bar{e}^{16} \}$ such that 
$V_{E_8}$ has the following decomposition: 
$$  V_{E_8}=M_{D^1}\oplus {\Ind}_E^{D^1}(H(\hf,0)H(\hf,0)) \mbox{ w.r.t. }
<\bar{e}^i:i=1,...,16>. \eqno{(5.11)} $$
Here $E=H_8\oplus H_8$. 

Choose $g\in {\Aut}(S_{E_8})$.  It is easy to see that $g\in A_{16}$ and so 
$g(e^{(1^{16})})=e^{(1^{16})}$.  By Lemma 3.2, we may assume 
$g\in {\Aut}(D_{E_8})$.  
For an $M_{D_{E_8}}$-module $W$, $g(W)$ denotes an 
$M_{D_{E_8}}$-module defined by 
$v_ng(u)=g(v^g_n(u))$ for $v\in M_{D_{E_8}}$ and $u\in W$.  
Since 
$$g(M_{D_{E_8}})\oplus g(V^{\al}_{E_8})\oplus g(V^{\be}_{E_8})
\oplus g(V^{\al+\be}_{E_8})$$ has a simple VOA structure with a positive 
definite invariant bilinear form, so does 
$$ g(V_{E_8})=\oplus_{\al\in S_{E_8}} g(V_{E_8}^{\al})   $$
by Theorem 3.8.
We note that $g(V_{E_8})$ contains $M_{D_{E_8}}$. 
Using the backward processes according to the sequence 
$$ S^5=g(S^5) \supseteq g(S^4) \supseteq g(S^3) \supseteq g(S^2) \supseteq 
g(S^1),   $$
we obtain a set $\{ \tilde{e}^1,...,\tilde{e}^{16}\}$  
of mutually orthogonal conformal vectors such that 
$g(V_{E_8})$ has the decomposition 
$$ g(V_{E_8}) \cong  M_{D^1}\oplus W  \mbox{ w.r.t. } 
<\tilde{e}^i,i=1,...,16>, $$
where $D^1$ is the set of all even words of length $16$ 
and $W$ is an irreducible 
$M_{D^1}$-module with $\tilde{\tau}(W)=(1^{16})$.  
Since the signs of the actions of $M_{(1^{16})}$ on a module $U$ 
with $\tilde{\tau}(U)$ are changing in each step, we can conclude 
that $W\cong {\Ind}_E^{D^1}(H(\hf,(0^8))H(\hf,(0^8)))$, which 
coincides with (5.11). 
Therefore, there is a VOA isomorphism 
$$ \phi: V_{E_8} \to g(V_{E_8})   $$
such that $\phi(\bar{e}^i)=\tilde{e}^i$  for $i=1,...,16$. 
By the process of changing the coordinate sets according to 
$$ S^1 \supseteq S^2 \subseteq \cdots \subseteq S^5 $$
and
$$g(S^1) \subseteq g(S^2) \subseteq \cdots \subseteq g(S^5),  $$ 
respectively, we have the desired automorphism of $V_{E_8}$. 
\prend

\section{Holomorphic VOA}
Let $V$ be a simple VOA containing a set of mutually orthogonal 
rational conformal vectors $\{e^i:i=1,...,n\}$ with central charge $\hf$ 
such that 
the sum of them is the Virasoro 
element.  Set $P=<\tau_{e^i}:i=1,...,n>$ and let
$V=\oplus_{\chi\in Irr(P)}V^{\chi}$ be the decomposition of $V$ into 
the eigenspaces of $P$.  From Proposition 3.1, 
the space $V^P$ of $P$-invariants is isomorphic to 
$M_D$ for some even linear code $D$ of length $n$. Assign a binary 
word $\al_{\chi}=(a_i)$ by $\chi(e^i)=(-1)^{a_i}$ to $\chi$, 
we can identify $Irr(P)$ and a linear code $S=\{\al_{\chi}:\chi\in Irr(P)\}$.  
As we showed in \cite{Mi5}, $S$ is orthogonal to $D$.  We will treat 
the case $S=D^{\perp}$ in this section.  

\begin{thm}  If $S=D^{\perp}$, then $V$ is the only irreducible $V$-module. 
\end{thm}

\pr
Let $U$ be an irreducible $V$-module.  Since $M_D$ is rational, 
$U$ is a direct sum  
of irreducible $M_D$-modules. Decompose $U$ into the direct sum 
$\oplus U^{\be}$ of $M_D$-modules 
such that $\tilde{\tau}(U^{\be})=\be$. Since $U^{\be}$ is a $M_D$-module, 
$\be\in D^{\perp}=S$ and so $V^{\be}\not=0$. 
Since $U=<v_nu: v\in V^{\al},n\in \Z, \al\in S>$ for any $0\not=u\in U^{\be}$
by \cite{DM2}, $U^{\be}=<v_nu:v\in M_D, n\in \Z>$ and so $U^{\be}$ is 
irreducible 
$M_D$-module.  
Since the restriction 
$$ I\pmatrix{U\cr V\quad U} \to I\pmatrix{U \cr V^{\be} \quad U^{\be}} $$
is injective, $U^{(0^n)}\not=0$. 
So we may assume $\be=(0^n)$.   
Hence $U^{\be}$ is isomorphic to 
a coset module $M_{D+\al}$ for some word $\al\in {\Z}_2^n$.  
Using the skew symmetry, we can define 
a nonzero intertwining operator $I(v,z)\in I_{M_D}\pmatrix{U \cr U\quad V}$ 
with integer powers of $z$. 
By restricting it to $U^{\be}$, we have 
a nonzero intertwining operator 
$I^{\ga}(v,z)\in I_{M_D}\pmatrix{ U^{\ga} \cr M_{\al+D}\quad V^{\ga}}$ for 
$\ga\in S$.
Since its vertex operator has integer 
powers of $z$, $\al$ is orthogonal to $S$ and so 
$\al\in S(P)^{\perp}=D$. Hence $U^{(0^n)}$ is isomorphic to $M_D$.  
Let $q$ be a highest weight vector of $U^{(0^n)}$ corresponding to the 
Vacuum.  Since $L(-1)q=0$, $I(q,z)$ is a scalar and so 
$I(q,z) \in I\pmatrix{U \cr U\quad V}$ gives an $M_D$-
isomorphism of $U$ to $V$. 
This completes the proof of Theorem 6.1. 
\prend

\section{Construction of the moonshine VOA}
In this section, we will construct a VOA 
$V^{\natural}$, which will be proved to be equal to the moonshine 
VOA constructed in \cite{FLM2} in the next section. 
In the section 5, we found 
a set of 16 mutually orthogonal conformal vectors 
$\{e^i:i=1,...,16\}$ of $V_{E_8}$ satisfying the following conditions: \\
(1) \quad $D_{E_8}=RM(4,2)$ \\
(2) \quad $P=<\tau_{e^i}:i=1,...,16>$ has the order $2^5$. \\
(3) \quad $V^P\cong M_{D_{E_8}}$ and $S_{E_8}=D_{E_8}^{\perp}$ 
is generated by 
$$\{ (1^{16}), (0^81^8), (\{0^41^4\}^2), (\{0^21^2\}^4), (\{01\}^8)\} .
\eqno{(7.1)}$$ 
To simplify  the notation, we denote $D_{E_8}$ and $S_{E_8}$ by $D$ and $S$ in 
this section, respectively. 
We note that $D$ and $S$ are $D(5)$ and $S^5$ in the section 4, 
respectively. 
For each codeword $\al\in S$, $V_{E_8}$ 
contains an irreducible $M_D$-module ${V_{E_8}}^{\al}$ such that 
$${V_{E_8}}\cong \boplus_{\al\in S}{V_{E_8}}^{\al} \eqno{(7.2)}$$ 
and ${V_{E_8}}^{(0^{16})}=M_D$.  Since $V_{E_8}$ is a simple VOA, 
Theorem 3.6 implies 
$$V_{E_8}^{\al}\times V_{E_8}^{\be}=V_{E_8}^{\al+\be} $$
for $\al,\be\in S$.

We note that all codewords of $S$ except $(0^{16})$ 
and $(1^{16})$ have weight 8.  
We define a new code $S^{\natural}$ of length 48 by 
$$
S^{\natural}=<(1^{16}0^{16}0^{16}), (0^{16}1^{16}0^{16}), 
(0^{16}0^{16}1^{16}), (\al,\al,\al): \al\in S> . \eqno{(7.3)}$$ 
The weight enumerator of $S^{\natural}$
 is $X^{48}+3X^{32}+120X^{24}+3X^{16}+1$ and there is 
 another expression:
$$S^{\natural}=\{
(\al,\al,\al), \ (\al,\al,\al^c), \ (\al,\al^c,\al), \ 
(\al^c,\al,\al):\al\in S\}.\eqno{(7.4)}$$ 
Set $D^{\natural}=(S^{\natural})^{\perp}$ and call it ``the moonshine code.'' 
Let's explain our choice of the codes $D^{\natural}$ and $S^{\natural}$. 
We may be able to construct the moonshine VOA 
from another pair $(D',S')$, but $(D^{\natural},S^{\natural})$ is very easy 
to handle when we calculate the characters of the elements of the Monster. 
Let's continue the construction. 
$D^{\natural}$ contains 
$D^3=\{(\al, \be, \ga ):\al, \be, \ga \in D\}$ 
and it is easy to see 
$$D^{\natural}=\{(\al,\be,\ga): \al+\be+\ga\in D, \al,\be,\ga 
\mbox{ is even }\}. \eqno{(7.5)}$$ 
Hence $D^{\natural}$ is of dimension 41 and has no codewords of weight 2.
We note that a pair $(D^3,S^{\natural})$ satisfies the conditions (1) and (2) 
in Hypotheses I. 
Denote $(10^{15})$ by $\xi_1$ and set 
$$Q=<(\xi_1\xi_10^{16}),(0^{16}\xi_1\xi_1)>.\eqno{(7.6)} $$ 
To simplify the notation, let $R$ denote a coset module 
$M_{\xi_1+D}$ and $RW$ denote 
the fusion product (tensor product) $R\times W$.  
As we explained in the introduction, our construction 
consists of the following steps. \\
At first, ${V_{E_8}}\ots {V_{E_8}}\ots {V_{E_8}}$ contains a set 
of 48 coordinate conformal vectors 
$$\{e^i\ots {\bf 1}\ots {\bf 1}, \quad 
{\bf 1}\ots e^j\ots {\bf 1}, \quad
{\bf 1}\ots {\bf 1}\ots e^k\quad : i,j,k=1,...,16\},  $$
where ${\bf 1}$ is the Vacuum of ${V_{E_8}}$. 
Decompose it into 
$$ {V_{E_8}}\ots {V_{E_8}}\ots {V_{E_8}}
=\boplus_{\al,\be,\ga\in S}({V_{E_8}}^{\al}\ots {V_{E_8}}^{\be}\ots 
{V_{E_8}}^{\ga}), \eqno{(7.7)}$$
By the fusion rules, 
$$V^1= \boplus_{(\al,\be,\ga)\in S^{\natural}}
({V_{E_8}}^{\al}\ots {V_{E_8}}^{\be}\ots {V_{E_8}}^{\ga}) \eqno{(7.8)}$$
is a sub VOA. 
Let's induce it to  
$$ V^2={\Ind}_{D^3}^{D^3+Q}(V^1). \eqno{(7.9)} $$
We note that since $\langle Q,S^{\natural}\rangle\not=0$, 
a vertex operator of some element in $V^2$ 
does not have integer powers of $z$. 
In particular, $V^2$ is not a VOA. However,  
as $M_{D^3}$-modules, we have 
$$\begin{array}{rl}
&{\Ind}_{D^3}^{D^3+Q}(V_{E_8}^\al\ots V_{E_8}^\be\ots V_{E_8}^\ga) \cr
=&(V_{E_8}^\al\otimes V_{E_8}^\be\otimes V_{E_8}^\ga)
\oplus (RV_{E_8}^\al\ots RV_{E_8}^\be\ots V_{E_8}^\ga) \cr
&\oplus (V_{E_8}^\al\ots RV_{E_8}^\be\ots RV_{E_8}^\ga) \oplus
(RV_{E_8}^\al\ots V_{E_8}^\be\ots RV_{E_8}^\ga). 
\end{array} $$
Using (7.4), define $W^{(\al,\be,\ga)}$ for 
$(\al,\be,\ga)\in S^{\natural}$ as follows: 
$$ \begin{array}{l}
W^{(\al,\al,\al)}={V_{E_8}}^{\al}\ots {V_{E_8}}^{\al}\ots 
{V_{E_8}}^{\al}, \cr
W^{(\al,\al,\al^c)}=(R{V_{E_8}}^{\al}) \ots (R{V_{E_8}}^{\al})
\ots {V_{E_8}}^{\al^c}, \cr
W^{(\al,\al^c,\al)}=(R{V_{E_8}}^{\al}) \ots V_{E_8}^{\al^c}\ots 
(R{V_{E_8}}^{\al}),\cr
W^{(\al^c,\al,\al)}=V_{E_8}^{\al^c}\ots 
(R{V_{E_8}}^{\al}) \ots (R{V_{E_8}}^{\al}). 
\end{array} \eqno{(7.10)}$$
Since  all $RV_{E_8}^{\al}$ are irreducible $M_D$-modules by Theorem 3.3,
$W^{(\al,\be,\ga)}$ are all irreducible $M_{D^3}$-modules.  
Induce them into 
$$V^{\chi}={\Ind}_{D^3}^{D^{\natural}}(W^{\chi}) \eqno{(7.11)}$$ 
for $\chi\in S^{\natural}$. 
Finally, set 
$$V^{\natural}=\boplus_{\chi\in S^{\natural}}(V^{\chi}).\eqno{(7.12)}$$
This is the desired Fock space. 

Since $(D^{\natural},S^{\natural})$ satisfies the conditions (1) and (2) 
of Hypotheses I, the remaining thing we have to do is to prove that 
$$V^{\chi,\mu}=
M_{D^{\natural}}\oplus V^{\chi}\oplus V^{\mu}\oplus V^{\chi+\mu} $$  
has a simple VOA structure with a positive definite invariant 
bilinear form for any $\mu,\chi\in S^{\natural}$ with $\dim <\mu,\chi>=2$. 
We note that since $M_{D^3}\oplus W^{(\al,\al,\al)}$ and 
$M_{D^{\natural}}\oplus W^{(\al,\al,\al^c)}$ are sub VOAs of 
${\Ind}_{D^3}^{<D^3,(\xi_1\xi_10^{16})>}(M_{D^3}\oplus W^{(\al,\al,\al)})$, 
they have simple VOA structures with positive definite invariant bilinear 
forms. 
Take a sub VOA 
$$(V^1)^{\chi,\mu}=M_{D^3}\oplus (V^1)^{\chi}\oplus (V^1)^{\mu}
\oplus (V^1)^{\chi+\mu}$$ 
of $V^1$ and set
$$W^{\chi,\mu}=M_{D^3}\oplus W^{\chi}\oplus W^{\mu}\oplus W^{\chi+\mu}$$ 
using (7.10). 
If $<\chi,\mu>$ is orthogonal to $(\xi_1\xi_10^{16})$, then  
${\Ind}_{D^3}^{D^3+<(\xi_1\xi_10^{16})>}((V^1)^{\chi,\mu})$ is a VOA 
with the desired properties. Moreover it contains 
$W^{\chi,\mu}$ as a sub VOA. 
Similarly, if $<\chi,\mu>$ is orthogonal to $(0^{16}\xi_1\xi_1)$ or 
$(\xi_10^{16}\xi_1)$, then we have the desired properties. 
Therefore we may assume that 
$\chi=(\al,\al,\al^c)$ and $\mu=(\be,\be^c,\be)$. 
Set $\ga=\al^c+\be$ and assume that $Supp(\al)\cap Supp(\be)\not=\emptyset$. 
Choose $t\in Supp(\al)\cap Supp(\be)$. Then $t\in Supp(\ga)$.  
Set $\xi_t=(0^{t-1}10^{15-t})$ and $R^t=M_{D+\xi_t}$. 
Since 
$$ 
(\xi_t\xi_t0^{16})+(\xi_1\xi_10^{16})
\in D^{\natural}, \quad 
(\xi_t0^{16}\xi_t)+(\xi_10^{16}\xi_1)
\in D^{\natural}, \quad
(0^{16}\xi_t\xi_t)+(0^{16}\xi_1\xi_1)
\in D^{\natural}, $$ 
we have 
$$ \begin{array}{rl}
{\Ind}_{D^3}^{D^{\natural}}(R^t{V_{E_8}}^{\al}\ots R^t{V_{E_8}}^{\al}
\ots {V_{E_8}}^{\al^c})
&={\Ind}_{D^3}^{D^{\natural}}(R{V_{E_8}}^{\al}\ots R{V_{E_8}}^{\al}
\ots {V_{E_8}}^{\al^c}), \cr
{\Ind}_{D^3}^{D^{\natural}}(R^t{V_{E_8}}^{\be}\ots {V_{E_8}}^{\be^c}
\ots R^t{V_{E_8}}^{\be})
&={\Ind}_{D^3}^{D^{\natural}}(R{V_{E_8}}^{\be}\ots {V_{E_8}}^{\be^c}\ots 
R{V_{E_8}}^{\be}), \cr
{\Ind}_{D^3}^{D^{\natural}}({V_{E_8}}^{\ga^c}\oplus R^t{V_{E_8}}^{\ga}\ots 
R^t{V_{E_8}}^{\ga})
&={\Ind}_{D^3}^{D^{\natural}}({V_{E_8}}^{\ga^c}\otimes 
R{V_{E_8}}^{\ga}\ots R{V_{E_8}}^{\ga}). \end{array} $$
Set
$$ \ga_1=(\xi_t\xi_t0^{16}), \quad \ga_2=(\xi_t0^{16}\xi_t), \quad 
\ga_3=(0^{16}\xi_t\xi_t). $$

Since $Supp(\ga_1)\subseteq Supp(\chi)$, $Supp(\ga_2)\subseteq Supp(\mu)$ and 
$Supp(\ga_3)\subseteq Supp(\chi+\mu)$, it follows from Lemma 3.13 that 
$$ \begin{array}{l}
R^t(V_{E_8})^{\al}\ots R^t(V_{E_8})^{\al}\ots (V_{E_8})^{\al^c} 
\cong \sigma_{\ga_1}((V^1)^{(\al,\al,\al^c)}), \cr
R^t(V_{E_8})^{\be}\ots (V_{E_8})^{\be^c}
\ots R^t(V_{E_8})^{\be} \cong \sigma_{\ga_2}(V^1)^{(\be,\be^c,\be)}, \cr
(V_{E_8})^{\ga^c}\ots R^1(V_{E_8})^{\ga}
\ots (V_{E_8})^{\ga} \cong \sigma_{\ga_3}(V^1)^{(\ga^c,\ga,\ga)}. 
\end{array} $$
Since $M_{D^3}\oplus (V^1)^{(\al,\al,\al^c)}\oplus (V^1)^{(\be,\be^c,\be)}
\oplus (V^1)^{(\ga^c,\ga,\ga)}$ has a simple VOA structure with a positive 
definite invariant bilinear form, so does 
$M_{D^3}\oplus \sigma_{\ga^1}(V^1)^{(\al,\al,\al^c)}\oplus 
\sigma_{\ga^2}((V^1)^{(\be,\be^c,\be)})
\oplus \sigma_{\ga^1+\ga^2}((V~1)^{(\ga^c,\ga,\ga)})$. 
Since $\ga^1+\ga^2+\ga^3=0$, 
$\sigma_{\ga^1+\ga^2}(V^1)^{(\ga^c,\ga,\ga)}=
\sigma_{\ga^3}(V^1)^{(\ga^c,\ga,\ga)}$. Hence  
$W^{\al,\be}=M_{D^3}\oplus W^{(\al,\al,\al^c)}\oplus 
W^{(\be,\be^c,\be)}
\oplus W^{(\ga^c,\ga,\ga)}$ has the desired VOA structure and so does 
$(V^{\natural})^{\chi,\mu}$.

Hence we assume $Supp(\al)\cap Supp(\be)=\emptyset$. Then 
one of $\{\al, \be, \al+\be^c\}$ is at least $(0^{16})$ since 
$\al,\be\in S$. Set $\ga=\al+\be^c$. 
So we may assume $\al=(0^{16})$ and $\ga^c=\be$. 
It follows from the structure of $D$ that 
there is a self dual subcode $E$ of $D^3$ which is a 
direct sum $\boplus_{i=1}^6E^i$ of 6 $[8,4,4]$-Hamming codes $E^i$ such that 
$E_{\delta}=\{\al\in E|Supp(\al)\subseteq Supp(\delta)\}$ 
is a direct factor of $E$ for any $\delta\in <\be,\ga>$. 
In particular, there are $M_E$-modules 
$U^{\al}$, $U^{\be}$, $U^{\ga}$ such that 
$$ \begin{array}{l}
{\Ind}_{E}^{D^{\natural}}(U^{\al})=(V^{\natural})^{(\al,\al,\al^c)},\cr 
{\Ind}_{E}^{D^{\natural}}(U^{\be})=(V^{\natural})^{(\be,\be^c,\be)}, \cr 
{\Ind}_{E}^{D^{\natural}}(U^{\ga})=(V^{\natural})^{(\ga^c,\ga,\ga)}.  
\end{array} $$
In the following, we assume $|\be|=8$. 
We can prove the assertion for $\be=(0^{16})$ or $\be=(1^{16})$ by 
the similar arguments. We may assume $\be=(1^80^8)$. 
As we showed in \S 5, we have a VOA  
$U={V_{E_8}}^{(0^{16})}\oplus {V_{E_8}}^{(1^80^8)}\oplus 
{V_{E_8}}^{(0^81^8)}\oplus {V_{E_8}}^{(1^{16})}$ with a positive 
definite invariant bilinear form such that 
$$ \begin{array}{l}
{V_{E_8}}^{(0^{16})}={\Ind}_{F}^D(H(\hf,(0^8))\otimes H(\hf,(0^8))) \cr
{V_{E_8}}^{(1^80^8)}={\Ind}_{F}^D(H(\st,\xi_1)\otimes H(\hf,\xi_1)) \cr
{V_{E_8}}^{(0^81^8)}={\Ind}_{F}^D(H(\hf,\xi_1)\otimes H(\st,\xi_1))\cr
{V_{E_8}}^{(1^{16})}={\Ind}_{F}^D(H(\st,(0^8))\otimes H(\st,(0^8))), 
\end{array} $$
where $F=D_{(1^80^8)}\oplus D_{(0^81^8)}$ is a direct sum of two 
Hamming codes. 
In order to simplify the notation, we omit the notation $"\ots"$ 
between $H(\ast,\ast)$ and $H(\ast,\ast)$. 
In particular, \\ 
$\tilde{U}=H(\hf,(0^8))H(\hf,(0^8))\oplus H(\st,\xi_1)H(\hf,\xi_1) 
\oplus H(\hf,\xi_1)H(\st,\xi_1) \oplus H(\st,(0^8))H(\st,(0^8)) $ 
has a VOA structure with a positive definite invariant bilinear form. 
Since $W^{(\al,\al,\al^c)}$ is given by 
$R{V_{E_8}}^{\al}\otimes 
R{V_{E_8}}^{\al}\otimes {V_{E_8}}^{\al^c}$, 
$$ \begin{array}{l}
U^{\al}=H(\hf,\xi_1)H(\hf,(0^8))H(\hf,\xi_1)H(\hf,(0^8))
H(\st,(0^8))H(\st,(0^8)).
\end{array}$$ 
We similarly obtain 
$$ \begin{array}{l}
U^{\be}
=H(\st,(0^8))H(\hf,\xi_1)H(\hf,\xi_1)H(\st,\xi_1)H(\st,(0^8))H(\hf,\xi_1) \cr
U^{\ga}=
H(\st,\xi_1)H(\hf,\xi_1)H(\hf,\xi_1)H(\st,\xi_1)H(\hf,\xi_1)H(\st,\xi_1).
\end{array}$$
By changing the order of the components, 
$(123456) \to (243516)$, we have 
$$ \begin{array}{ll}
M_{E}=&
H(\hf,(0^8))H(\hf,(0^8))H(\hf,(0^8))H(\hf,(0^8))H(\hf,(0^8))H(\hf,(0^8)), \cr
U^{\al}=&
H(\hf,(0^8))H(\hf,(0^8))H(\hf,(\xi_1))H(\st,(0^8))H(\hf,(\xi_1))H(\st,(0^8)),
\cr
U^{\be}=&
H(\hf,(\xi_1))H(\st,(\xi_1))H(\hf,(\xi_1))H(\st,(0^8))H(\st,(0^8))
H(\hf,(\xi_1)), \cr
U^{\ga}=&
H(\hf,(\xi_1))H(\st,(\xi_1))H(\hf,(0^8))H(\hf,(\xi_1))H(\st,(\xi_1))
H(\st,(\xi_1)).
\end{array} $$

By Lemma 3.8, there is another coordinate set of conformal vectors 
$\{d^1,...,d^8\}$ in $M_{H_8}$  such that 
$$\begin{array}{lcl}
H(\hf,(\xi_1))\mbox{ w.r.t. }<e^i> &\cong &H(\st,(0^8)) \mbox{ w.r.t. }
<d^i> \cr
H(\st,\xi_1) \mbox{ w.r.t. }<e^i> &\cong & H(\hf,(\xi_1)) \mbox{ w.r.t. }
<d^i> \cr
H(\st,(0^8)) \mbox{ w.r.t. }<e^i> &\cong &H(\st,(\xi_1))\mbox{ w.r.t. }
<d^i>. 
\end{array}$$
Changing the coordinate sets, we have 
$$ \begin{array}{l}
{\tilde{V}_{E_8}}^{(0^{16})}={\Ind}_F^D(H(\hf,(0^8))\oplus H(\hf,(0^8))) \cr
{\tilde{V}_{E_8}}^{(1^80^8)}={\Ind}_{F}^D(H(\hf,(\xi_1))\ots H(\st,(0^8))) \cr
{\tilde{V}_{E_8}}^{(0^81^8)}={\Ind}_{F}^D(H(\st,(0^8))\ots H(\hf,(\xi_1))) \cr
{\tilde{V}_{E_8}}^{(1^{16})}={\Ind}_{F}^D(H(\st,(\xi_1))\ots H(\st,(\xi_1)))
\end{array} $$
with respect to $\{d^1,...,d^8,d^9,...,d^{16}\}$. 
Therefore, 
$\bar{U}=M_{E}\oplus U^{\al}\oplus U^{\be}\oplus U^{\ga}$ 
is a subset of a VOA 
$\tilde{V}_{E_8}\otimes \tilde{V}_{E_8}\otimes \tilde{V}_{E_8}$. It is also 
easy to check that $\bar{U}$ is closed under the products. Hence, $\bar{U}$ is a VOA with a positive 
definite invariant bilinear form and so does $(V^{\natural})^{\chi,\mu}=
{\Ind}_{E}^{D^{\natural}}(\bar{U})$. 
This completes the construction of $V^{\natural}$. 
\prend

\begin{cry}  
$V^{\natural}$ has a positive definite invariant bilinear 
form. 
\end{cry}

\begin{rmk}
Because of our construction, a VOA satisfying Hypotheses I is a direct sum of 
the tensor product of $L(\hf,0),L(\hf,\hf),L(\hf,\st)$ and we know 
the multiplicities of irreducible $L(\hf,0)^{\otimes n}$-modules by 
Theorem 3.2, (c.f. Corollary 5.2 in \cite{Mi3}). Hence it is not difficult 
to calculate its character 
$$  ch_V(z)= e^{2\pi iz/(rank(V))}(\sum_{n=0}\dim V_n\ e^{2\pi iz} ).  $$
\end{rmk}

For example, let's show $(V^{\natural})_1=0$. 
We first have $(M_{{D^{\natural}}})_1=0$ 
since ${D^{\natural}}$ has no codewords of weight $2$. 
Also, if $(V^{\natural})^{\chi}_1\not=0$, then 
the weight of $\chi$ is equal to $16$ and so 
$\chi$ is one of 
$(1^{16}0^{16}0^{16})$, $(0^{16}1^{16}0^{16})$ or $(0^{16}0^{16}1^{16})$.  
Say $\chi=(1^{16}0^{16}0^{16})$. 
Since $(V^{\natural})^{\chi}=
{\Ind}_{D_{E_8}^3}^{D^{\natural}}(V_{E_8}^{(1^{16})}\otimes 
M_{D_{E_8}+\xi_1}\otimes M_{D_{E_8}+\xi_1})$ and $D^{\natural}$ 
does not contains any words of the form $(\al,\xi_1,\xi_1)$, 
the minimal weight of 
$(V^{\natural})^{\chi}$ is greater than $1$. 
Therefore, we obtain $V^{\natural}_1=0$.

\section{Conformal vectors}
Since each rational conformal vector $e\in V$ with central charge 
$\hf$ offers an automorphism $\tau_e$, 
it is very important to find such conformal vectors for studying 
the automorphism group ${\Aut}(V)$.
Therefore, we will construct several conformal vectors of $V^{\natural}$ 
explicitly. 

\subsection{Case I}
Set $D_1=<H_8\oplus H_8,(\xi_1\xi_1)>$ and 
$S=<(1^{16})>$, where $\xi_1=(10^7)$. 
Then the pair $(D_1,S)$ satisfies the conditions (1) and (2) 
of Hypotheses I. 
Set 
$$U=H(\hf,0)H(\hf,0)\oplus H(\hf,\xi_1)H(\hf,\xi_1)\oplus 
H(\st,\xi_1)H(\st,0)\oplus H(\st,0)H(\st,\xi_1). $$ 
$U$ is a sub VOA of $V_{E_8}$. It is easy to see that 
$\dim (H(\hf,0)H(\hf,0))_1=0$ and \\
$\dim (H(\hf,\xi_1)H(\hf,\xi_1))_1
=\dim (H(\st,\xi_1)H(\st,0))_1=\dim (H(\st,0)H(\st,\xi_1))_1=1$. Hence   
$U_1$ is isomorphic to $sl(2)$.  Viewing  
$(H(\hf,\xi_1)H(\hf,\xi_1))_1$ 
as a Cartan subalgebra of $sl(2)$, $H(\st,\xi_1)H(\st,0)\oplus 
H(\st,0)H(\st,\xi_1)$ contains two roots $\iota(x)$ and $\iota(-x)$. 
Take a sub lattice VOA of type $A_1$ generated by $U_1$, we may 
obtain the following elements: 
$$ \begin{array}{l}
x(-1){\bf 1}\in (H(\hf,\xi_1)H(\hf,\xi_1))_1, \cr
\iota(x)+\iota(-x)\in (H(\st,\xi_1)H(\st,0))_1,\mbox{ and } \cr
\iota(x)-\iota(-x)\in (H(\st,\xi_1)H(\st,0))_1. 
\end{array} $$
Take another copy of them and set 
$$ \begin{array}{l}
y(-1){\bf 1}\in (H(\hf,\xi_1)H(\hf,\xi_1))_1, \cr
\iota(y)+\iota(-y)\in (H(\st,\xi_1)H(\st,0))_1,\mbox{ and } \cr
\iota(y)-\iota(-y)\in (H(\st,\xi_1)H(\st,0))_1.
\end{array} $$
Then we have  
$$ \begin{array}{lll}
\iota(\pm x)\ots \iota(\pm y)+\iota(\mp x)\ots 
\iota(\mp y) &\in&
H(\st,0)H(\st,\xi_1)H(\st,0)H(\st,\xi_1)\cr
&&\oplus H(\st,\xi_1)H(\st,0)H(\st,\xi_1)H(\st,0)\cr
x(-1)y(-1)&\in& 
H(\hf,\xi_1)H(\hf,\xi_1)H(\hf,\xi_1)H(\hf,\xi_1), \mbox{ and }\cr
x(-1)^2{\bf 1},\  y(-1)^2{\bf 1}&\in& 
H(\hf,0)H(\hf,0)H(\hf,0)H(\hf,0).
\end{array} $$
Since $\langle x\pm y, x\pm y\rangle=2$, 
$e^+(x\pm y)={1\over 16}(x\pm y)(-1)^2{\bf 1}+ 
{1\over 4}(\iota(x\pm y)+\iota(-x\mp y))$ and 
$e^-(x\pm y)={1\over 16}(x\pm y)(-1)^2{\bf 1}- 
{1\over 4}(\iota(x\pm y)+\iota(-x\mp y))$ 
are rational conformal vectors with central charge $\hf$. 
Hence, we obtain 
four rational conformal vectors $e^{\pm}(x\pm y)$ in 
$$\begin{array}{rl}
&H(\hf,0)H(\hf,0)H(\hf,0)H(\hf,0)\oplus
H(\hf,\xi_1)H(\hf,\xi_1)H(\hf,\xi_1)H(\hf,\xi_1) \cr
\oplus &H(\st,0)H(\st,\xi_1)H(\st,0)H(\st,\xi_1)\oplus 
H(\st,\xi_1)H(\st,0)H(\st,\xi_1)H(\st,0).
\end{array}$$ 

\subsection{Case II}
We first treat the first component $V_{E_8}\otimes {\bf 1}\otimes {\bf 1}$ 
of $V_{E_8}\ots V_{E_8} \ots V_{E_8}$. We denote $D_{E_8}$, $S_{E_8}$ 
and $V_{E_8}$ 
by $D$, $S$, $V$ here, respectively. 
Let $\al,\be\in S$ so that $|\al|=|\be|=|\al+\be|$. By rearranging 
the coordinate sets, we may assume $\al=(1^80^8), \be=(1^40^41^40^4)$. 
As we showed, $V$ contains a sub VOA 
$$U=M_D\oplus V^{\al^c}\oplus V^{\be^c}\oplus 
V^{\al+\be} $$
for $\al,\be\in S$. 
Since $D_{\al}$, $D_{\be}$ and $D_{\al+\be}$ are all isomorphic 
to $H_8$, 
the multiplicities 
of the irreducible $L(\hf,0)^{\ots 8}$-modules in 
$V^{\al^c}\oplus V^{\be^c}\oplus V^{\al+\be}$ 
are all one by Theorem 3.2. Hence $\dim (V^{\al^c})_1=
\dim (V^{\be^c})_1=\dim (V^{\al+\be})_1=8$. 
Since $D$ does not contain any words of weight 2, 
$(M_D)_1=0$ and so $(V^{\al^c})_1$, $(V^{\be^c})_1$ and 
$(V^{\al+\be})_1$ are all commutative Lie algebras. 
Since $U$ is a sub VOA of a lattice VOA $V$ of rank 8 and so 
$U_1$ is isomorphic to $sl(2)^8$. 
Viewing $(V^{\al+\be})_1$ as a Cartan subalgebra and 
embedding it into a lattice VOA $V_{A_1^8}$ of root lattice $A_1^8$, 
we denote the positive roots by $\iota(x_1),...,\iota(x_8)$ 
and the negative roots by $\iota(-x_1),...,\iota(-x_8)$.
In addition, we may assume 
$$ \begin{array}{lll}
x_i(-1) &\in &V^{\al+\be}_1 \cr
\iota(x_i)+\iota(-x_i) &\in &V^{\al^c}_1 \cr
\iota(x_i)-\iota(-x_i) &\in &V^{\be^c}_1
\end{array}  $$
for $i=1,...,8$. 

We next treat the second and third components of $V_{E_8}\ots 
V_{E_8}\ots V_{E_8}$. 
By the similar arguments as in the construction of the moonshine VOA, 
$M_{D^2}\oplus W^{(\al,\al)}\ots W^{(\be,\be)}\oplus W^{(\al+\be,\al+\be)}$ 
has a simple VOA structure, where 
$W^{(\al,\al)}=RV_{E_8}^{\al}\ots RV_{E_8}^{\al}$, 
$W^{(\be,\be)}=RV_{E_8}^{\be}\ots RV_{E_8}^{\be}$ and  
$W^{(\al+\be,\al+\be)}=V_{E_8}^{\al+\be}\ots V_{E_8}^{\al+\be}$. 
Set $F=\{(\al',\be'):\al'+\be'\in D, \al',\be' \mbox{ even }\}$. 
Then $M_F$ does not contain any roots and $D\oplus F\subseteq D^{\natural}$.
Set $U^{\ga,\ga}={\Ind}_{M_{D^2}}^{M_F}(W^{\ga,\ga})$ for $\ga\in \{\al,\be,\al+\be\}$. 
By Theorem 3.7, we have a VOA 
$U=M_{F}\oplus U^{\al,\al}\ots U^{\be,\be}\oplus U^{\al+\be,\al+\be}$. 
Since $|F_{(\al\al)}|=|F_{(\be\be)}|=|F_{((\al+\be)(\al+\be))}|=2^{11}$, 
the multiplicities of irreducible $L(\hf,0)^{\ots 16}$-submodules 
is 8.  Hence, $\dim (U^{\ga,\ga}_1)=8$ for $\ga\in \{\al,\be,\al+\be\}$. 
Set 
$U^{(1^{32})}={\Ind}_{M_{D^2}}^{M_F}(R{V_{E_8}}^{(1^{16})}
\ots R{V_{E_8}}^{(1^{16})})$, 
$U^{\al^c,\al^c}={\Ind}_{M_{D^2}}^{M_F}({V_{E_8}}^{\al}\ots {V_{E_8}}^{\al})$, 
$U^{\be^c,\be^c}
={\Ind}_{M_{D^2}}^{M_F}({V_{E_8}}^{\be}\ots {V_{E_8}}^{\be})$ and  
$U^{\al^c+\be,\al^c+\be}={\Ind}_{M_{D^2}}^{M_F}(R{V_{E_8}}^{\al+\be^c}\ots 
R{V_{E_8}}^{\al+\be^c})$. Then, 
$X=M_{F}\oplus U^{(1^{16}),(1^{16})}\oplus U^{\al,\al}\oplus  
U^{\be,\be}\oplus U^{\al+\be,\al+\be}
\oplus U^{\al^c\al^c}\oplus U^{\be^c,\be^c}\oplus U^{\al+\be^c,\al+\be^c}$ 
has a VOA structure.  Since $(M_F\oplus U^{(1^{16}),(1^{16})})_1=0$, 
$(U^{\al+\be^c,\al+\be^c}\oplus U^{\al+\be,\al+\be})_1$ is of dimension 
16. 
Since $X$ is a sub VOA of a lattice VOA $V$ of rank 16, 
$X_1$ is isomorphic to $sl(2)^{16}$ and 
$U_1$ is isomorphic to $sl(2)^8$.  Viewing 
$(U^{\al+\be,\al+\be})_1$ as a Cartan subalgebra and embedding them 
in a lattice VOA $V_{A_1^8}$ of the root lattice $A_1^8$, 
we denote the positive roots by $\iota(y_1),...,\iota(y_8)$ 
and the negative roots by $\iota(-y_1),...,\iota(-y_8)$. 
We may also assume that  
$$ \begin{array}{lll}
y_i(-1) &\in &(U^{\al+\be,\al+\be})_1 \cr
\iota(y_i)+\iota(-y_i) &\in &(U^{\al^c,\al^c})_1 \cr
\iota(y_i)-\iota(-y_i) &\in &(U^{\be^c,\be^c})_1
\end{array}  $$
for $i=1,...,8$. 

Set 
$$ \begin{array}{l}
W^{\al}=V_{E_8}^{\al^c}\ots U^{\al,\al}, \cr
W^{\be}=V_{E_8}^{\be^c}\ots U^{\be,\be} \mbox{ and }\cr
W^{\al+\be}=V_{E_8}^{\al+\be}\ots U^{\al+\be,\al+\be}.
\end{array} $$ 
Then 
$$V^1=M_{D\oplus F}\oplus W^{\al}\oplus W^{\be}\oplus W^{\al+\be} $$
is a sub VOA of $V^{\natural}$.  
We have 
$$ \begin{array}{l}
x_i(-1)^2 \in M_D, \cr
y_i(-1)^2 \in M_F, \cr
x_i(-1)y_i(-1)\in W^{\al+\be}, \cr
(\iota(x^i)+\iota(-x^i))(\iota(y^i)+\iota(-y^i))\in W^{\al} \mbox{ and }, \cr
(\iota(x^i)-\iota(-x^i))(\iota(y^i)-\iota(-y^i))\in W^{\be}.
\end{array}$$ 
By the same arguments as in the case I, 
we have 32 mutually orthogonal conformal vectors
$$\begin{array}{l}
d^{4i-3}={1\over 16}(x^i+y^i)(-1)^2{\bf 1}+{1\over 4}(\iota(x^i+y^i)+
\iota(-x^i-y^i)) \cr
d^{4i-2}={1\over 16}(x^i+y^i)(-1)^2{\bf 1}-{1\over 4}
(\iota(x^i+y^i)+\iota(-x^i-y^i))\cr
d^{4i-1}={1\over 16}(x^i-y^i)(-1)^2{\bf 1}+{1\over 4}(\iota(x^i-y^i)+
\iota(-x^i+y^i)) \cr
d^{4i}={1\over 16}(x^i-y^i)(-1)^2{\bf 1}-{1\over 4}
(\iota(x^i-y^i)+\iota(-x^i+y^i))
\end{array}$$ 
in $V^1$, where $\iota(x^i+y^i)$ denotes $\iota(x^i)\otimes \iota(y^i)$.

\section{Automorphism group}
In this section, we will prove that the full automorphism 
group of $V^{\natural}$ is the Monster simple group. 
We first quote the following two theorems about the 
finiteness of automorphism group from \cite{Mi4}. \\

\noindent
{\bf Hypotheses II} \\
(1)  $V=\sum_{i=0}^{\infty} V_i$ is a VOA over ${{\R}}$.  \\
(2)  $\dim V_0=1$.  \\
(3)  $V_1=0$.  \\
(4)  $V$ has a positive definite invariant bilinear form 
$\langle \  ,\  \rangle$. \\
(5)  The Virasoro element is a sum of mutually 
orthogonal conformal vectors with central charge ${1\over 2}$. \\

Under the above Hypotheses II, we recall the following results 
from \cite{Mi4}. 

\begin{thm} 
Let $e$, $f$ be two distinct conformal vectors with 
central charge ${1\over 2}$. Then we have 
$$\langle e,f\rangle\leq {1\over 12} \quad \mbox{ and } \quad 
\langle e-f,e-f \rangle \geq {1\over 3}. $$
In particular, 
there are only finitely many conformal vectors with central charge ${1\over 2}$.\end{thm}

\pr
Using the product $ab=a_1b$ and the inner product 
$\langle a,b\rangle{\bf 1}=a_3b$ for $a,b\in V_2$, 
$V_2$ becomes a commutative algebra called Griess algebra. 
Let $V_2={{{\R}}}e\oplus {{{\R}}}e^{\perp}$ be the decomposition of $V_2$, 
where ${\R}e^{\perp}=\{v\in V_2|\langle v,e\rangle=0\}$.
For $f$, there are $r\in {\R}$ and $w\in {{{\R}}}e^{\perp}$ such that 
$$f=re+w. $$ 
Since $\langle ew,e\rangle=\langle w,e^2\rangle=\langle w,2e\rangle=0$, 
we have $ew\in {{{\R}}}e^{\perp}$ and so
$$2re+2w=2f=f^2=\{r^22e+ w^2_e\} +\{(w^2-w^2_e)+2rew\}, $$
where ${w^2}_e$ denotes the first entry of $w^2$ in the decomposition 
${{{\R}}}e\oplus {{{\R}}}e^{\perp}$. Hence, 
$$r^2/2+\langle e,w_e^2\rangle=\langle e,2r^2e+w^2_e \rangle 
=\langle e, f^2\rangle=\langle e,2f\rangle=\langle e,2re\rangle =r/2$$ 
and so $\langle e,w_e^2\rangle=r(1-r)/2$.
On the other hand, 
$${1\over 4}=\langle f,f\rangle=r^2{1\over 4}+\langle w,w\rangle, $$
and so $\langle w,w\rangle={1\over 4}(1-r^2)$. Since 
$<e>\cong L(\hf,0)$ and every irreducible $L(\hf,0)$-module is isomorphic 
to one of $L(\hf,0), L(\hf,\hf), L(\hf,\st)$ and ${\bf w}-e$ 
is a sum of rational conformal vectors with central charge $\hf$, 
the eigenvalues of 
$L(0)-e_1$ is nonnegative. 
Hence, the eigenvalues of $e_1$ on $V_2$ is $0,1,2,\hf,
\hf+1,\st,\st+1$. If $e_1v=(\hf+1)v$ or $(\st+1)v$, then 
then $e_2v\not=0$, which contradicts to $e_2v\in V_1=0$. If $e_1v=2v$, 
then $e_3v\not=0$, which contradicts to $v\in {\R}e^{\perp}$. 
If $e_1v=v$, then $v\in (L(\hf,0))_1=0$. Hence, 
the eigenvalues of $e$ on 
${\R}e^{\perp}$ are $0$, ${1\over 2}$, or ${1\over 16}$.  
Hence, we obtain 
$$r/2-r^2/2=\langle e,{w^2}_e\rangle
=\langle e,w^2\rangle=\langle we,w\rangle 
\leq {1\over 2}\langle w,w\rangle={1\over 8}(1-r^2)$$
and so $3r^2-4r+1\geq 0$. This implies 
$r\geq 1$ or $r\leq {1 \over 3}$.
If $r\geq 1$, then it contradicts $\langle w,w\rangle > 0$.  
We hence have $r\leq {1\over 3}$ and so 
$\langle e,f\rangle \leq {1\over 12}$, which implies 
$\langle e-f,e-f\rangle\geq {1\over 3}$. Hence, 
there are only finitely many conformal vectors with 
central charge ${1\over 2}$ since $\{v\in V_2|\langle v,v\rangle=4\}$ 
is a compact space. 
\prend

\begin{thm}
If $V$ satisfies Hypothesis II, then 
${\Aut}(V)$ is finite.
\end{thm}

\pr
Suppose false and let $G$ be an automorphism group of $V$ of infinite order. 
Since $G$ acts on the set $J$ of all conformal vectors with 
central charge ${1\over 2}$ and $J$ is a finite set by Theorem 9.1, 
we may assume that $G$ fixes all conformal vectors with central 
charge $\hf$. 
In particular, $G$ fixes all coordinate conformal vectors $e^i$ for 
$i=1,...,n$. 
Set $P=<\tau_{e^i}:i=1,...,n>$. By the definition of $\tau_{e^i}$, 
$P$ is an elementary abelian 2-group.  Let $V=\oplus_{\chi\in Irr(P)}V^{\chi}$ 
be the decomposition of $V$ into the eigenspaces of $P$, where 
$Irr(P)$ is the set of all linear characters of $P$ and 
$V^{\chi}=\{v\in V: gv=\chi(g)v\ \forall g\in P\}$.  As we mentioned 
in the introduction, $\tilde{\tau}(V^{\chi})=(a_1,...,a_n)$ is given by 
$(-1)^{a_i}=\chi(e^i)$.  Since $G$ fixes all $e^i$ and 
$g^{-1}\tau_{e^i}g=\tau_{g(e^i)}$ for $g\in {\Aut}(V)$ by the definition,  
$[G,P]=1$ and so $G$ leaves all $V^{\chi}$ 
invariant. In particular, $G$ acts on $V^{1_G}$. 
We think over the action of $G$ on $V^{1_G}\ (=V^P)$ for a while. 
Set $T=<e^1,...,e^n>\cong L(\hf,0)^{\otimes n}$.  Since $\dim V_0=1$, 
$T$ is the only irreducible $T$-submodule of $V$ isomorphic to 
$L(\hf,0)^{\otimes}$.  
By the hypotheses, $V$ has a positive definite invariant bilinear form and 
so $V^P$ is 
simple. Hence, $V^P$ is isomorphic to a code VOA $M_D$ for some 
even linear code $D$. 
In particular, $V^P$ is a direct sum of 
finite distinct irreducible $T$-modules $M_{\al}$.
Since $T$ is generated by $\{e^i:i=1,...,n\}$ and $G$ fixes all $e^i$, 
$G$ fixes all vectors of $T$ and so the action of $g\in G$ on 
$M_{\al}$ is a scalar $\lambda_{\al}$.  Since $V$ has a positive 
definite invariant bilinear form, we have
$0\not=\langle v,v\rangle=\langle g(v),g(v)\rangle
=\la_{\al}^2\langle v,v\rangle$ 
and so $\la_{\al}=\pm 1$. Since $|D|$ is finite, 
a finite index subgroup of $G$ fixes all vectors 
of $V^P$.  So we may assume that $G$ fixes all vectors in $V^P$. 
Since $V^{\chi}$ is a irreducible $V^P$-module by 
\cite{DM2}, $g\in G$ acts on 
$V^{\chi}$ as a scalar $\lambda_{\chi}$. By the same arguments as above, 
we have a contradiction. 
\prend

In \S 3, we proved that we can induce every automorphism of $D$ 
into an automorphism of $M_D$.  We will here show that we can induce 
every automorphism of $S^{\natural}$ into an automorphism of $V^{\natural}$. 

\begin{lmm}  For any $g\in {\Aut}(S^{\natural})$, there is  
an automorphism $\tilde{g}$ of $V^{\natural}$ such that 
$\tilde{g}(e^i)=e^{g(i)}$. 
\end{lmm}

\pr 
By Lemma 3.2, we may assume that $g$ is an automorphism 
of $M_{D^{\natural}}$. 
Let $g((V^{\natural})^{\chi})$ be an $M_{D^{\natural}}$-module defined by 
$v_m(g\cdot u))=g\cdot(g^{-1}(v)_mu)$ for $v\in M_{D^{\natural}}$ and 
$u\in (V^{\natural})^{\chi}$.  Clearly, 
$\tilde{\tau}(g((V^{\natural})^{\chi})=
g^{-1}(\chi)$ and 
$$g(V^{\natural})=\oplus_{\chi\in S^{\natural}}g((V^{\natural})^{\chi}) $$
has a VOA structure containing $g(M_{D^{\natural}})=M_{D^{\natural}}$ 
by Theorem 3.8. 
We will prove that there is an $M_{D^{\natural}}$-isomorphism 
$$\pi_{\chi}: g((V^{\natural})^{\chi}) \to (V^{\natural})^{g(\chi)} $$ 
for $\chi\in S^{\natural}$.  
Then, by the uniqueness theorem (Theorem 3.3), 
there are scalars $\la_{\chi}$ such that 
$$\phi:g(V^{\natural}) \to V^{\natural} $$ 
given by $\phi=\la_{\chi}\pi{\chi}$ on $g((V^{\natural})^{\chi})$ 
is a VOA-isomorphism. Hence, 
$\tilde{g}(v)=\phi (g\cdot v)$ for $v\in V^{\natural}$ is one of 
the desired automorphisms of $V^{\natural}$. 

Since 
$S^{\natural}=\{(\al,\be,\ga):\al,\be,\ga\in S_{E_8}, \be,\ga=\al 
\mbox{ or }\al^c \}$, 
${\Aut}(S^{\natural})=S_3\times {\Aut}(S_{E_8})$, where $S_3$ is the symmetric 
group on three letters.  As we showed in \S 5,  
${\Aut}(S_{E_8})\cong GL(5,2)_1=\{ g\in GL(5,2): g{}^t(10000)={}^t(10000) \}$. 
In particular, $g$ leaves 
$D^3=D_{E_8}\oplus D_{E_8}\oplus D_{E_8}$ and $D^{\natural}$ invariant. 
Set $\chi=(\al,\be,\ga)$. 
We first assume that $g\in S_3$. 
Since $(V^{\natural})^{\chi}={\Ind}_{D^3}^{D^{\natural}}(W^{(\al,\be,\ga)})$ 
and $W^{(\al,\be,\ga)}$ is given by (7.10), we have
$g(W^{(\al,\be,\ga)})\cong W^{g(\al,\be,\ga)}$ 
as $M_{D^3}$-modules and so we have the desired isomorphism for $g\in S_3$.  
Now, assume $g=(h,h,h)$ with $h\in {\Aut}(S_{E_8})$. Set $j=h(1)$. 
By Lemma 5.1, $h(V_{E_8}^{\al})\cong V_{E_8}^{h(\al)}$ and so 
$g(W^{(\al,\al,\al)})\cong W^{(h(\al),h(\al),h(\al))}$. 
Hence, we may assume $\chi=(\al,\al,\al^c)$.  By the definition,  
$$ \begin{array}{rl}
g(W^{(\al,\al,\al^c)})
=&h(RV_{E_8}^{\al})\otimes h(RV_{E_8}^{\al}) 
\otimes h(V_{E_8}^{\al^c})\cr
\cong &(h(R))V_{E_8}^{h(\al)}\otimes (h(R))V_{E_8}^{h(\al)}
\otimes V_{E_8}^{h(\al^c)} \cr
\end{array} $$
as $M_{D_{E_8}}\otimes M_{D_{E_8}}\otimes M_{D_{E_8}}$-modules. 
Since $R=M_{D_{E_8}+\xi_1}$, $h(R)=M_{D_{E_8}+\xi_j}$, where 
$\xi_j=(0^{j-1}10^{16-j})$.   
Since $(\xi_1+\xi_j,\xi_1+\xi_j,0^{16})\in D^{\natural}$, 
$(R\times h(R))\otimes (R\times h(R))\otimes M_{D_{E_8}}$ 
is a submodule $M_{D^3+(\xi_1+\xi_j,\xi_1+\xi_j,0^{16})}$ 
of $M_{D^{\natural}}$ and so  
we have 
$$ \begin{array}{rl}
g(V^{\natural})^{\chi}=
&g({\Ind}_{D_{E_8}^3}^{D^{\natural}}W^{(\al,\al,\al^c)}) \cr
=&{\Ind}_{D_{E_8}^3}^{D^{\natural}}
(h(R))V_{E_8}^{h})\otimes (h(R))V_{E_8}^{h(\al)}) 
\otimes (V_{E_8}^{h(\al^c)}) \cr
\cong &{\Ind}_{D_{E_8}^3}^{D^{\natural}}
RV_{E_8}^{h(\al)}\otimes RV_{E_8}^{h(\al)}
\otimes V_{E_8}^{h(\al)^c} \cr
\cong &(V^{\natural})^{g(\chi)}. 
\end{array} $$  
\prend

Let $\Lambda$ be the Leech lattice and let $V_{\Lambda}$ be a lattice 
VOA constructed from $\Lambda$. 
The following result easily comes from the 
construction of $V_{\Lambda}$ in \cite{FLM2}. 

\begin{lmm}  ${\Aut}(V_{\Lambda})\cong ((\R^{\times})^{\oplus 24})Co.0$, 
where ${\R}^{\times}={\R}-\{0\}$ is the multiplicative group of ${\R}$. 
\end{lmm}

\pr 
Since $(V_{\Lambda})_1$ is a commutative Lie algebra ${\R}\Lambda$ 
of rank 24 and 
${\exp}(\al(0))=\sum_{i=0}^{\infty} {1\over i!}(\al(0))^i$ is an automorphism 
acting $\iota(x)$ as ${\exp}(\langle \al,x\rangle)\iota(x)$ 
for $\al\in (V_{\Lambda})_1$ and $x\in \Lambda$, we have 
an automorphism group ${\R^{\times}}^{\oplus 24}$, which is a normal 
subgroup of ${\Aut}(V_{\Lambda})$.  On the other hand, Frenkel, Lepowsky 
and Meurman \cite{FLM2} induced $g\in {\Aut}(\Lambda)$ into an automorphism 
of the group extension $\hat{\Lambda}=\{\pm \iota(x):x\in \Lambda\}$ and 
also into an automorphism of $V_{\Lambda}$ using cocycles.  Hence, 
$V_{\Lambda}$ has an automorphism group $({\R^{\times}}^{\oplus 24})Co.0$.  
We note that this is not split extension. Conversely, 
choose $g\in {\Aut}(V_{\Lambda})-({\R^{\times}}^{\oplus 24})Co.0$, 
then $g$ leaves $(V_{\Lambda})_1$ invariant 
and so it leaves a sub VOA $<(V_{\Lambda})_1>$ of free bosons and so 
$g$ acts on the lattice of highest weights of $<(V_{\Lambda})_1>$ in 
$V_{\Lambda}$, which is isomorphic to the Leech lattice.  
Multiplying an element of $Co.0$, 
we may assume that $g$ fixes all highest weights vectors 
$\iota(x):x\in \Lambda$ up to scalar multiple and so $g$ commutes with 
$x(0)$ for $x\in L$.  Hence, $g$ fixes 
all elements of $(V_{\Lambda})_1$ and so  
$g\in ({\R^{\times}}^{\oplus 24})$. 
\prend

\begin{thm}  ${\Aut}(V^{\natural})$ is the Monster simple group. 
\end{thm}

\pr  
As we proved, the full automorphism group of $V^{\natural}$ is 
finite.  
Set $\delta=\tau_{e^1}\tau_{e^2}$ and decompose 
$V^{\natural}$ into the direct sum 
$$ V^{\natural}=V^+\oplus V^-$$
of the eigenspaces of $\delta$,  
where $V^{\pm}=\{v\in V^{\natural}:\delta(v)=\pm v\}$. 
By the definition of $\tau_{e_i}$, 
$$V^+=\sum_{\al\in S^{\natural},\ 
\langle \al,(110^{46})\rangle=0} (V^{\natural})^{\al}. $$ 
Set $S_{\Lambda}=<(110^{46})>^{\perp}\cap S^{\natural}$ and 
$D_{\Lambda}=S_{\Lambda}^{\perp}$.  
Since $$S^{\natural}=\{(\al,\be,\ga): \al,\be,\ga\in S_{E_8}, 
\be,\ga\in \{\al,\al^c\} \}$$ and 
$$S_{E_8}=<(1^{16}),(1^80^8),(1^40^4)^2,(1^20^2)^4,(10)^8>, $$
we have the expression: 
$$S_{\lambda}=\{(a_1,...,a_{24})\in S^{\natural}: a_i=(00),(11) \}.$$
In particular, 
$\delta$ is equal to $\tau_{e^{2m-1}}\tau_{e^{2m}}$ for $m=1,...,24$. 
It is straightforward to check that 
$V^+$ has the structure given in Hypotheses I for 
$S=S_{\Lambda}$ and $D=D^{\natural}$. 
Since $S_{\Lambda}^{\perp}$ is larger than $D$, 
we can construct an induced VOA 
$$\tilde{V}_{\Lambda}=
{\Ind}_{D^{\natural}}^{D_{\Lambda}}(V^+).$$
Since $(S_{\Lambda})^{\perp}=D_{\Lambda}$,  
$\tilde{V}_{\Lambda}$ is a holomorphic VOA of rank $24$ by Theorem 6.1.  
It follows from the direct 
calculation that the codewords of $D_{\Lambda}$ of weight $2$ are 
$$\{ (110^{46}), (00110^{44}),...,(0^{46}11)\}. $$ 
We assert that $({\Ind}_{D^{\natural}}^{D_{\Lambda}}(V^{\natural})^{\al})_1=0$ 
for $\al\not=0$. Suppose 
$({\Ind}_{D^{\natural}}^{D_{\Lambda}}(V^{\natural})^{\al})_1\not=0$ for 
some $\al$. Then  
$|\al|=16$ and so $\al$ is one of 
$(1^{16}0^{32}),(0^{16}1^{16}0^{16}),(0^{32}1^{16})$, 
say $\al=(1^{16}0^{32})$. Since $(V^{\natural})^{\al}$ is given by 
${\Ind}_{D_{E_8}^3}^{D^{\natural}}(V_{E_8}^{(1^{16})}
\otimes M_{D_{E_8}+\xi_1}\otimes M_{D_{E_8}+\xi_1})$ and 
$D_{\Lambda}$ does not contain any word of form $(\ast\xi_1\xi_1)$, 
$({\Ind}_{D^{\natural}}^{D_{\Lambda}}(V^{\natural})^{\al})_1=0$.  
Consequently, 
$$\CG=(\tilde{V}_{\Lambda})_1=(M_{D_{\Lambda}})_1
=\oplus_{\al\in D_{\Lambda}, |\al|=2}
(M_{\al})_1 $$
is a commutative Lie algebra of rank 24.  
and $<(\tilde{V}_{\Lambda})_1>$ is a VOA 
of free bosons of rank 24. We note that $\CG$ has a positive definite 
invariant bilinear form $\langle\cdot,\cdot \rangle$ given by 
$v_1u=\langle v,u\rangle{\1}$ since $\tilde{V}_{\Lambda}$ has a positive 
definite invariant bilinear form. 
Hence, ${\C}\tilde{V}_{\Lambda}$ is isomorphic to a lattice 
VOA  ${\C}V_{\Lambda}$ of the Leech lattice $\Lambda$ by \cite{Mo}. 
More precisely, we will show the following lemmas in order to 
continue the proof of the theorem. 

\begin{lmm} 
$\tilde{V}_{\Lambda}$ is isomorphic to the lattice VOA $\tilde{V}_{\Lambda}$ 
of Leech lattice given in Proposition 2.2.  In particular, 
we can choose a set of mutually orthogonal vectors $\{x^1,...,x^{24}\}$ 
in $\Lambda$ of squared length $4$ such that 
$$ e^{2j-i}={1\over 16}x^j(-1)^2e^0+(-1)^i(\iota(x^j)+\iota(-x^j))  $$
for $j=1,...,24$ and $j=0,1$. 
Moreover, $(b_1b_1b_2b_2\cdots b_{24}b_{24})\in S_{\Lambda}$ if 
and only if there 
is $(a_i)\in {\Z}_2^{24}$ such that 
$x={1\over 2}\sum a_ix^i+{1\over 4}\sum b_ix^i\in \Lambda$. 
\end{lmm}

\pr
Set 
$$W=\{v\in \tilde{V}_{\Lambda}: x(n)v=0 \mbox{ for all }
x\in \CG \mbox{ and } n>0\}.$$ 
Then the actions of $\{x(0):x\in \CG \}$ on ${\C}W$ 
is diagonalizable since $\CG$ is commutative. Let $L$ be the set  
of highest weights of $\CG$ in ${\C}W$. It is 
easy to see that $L$ is an even unimodular positive definite 
lattice without roots since 
$(\tilde{V}_{\Lambda})_1=0$.  Hence, 
$L$ is the Leech lattice and 
${\C}\tilde{V}_{\Lambda}\cong {\C}V_{\Lambda}$. 

On the other hand, 
from Theorem 4.1, $\tilde{V}_{\Lambda}$ has a 
positive definite invariant bilinear form and 
it also has a ${\Z}_2$-grading 
$$\tilde{V}_{\Lambda}=(V^{\natural})^{<\delta>}\oplus V_{\Lambda}^-$$ 
by the definition of induced VOAs, where $V_{\Lambda}^-=
M_{(110^{46})+D^{\natural}}\times (V^{\natural})^{<\delta>}$. 

Let $\theta$ be an automorphism defined by 
$1$ on $(V^{\natural})^{<\delta>}$ and $-1$ on $\tilde{V}_{\Lambda}^-$.  
Since $\theta$ is acting on $(\tilde{V}_{\Lambda})_1$ as $-1$ and so 
it is equal to the automorphism of ${\C}V_{\Lambda}$ 
induced from $-1$ on $\Lambda$.  
Set $V=(V^{\natural})^{<\delta>}\oplus \sqrt{-1}\tilde{V}_{\Lambda}^-$.  
It is also a subVOA of ${\C}\tilde{V}_{\Lambda}.$
Let $\iota(x)$ denote a highest weight vector of $\CG$  
in ${\C}\tilde{V}_{\Lambda}$ with a highest weight $x\in \Lambda$. 
Namely, $u(0)\iota(x)=\langle u,x\rangle\iota(x)$ for $u\in \CG$.  
We note that $\theta(\iota(x))=(-1)^k\iota(x)$ for $\langle x,x\rangle=2k$. 
As a $\CG$-module, the space $W$ of highest weight vectors 
is a direct sum of irreducible $\CG$-modules $W^i$
whose dimension are less than or equal to $2$. If $\dim W^i=1$, then 
${\C}W^i={\C}\iota(x)$ for some $x\in \Lambda$. If $\dim W^i=2$, then 
${\C}W^i={\C}\iota(x)+{\C}\iota(y)$.  Since $W^i$ is irreducible, 
$\iota(x)$ and $\iota(y)$ are in the same homogenous space 
${\C}(\tilde{V}_{\Lambda})_k$ for some $k$. 
Since ${\C}\CG=\C\Lambda$, 
we have ${\Z}x={\Z}y$ and so $y=-x$. 
So $W^i$ has a basis $\{a\iota(x)+b\iota(-x), c\iota(x)+d\iota(-x)\}$ 
for some $a,b,c,d\in \C$. 
We may assume that $a\in {\R}$.  
Since $\tilde{V}_{\Lambda}$ has a positive definite 
invariant bilinear form, we have assume that 
$\{{1\over \sqrt{2}}(a\iota(x)+b\iota(-x)), 
{1\over\sqrt{2}}(c\iota(x)+d\iota(-x))\}$ is an orthonormal 
basis. Therefore, 
$b=(-1)^ka^{-1}$, $d=(-1)^kc^{-1}$ and $ad+bc=(-1)^k(ac^{-1}+a^{-1}c)=0$.  
Hence, $a^2=-c^2>0$ and so we have $c=\sqrt{-1}a$ and $d=-\sqrt{-1}b$.  
Since ${\C}W^i={\C}\iota(x)+{\C}\iota(-x)$ and $W^i={\C}W^i\cap 
\tilde{V}_{\Lambda}$, $\theta$ keeps $W^i$ invariant. 
Hence, $\theta(a\iota(x)+(-1)^ka^{-1}\iota(-x))
=a^{-1}\iota(x)+(-1)^ka\iota(-x)\in W^i$ and so we have $a=\pm 1$.  
Hence, $\iota(x)+(-1)^k\iota(-x), 
\sqrt{-1}(\iota(x)-(-1)^k\iota(-x))\in W$ 
and $\sqrt{-1}x(0)\in \CG$ for $x\in \Lambda$. Consequently, 
$\tilde{V}_{\Lambda}$ coincides with the lattice VOA $V_{\Lambda}$ 
defined in Proposition 2.2. 
We recall the structure 
$V_{\Z x}\cong L(\hf,0)\otimes L(\hf,0)\oplus L(\hf,\hf)\otimes L(\hf,\hf)$  
and  $(L(\hf,\hf)\otimes L(\hf,\hf))_1={\R}\sqrt{-1}x(-1){\1}$ for 
a VOA $V_{\Z x}$ with $\langle x,x\rangle=4$.  
Since $(\tilde{V}_{\Lambda})_1=(M_{D_{\Lambda}})_1
=\oplus_{i=1}^{24} (M_{\xi_{2i-1}+\xi_{2i}})_1$, 
$e^{2j}-e^{2j-1}\in L=\{ v\in V_{\Lambda}| x(n)v=0 \mbox{ for all 
} x\in (V_{\Lambda})_1 \mbox{ and } n>0 \}$ and 
${\R}(e^{2j}-e^{2j-1})+\sqrt{-1}{\R}x^j(0)(e^{2j}-e^{2j-1})$ is irreducible 
$\CG$-submodule of $L$. Hence, by the above arguments we have 
$$ e^{2j-i}={1\over 16}x^j(-1)^2e^0
+(-1)^i{1\over 4}(\iota(x^j)+\iota(-x^j))  $$
for some $x^j\in \Lambda$.  Since 
$$0=(e^{2j-1}+e^{2j})_1(e^{2k}-e^{2k-1})={1\over 64}\langle x^j,x^k\rangle^2
(\iota(x^k)+\iota(-x^k))$$
 for $k\not=j$, we have $\langle x^j,x^k\rangle=0$.  
Namely, $\{x^1,...,x^{24}\}$ is a set of mutually orthogonal 
vectors of $\Lambda$ with squared length 4.  
If $y=\sum c_ix^i\in \Lambda$, then $c_i\in {1\over 4}{\Z}$. 
Assume that $y={1\over 4}\sum b_ix^i+{1\over 2}\sum a_ix^i$ is in 
$\Lambda$ and 
set $W=V_{<x^1,...,x^{24}>+y}$ and $T^j=<e^{2j-1},e^{2j}>$.  
As we showed in \S 2, \\
(1) $b_j=1$ if and only if irreducible $T^j$-submodule of $W$ 
is isomorphic to $L(\hf,\st)\otimes L(\hf,\st)$. \\
(2) $b_j=0$ and $a_j=1$ if and only if irreducible $T^j$-submodule of $W$ 
is isomorphic to 
$L(\hf,\hf)\otimes L(\hf,0)$ or $L(\hf,0)\otimes L(\hf,\hf)$. \\
(3) $b_i=0$ and $a_j=0$ if and only if irreducible $T^j$-submodule of $W$ 
is isomorphic to $L(\hf,0)\otimes L(\hf,0)$ or 
$L(\hf,\hf)\otimes L(\hf,\hf)$. \\
In particular, we have $(b_1b_1b_2b_2\cdots b_{24}b_{24})\in S_{\Lambda}$.  

Conversely, if $\ga=(b_1b_1b_2b_2\cdots b_{24}b_{24})\in S_{\Lambda}$, 
then $<(\tilde{V}_{\Lambda})_1>$ acts
on $(\tilde{V}_{\Lambda})^{\ga}$ and so 
$(\tilde{V}_{\Lambda})^{\ga}\cap L\not=0$. 
Hence, by the above arguments, there is an element $x\in 
\Lambda$ such that $\iota(x)\in \tilde{V}_{\Lambda}$ or 
$\iota(x)+(-1)^{|x|/2}\iota(-x)\in (\tilde{V}_{\Lambda})^{\ga}$.  
We can also find $(a_i)\in {\Z}_2^{24}$ 
such that $x={1\over 2}\sum a_ix^i+{1\over 4}\sum b_ix^i$. 
\prend

\begin{lmm}  
For any $y\in \Lambda$ with squared length 4, 
$\tau_{e(y)^+}=\tau_{e(y)^-}$ in ${\Aut}(V_{\Lambda})$ and 
$\tau_{e(y)^+}\in {<\pm 1>}^{\oplus 24}\subseteq ({\R}^{\times})^{\oplus 24}$.
\end{lmm}

\pr
Since $Co.0$ acts on the set of all vectors in 
$\Lambda$ with squared length 4 transitively, we may assume that $y=x^1$ and 
$e(y)^+=e^1$ and $e(y)^-=e^2$, where $\{x^1,...,x^{24}\}$   
is the set defined in the above lemma.  
By the arguments in the proof of the above lemma, 
it is clear that $\tau_{e(y)^+}=\tau_{e(y)^-}$. 
Since 
$\tau_{e^1}\iota(x)=(-1)^{\langle x^1,x\rangle}\iota(x)$ and 
$[\tau_{e^1},x(-n)]=0$, 
we have $\tau_{e^1}\in {<\pm 1>}^{\oplus 24}$. 
\prend

Let's go back to the proof of the theorem 9.3. 
Set $V_{\Lambda}=(V^{\natural})^{\delta}\oplus \sqrt{-1}V^-$.  
By the proof of Proposition 2.2 and the above lemma, 
$V_{\Lambda}$ is isomorphic to a lattice VOA of the Leech lattice which 
is given by the ordinary construction.  
Let $\theta$ be an automorphism of $V_{\Lambda}$ defined by 
1 on $(V^{\natural})^{\delta}$ and $-1$ on $\sqrt{-1}V^-$. 
We identify $(V^{\natural})^{<\delta>}$ and 
$V_{\Lambda}^{\theta}$. 
Let $J$ be the set of all rational conformal vectors in 
$(V^{\natural})^{<\delta>}$ with central charge $\hf$.  Set \\ 
$K^{\natural}=<\tau_e: e\in J>\subseteq {\Aut}(V^{\natural})$,   \\
$K=<\tau_e: e\in J>\subseteq {\Aut}((V^{\natural})^{<\delta>})$ and  \\
$K_{\Lambda}=<\tau_e: e\in J>\subseteq {\Aut}(V_{\Lambda})$. \\
Set $G={\Aut}(V^{\natural})$ and $H={\Aut}(V_{\Lambda})$.  By Lemma 9.2, 
$H\cong ({\R^{\times}}^{\oplus 24})Co.0$ and  \\
$C_H(<\theta>)\cong 2^{24}Co.0$.  Clearly, 
$K^{\natural}\subseteq C_G(<\delta>)$ and $K_{\Lambda}\subseteq 
C_H(<\theta>)$. 

By the restrictions from $V^{\natural}$ to $(V^{\natural})^{<\delta>}$ 
and from $V_{\Lambda}$ to $V_{\Lambda}^{<\theta>}$, we have epimorphisms 
$\pi^{\natural}: K^{\natural} \to K$ and 
$\pi_{\Lambda}:K_{\Lambda}\to K$.  By \cite{DM2}, 
$Ker(\pi^{\natural})=<\delta>$ and 
$Ker(\pi_{\Lambda})=<\theta>\cap K_{\Lambda}$. 
So we have the following diagram:

\setlength{\unitlength}{1mm}
\begin{picture}(120,55)
\put(10,50){$G={\Aut}(V^{\natural})$}
\put(14,44){\line(0,1){4}}
\put(12,40){$C_G(\delta)$}
\put(14,24){\line(0,1){14}}
\put(12,20){$K^{\natural}$}
\put(24,42){\line(4,-1){20}}
\put(20,22){\line(4,-1){40}}
\put(14,14){\line(0,1){4}}
\put(10,10){$<\delta>$}
\put(60,50){${\Aut}((V^{\natural})^{\delta})$}
\put(62,48){\line(-1,-1){10}}
\put(68,48){\line(1,-1){10}}
\put(68,12){\line(1,1){10}}
\put(62,12){\line(-1,1){10}}
\put(65,14){\line(0,1){34}}
\put(45,30){$\ol{C_G(\delta)}$}
\put(74,30){$\ol{C_H(\theta)}$}
\put(63,10){$K$}
\put(65,5){\line(0,1){4}}
\put(64,0){$1$}
\put(110,50){$H={\Aut}(V_{\Lambda})$}
\put(114,44){\line(0,1){4}}
\put(112,40){$C_H(\theta)$}
\put(114,24){\line(0,1){14}}
\put(112,20){$K_{\Lambda}$}
\put(110,42){\line(-4,-1){20}}
\put(110,22){\line(-4,-1){40}}
\put(114,14){\line(0,1){4}}
\put(110,10){$<\theta>\cap K_{\Lambda}$}
\end{picture}
\newline

First, we will show that 
$K_{\Lambda}\not\subseteq 2^{24}\!<\!\theta\!>$.  
Let g be a permutation on 48 letters $\{1,...,48\}$ 
such that $g$ fixes all $1+2m$ and $3+2m$ and switches $2+4m$ and $4+4m$ 
for $m=0,...,11$. 
It is straightforward to check that $g$ is an automorphism of $S^{\natural}$.  
By Lemma 9.1, there is an automorphism 
$\tilde{g}\in {\Aut}(V^{\natural})$ such that 
$\tilde{g}(e^i)=e^{g(i)}$.  
Set $\delta'=\tau_{e^1}\tau_{e^4}\ (=\tilde{g}(\delta))$ 
and $\tilde{L}'_{\Lambda}=g(\tilde{L}_{\Lambda})$ and then 
apply the above arguments. 
By the above lemma, there is a set of mutually orthogonal vectors 
$\{\tilde{x}^1,...,\tilde{x}^{24}\}$ in 
$\Lambda$ of squared length 4 such that 
$$ e^{2j-i}={1\over 16}\tilde{x}^j(-1)^2\iota(0)
+(-1)^i{1\over 4}(\iota(\tilde{x}^j)
+\iota(-\tilde{x}^j)).  $$ 
It is easy to see that 
$\ga=(0^81^80^81^80^81^8)\in S_{\Lambda}$. 
Since $((V^{\natural})^{\ga})_2\not=0$, there is $y\in \Lambda$ of squared 
length 4 such that $\langle y,x^i\rangle\equiv 1 \pmod{2}$ 
if and only if $i\in Supp(\ga)$. Then  
$e(y)={1\over 16}y(-1)^2\iota(0)+{1\over 4}(\iota(y)+\iota(-y))$ 
is a rational conformal vector in 
$(V_{\Lambda})^{<\theta, \tau_{e^1},\tau_{e^2},...,\tau_{e^8}>}$. 
In particular,  $g(e(y))\in (V^{\natural})^{<\delta>}$. 
Since $\langle y,x^5\rangle\equiv 1 \pmod{2}$, we have 
$\tau_{e(y)}(\iota(\pm x^5)=-\iota(\pm x^5)$ and so 
$\tau_{e(y)}$ switches $e^{9}$ and $e^{10}$. 
On the other hand, $\tilde{g}$ fixes $e^9$ and switches 
$e^{10}$ and $e^{12}$. 
Hence, $\tau_{\tilde{g}}(e(y))$ switches $e^9$ and $e^{12}$ and so 
$\tau_{\tilde{g}(e(y))}$ does not belong to $2^{24}<\theta>$. 

Since $K_{\Lambda}$ is generated by all conformal vectors in 
$(V_{\Lambda})^{<\theta>}$, $K_{\Lambda}$ is a normal subgroup of 
$C_H(<\theta>)\cong 2^{24}Co.0$ and so we have $K_{\Lambda}=C_H(<\theta>)$.  
Hence, $K\cong 2^{24}Co.1$ and so we have $K^{\natural}=2^{1+24}Co.1$.  
If $O_2(K^{\natural})$ is an Abelian 2-group, then 
$2^{1+24}=<\delta>\oplus {\Z}_2^{24}$ 
as a $Co.1$-module.  
Let $y$ be a vector of $\Lambda$ of squared length 4 and $\langle y, 
x^{24}\rangle=1$. Then $e^{\pm}(y)\in (V^{\natural})^{<\delta>}$ 
and $\tau_{e^+(y)}$ fixes 
$\delta=\tau_{e^1}\tau_{e^2}=\tau_{e^{47}}\tau_{e^{48}}$ and 
switches $e^{47}$ and $e^{48}$.

By Lemma 9.4, $\tau_{e^{47}},\tau_{e^{48}}\in 2^{1+24}$. 
Since $\delta=\tau_{e^{47}}\tau_{e^{48}}$, 
we may assume $e^{47}\in {\Z}_2^{24}$ and $e^{48}\not\in {\Z}_2^{24}$, 
which contradicts that $\tau_{e(y)}$ switches $e^{47}$ and $e^{48}$.  
Hence, $2^{1+24}$ is a non-abelian and so 
$2^{1+24}$ are isomorphic to a central 
extension of $\Lambda/2\Lambda$ using the inner product, 
since $Co.1$ acts on faithfully. 
By Lemma 9.1, ${\Aut}(V^{\natural})$ 
contains a subgroup whose restriction on $\{e^1,...,e^{48}\}$ 
is isomorphic to $GL(5,2)_1\times S_3$, 
where $S_3$ permutes 3 components of 
$V_{E_8}^{\ots 3}$ and $GL(5,2)_1=\{A\in GL(5,2):Av=v \mbox{ for }
v={}^t(10000)\}$.   Set 
$\delta_1=\tau_{e^1}\tau_{e^3}$ and $B^2=<\delta,\delta_1>$.  
Denote $\delta$ and $\delta\delta_1$ by $\delta_0$ and $\delta_2$, 
respectively. Since a subgroup of $GL(5,2)_1$ acts on 
$\{\delta_0,\delta_1,\delta_2\}$ transitively and 
$e^3$ is given by a vector of $\Lambda$ of squared length 4, we have 
$N_{{\Aut}(V^{\natural})}(B^2)\cong 2^{2+12+22}(S_3\times M_{24})$ 
from the structure of $C_{{\Aut}(V^{\natural})}(\delta)\cong 2^{1+24}Co.1$.  
Similarly, all nontrivial elements of  
$B^3=<\tau_{e^1}\tau_{e^2}, \tau_{e^1}\tau_{e^3}, 
\tau_{e^1}\tau_{e^5}>$ are conjugate 
by the actions of $GL(5,2)_1\subseteq {\Aut}(V^{\natural})$ and so 
$N_{{\Aut}(V^{\natural})}(B^3)\cong 2^{3+6+12+18}(3 S_6\times L_3(2))$. 
By the same arguments, we can calculate the normalizer of 
$B^4=<\tau_{e^1}\tau_{e^2}, 
\tau_{e^1}\tau_{e^3}, 
\tau_{e^1}\tau_{e^5},
\tau_{e^1}\tau_{e^9}>$.
We leave these calculation to the reader. 

We will next prove that ${\Aut}(V^{\natural})$ is a simple group. 
Let
$H$ be a nontrivial minimal normal subgroup of ${\Aut}(V^{\natural})$. Then 
$C_H(\delta_i)$ is a normal subgroup of $C(\delta_i)=2^{1+24}Co.1$ 
for $i=0,1,2$. 
Hence, 
$C_H(\delta_i)=2^{1+24}Co.1$ or $C_H(\delta_i)=2^{1+24}$ or 
$C_H(\delta_i)=<\delta_i>$. 
We note that $\delta_i$ $(i=0,1,2)$ are conjugate to each other in 
${\Aut}(V^{\natural})$ and so $C_H(\delta_i)\cong C_H(\delta_0)$ for $i=1,2$.  
In any cases, $\delta_i\in H$ and so 
$C_H(\delta_i)\not=<\delta_i>$ since $\delta_j\in <C_H(\delta_i): i=1,2,3>=H$.
If $P=C_H(\delta_1)=2^{1+24}$ then $P$ is a Sylow 2-subgroup 
of $H$. Since $|P:C_P(\de_2)|=2$ and $C_P(\de_2)$ is not abelian, 
$[C_P(\de_2),C_P(\de_2)]=<\de_1>$, which contradicts 
$[C_H(\de_2),C_H(\de_2)]=<\de_2>$.  
Hence we have $C_H(\delta_i)=2^{1+24}Co.1$. Since $<\delta_i>$ 
is a characteristic subgroup of a Sylow 2-subgroup of $H$, we have 
$H={\Aut}(V^{\natural})$ and so  
${\Aut}(V^{\natural})$ is a simple group.   
By the characterization of the Monster simple group and 
the above facts, we know ${\Aut}(V^{\natural})$ is the Monster simple group, 
see \cite{I},\cite{S},\cite{T}. 
\prend

Since $V^{\natural}$ is a holomorphic VOA with rank 24 with 
$(V^{\natural})_1=0$ 
and the Monster simple group 
acts on $B=V^{\natural}_2$ faithfully, $B$ is isomorphic to 
the Griess algebra constructed in \cite{Gr}. We have also proved that 
$(V^{\natural})^{\delta}$ is isomorphic to $(V_{\Lambda})^{\theta}$.  
Hence, $V^{\natural}$ is equal to the moonshine VOA constructed in 
\cite{FLM2}.

\section{Meromorphic VOAs}
In this section, we will construct an infinite series of holomorphic VOAs 
whose full automorphism groups are finite.  
We will adopt the notation from \S 7 and repeat the similar constructions 
as in \S 7.

For $n=1,2,\cdots $, set 
$$ S^{\natural}(n)=<(\{0^{16}\}^i1^{16}\{0^{16}\}^{2n-i}), (\{\al\}^{2n+1}):
\al\in S(P), i=1,\cdots,2n>.  $$
$S^{\natural}(n)$ is an even linear code of length $16+32n$ and 
$(S^{\natural}(n))^{\perp}$ contains a direct sum $D^{2n+1}$ of $2n+1$ 
copies of $D$ for each $n$.  Let $\ga$ be an element of 
$S^{\natural}(n)$, 
then there is $\al\in S(P)$ such that 
$$ \ga=(\be_1,...,\be_{2n+1}), $$
where $\be_i\in \{ \al,\al^c \}$.  We may assume that the number of 
$\be_i$ satisfying $\be_i=\al$ is odd.  Set 
$$W^{\gamma}=\ots_{i=1}^{2n+1} \tilde{W}^{\be_i},  $$
where 
$$ \begin{array}{ll}
\tilde{W}^{\be_i}={V_{E_8}}^{\al} & 
\mbox{ if } \be_i=\al \quad \mbox{ and } \cr
\tilde{W}^{\be_i}=R{V_{E_8}}^{\al^c} 
& \mbox{ if } \be_i=\al^c.
\end{array} $$
Set 
$$V^3(n)=\boplus_{\gamma\in S^{\natural}(n)}W^{\gamma}   $$
and 
$$  V^{\natural}(n)={\Ind}_{M_{D^{2n+1}}}^{M_{(S^{\natural}(n))^{\perp}}}
(V^3(n)).  $$
Then we can show that 
$V^{\natural}(n)$ has a VOA structure by exactly the same proof as in 
the construction of $V^{\natural}$.  It also satisfies 
$(V^{\natural}(n))_1=0$.  
Moreover, it is a holomorphic VOA by Theorem 6.1 and its full automorphism 
group is finite by Theorem 9.2.

\section{Characters}
In this section, we will calculate the characters of $3C$ element 
and $2B$ element of the 
Monster simple group. 
It follows from our construction that we can induce an automorphism of 
$D^{\natural}$ into an automorphism of $V^{\natural}$. 

\subsection{$3C$}
Clearly, $\hat{g}=(1,17,33)(2,18,34)...(16,32,48)$ is an  
automorphism of $D^{\natural}$. 
Let $g$ be an automorphism of $V^{\natural}$ induced from $\hat{g}$.  
By the definition, $g$ acts on $\{e^i:i=1,...,48\}$ as 
$(1,17,33)(2,18,34)...(16,32,48)$. 

$V^{\natural}$ contains 
$M_{D^3}=M_D\otimes M_D\otimes M_D$, where $D=D_{E_8}$. 
We view $V^{\natural}$ as an $M_D\otimes M_D\otimes M_D$-module. 
Since $g$ permutes $\{ V^{\chi}:\chi\in S^{\natural}\} $, 
we obtain 
$$
\begin{array}{rl}
\ch_{V^{\natural}}(g,z)=& \tr_{g,z}(V^{\natural}) \cr
=&\tr_{g,z}(\boplus_{\chi^g=\chi\in S^{\natural}}V^{\chi}) \cr
=&\tr_{g,z}(\boplus_{\al\in D_{E_8}}V^{(\al,\al,\al)}),
\end{array}$$
where $\tr_{g,z}(V)=\sum \tr(g)_{|V_n}e^{2\pi inz}$ for $V=\oplus V_n$. 

By the definition of $V^{(\al,\al,\al)}$, 
$$V^{\al,\al,\al}=
{\Ind}_{M_{D^3}}^{M_{D^{\natural}}}(
V_{E_8}^{\al}\ots V_{E_8}^{\al}\ots 
V_{E_8}^{\al}). $$ 
It follows from the definition of induced modules, 
$${\Ind}_{M_{D^3}}^{M_{D^{\natural}}}(U)
\cong \boplus_{\mu\in D^{\natural}/D^3}M_{D^3+\mu}\times U$$
as $M_{D^3}$-modules. 
Since $D^{\natural}=\{(\al,\be,\ga):\al+\be+\ga\in D, \al,\be,\de 
\mbox{ even }\}$, we obtain that 
$g(D^3+\mu)=D^3+\mu$ if and only if $\mu\in D^3$.  
Hence, 
$$\begin{array}{rl}
&\tr_{g,z}(V^{(\al,\al,\al)})\cr
=&\tr_{g,z}(V_{E_8}^{\al}\otimes V_{E_8}^{\al}\otimes V_{E_8}^{\al}) \cr
=&\tr_{g,3z}(V_{E_8}^{\al}).
\end{array}
$$
Therefore, we have 
$$\begin{array}{rl}
\ch_{V^{\natural}}(g,z)=&\sum_{\al\in D_{E_8}}\tr_{g,3z}(V_{E_8}^{\al}) \cr
=&\tr_{g,3z}(V_{E_8})=\ch_{V_{E_8}}(1,3z). 
\end{array}$$

\subsection{$1$ and $2B$}
Let $\delta=\tau_{e^1}\tau_{e^2}$.  
We proved that $(V^{\natural})^{<\de>}$ is 
isomorphic to $(V_{\Lambda})^{<\theta>}$.  
Hence, 
$$ \ch((V^{\natural})^{<\de>})=1+98580q^2+... . $$
So we will calculate 
the character of $(V^{\natural})^{-}=\{v\in V^{\natural}:\de(v)=-v\}$. 
It follows from the definition of $\tau_{e^i}$ that  
$$ 
\ch((V^{\natural})^-)
=\sum_{\langle \chi,(110^{46})\rangle=1} \ch((V^{\natural})^{\chi}). $$
Set $\chi=(\al,\be,\ga)$ with $\al,\be,\ga\in {\Z}_2^{16}$. 
Assume $\langle \chi,(110^{46})\rangle=1$. Then 
the weight of $\al$ is $8$ and so the weight of $\chi$ is 24.  
Hence, $\dim D^{\natural}_{\chi}=7+7+4$ and so the multiplicity 
of every irreducible $T$-submodule of $(V^{\natural})^{\chi}$ is $2^6$. 
Let $U$ be an irreducible $T$-submodule of $(V^{\natural})^{\chi}$.
It follows from the total degree that 
the number of $L(\hf,\hf)$ in $U=\otimes_{i=1}^{48} L(\hf,h^i)$ is odd. 
On the other hand, let $\ga$ be an odd word with 
$Supp(\ga)\cap Supp(\chi)=\emptyset$. 
By the action of $M_{D^{\natural}}$, 
there exists an irreducible $T$-submodule isomorphic to $\otimes L(\hf,h^i)$ 
with $h^i=\hf$ for $i\in Supp(\ga)$, $h^i=\st$ for $i\in Supp(\chi)$ 
and $h^i=0$ for $i\not\in Supp(\chi+\ga)$. 
Hence 
$$\begin{array}{rl}
 \ch((V^{\natural})^{\chi})
=&2^6 \ch\{ L(\hf,\st)^{\otimes 24}
\hf((L(\hf,0)+L(\hf,\hf))^{\otimes 24}-(L(\hf,0)-L(\hf,\hf))^{\otimes 24})\}\cr
=&32 q^{3/2}\prod_{n\in \N}(1+q^n)^{24}(\prod_{n\in \N+\hf}(1+q^n)^{24}
-\prod_{n\in \N+\hf}(1-q^n)^{24}). 
\end{array} $$
Since there are 64 codewords $\chi$ such that 
$\langle \chi,(110^{46})\rangle=1$, we have
$$ \begin{array}{rl}
ch((V^{\natural})^-)
=&2^{11} q^{3/2}\prod_{n\in \N}(1+q^n)^{24}(\prod_{n\in \N+\hf}(1+q^n)^{24}
-\prod_{n\in \N+\hf}(1-q^n)^{24}) \cr
=&2^{11}q^{3/2}(1+24q+...)(48q^{1/2}+...) \cr
=&2^{12}(24q^2+...). 
\end{array} $$

In particular, we obtain $(V^{\natural})_1=0$ and 
$\dim (V^{\natural})_2=196884$.

\end{document}